\documentclass[twoside, 10pt]{article}

\usepackage[hmargin= 0.8in, twosideshift=0in,height=9.3in]{geometry}
\usepackage{amssymb,amsthm, amsfonts,amsmath,mathrsfs, amsxtra, url, hyperref}
\usepackage{delarray,verbatim, natbib}
\renewcommand{\cite}{\citet}

\usepackage{fancyhdr}
\fancypagestyle{plain}{
 \fancyhf{}

}

 \numberwithin{equation}{section}

\newtheorem {thm}{Theorem}[section]
\newtheorem {prop}{Proposition}[section]
\newtheorem {lemm}{Lemma}[section]
\newtheorem {deff}{Definition}[section]

\newtheorem {rem}{Remark}[section]

\newtheorem {eg}{Example}[section]

\def\ba{\begin{array}}
\def\ea{\end{array}}
\def\bea{\begin{eqnarray}}
\def\eea{\end{eqnarray}}
\def\beas{\begin{eqnarray*}}
\def\eeas{\end{eqnarray*}}
\def\bi{\begin{itemize}}
\def\ei{\end{itemize}}

\def\a{\alpha}
\def\g{\gamma}
\def\d{\delta}
\def\e{\varepsilon}
\def\z{\zeta}
\def\k{\kappa}
\def\l{\lambda}

\def\si{\sigma}

\def\t{\tau}
\def\th{\theta}
\def\o{\omega}

\def\vf{\varphi}
\def\p{\psi}

\def\G{\Gamma}
\def\L{\Lambda}
\def\O{\Omega}

\def\P{\Psi}

\def\U{\Upsilon}

\def\bF{{\bf F}}

\def\cA{{\cal A}}

\def\cD{{\cal D}}

\def\cF{{\cal F}}
\def\cG{{\cal G}}

\def\cI{{\cal I}}

\def\cM{{\cal M}}
\def\cN{{\cal N}}

\def\cP{{\cal P}}
\def\cQ{{\cal Q}}

\def\cS{{\cal S}}

\def\cZ{{\cal Z}}

\def\hC{\mathbb{C}}

\def\hH{\mathbb{H}}

\def\hK{\mathbb{K}}
\def\hL{\mathbb{L}}

\def\hN{\mathbb{N}}

\def\hR{\mathbb{R}}

\def\hZ{\mathbb{Z}}

\def\sB{\mathscr{B}}

\def\sE{\mathscr{E}}

\def\sP{\mathscr{P}}

\def\({\textnormal{(}}
\def\){\textnormal{)}}
\def\[{[\neg[}
\def\]{]\neg]}
\def\lan{\langle}
\def\ran{\rangle}

\def\no{\noindent}

\def\ss{\smallskip}
\def\ms{\medskip}
\def\bs{\bigskip}
\def\q{\quad}
\def\qq{\qquad}

\def\neg{\negthinspace}
\def\dneg{\neg \neg}
\def\tneg{\neg \neg \neg}

\def\ol{\overline}
\def\ul{\underline}
\def\ua{\mathop{\uparrow}}
\def\da{\mathop{\downarrow}}

\def\wt{\widetilde}

\def\dtp{{\hbox{$dt \otimes dP-$a.e.}}}
\def\pas{{\hbox{$P-$a.s.}}}

\def\hb{\hbox}
\def\dis{\displaystyle}
\def\cd{\cdot}
\def\cds{\cdots}

\def\fa{\,\forall \,}
\def\pa{\partial}

\def\dfnn{\stackrel{\triangle}{=}}
\def\b1{{\bf 1}}
\def\qed{\hfill $\Box$ \medskip}

\def\essinf{\mathop{\rm ess\,inf}}
\def\esssup{\mathop{\rm ess\,sup}}
\def\liminf{\mathop{\ul{\rm lim}}}
\def\limsup{\mathop{\ol{\rm lim}}}
\newcommand{\esup}[1]{ \underset{#1}{\esssup}\,}
\newcommand{\einf}[1]{ \underset{#1}{\essinf}\,}

\newcommand{\linf}[1]{ \underset{#1}{\liminf}}
\newcommand{\lmt}[1]{ \underset{#1}{\lim}}
\newcommand{\lmtu}[1]{ \underset{#1}{\lim} \neg \ua \,}
\newcommand{\lmtd}[1]{ \underset{#1}{\lim} \neg \da \,}

\begin{document}

\title{\bf  Optimal Stopping for Dynamic Convex Risk Measures}
\author{
Erhan Bayraktar\thanks{ \noindent Department of
  Mathematics, University of Michigan, Ann Arbor, MI 48109; email:
{\tt erhan@umich.edu}. This author is supported in part by the National
Science Foundation, under grant number DMS-0906257.}$~$,
$~~$Ioannis Karatzas\thanks{INTECH Investment Management, One Palmer Square, Suite 441, Princeton, NJ 08542; e-mail: {\tt ik@enhanced.com}; and 
\newline Department of Mathematics, Columbia University, New York, NY 10027;
e-mail: {\tt ik@math.columbia.edu.} Research supported in part by the National Science Foundation under Grant NSF-DMS-09-05754.}$\,\,$,
$~~$Song Yao\thanks{
\noindent Department of
  Mathematics, University of Michigan, Ann Arbor, MI 48109; email: {\tt songyao@umich.edu.} } }
\date{ }

\maketitle

 \begin{abstract}

We use martingale and stochastic analysis techniques to study a continuous-time optimal stopping problem, in which the decision maker uses a dynamic convex risk measure to evaluate future rewards.
We also find a saddle point for an equivalent zero-sum game of control and stopping,     between an agent (the ``stopper") who chooses the termination time of the game, and an agent (the ``controller", or ``nature") who   selects the probability measure.
\end{abstract}


 \smallskip
\section{Introduction}

Let us consider a complete, filtered probability space $\, (\Omega, \mathcal{F}, P)$, $\bF=\{\mathcal{F}_t\}_{t \geq 0}\,$, and on it    a bounded,  adapted process $Y$ that  satisfies certain regularity conditions.  Given any stopping time $\nu$ of the filtration $\bF$,  our goal is to find a stopping time $\t_*(\nu) \in \cS_{\nu,T}$ that satisfies
  \bea \label{eq:defn-op-rm-intro}
  \underset{ \g  \in  \cS_{\nu,T}}{\essinf}\rho_{\nu, \g}\left(Y_\g\right)
  =\rho_{\nu, \t_*(\nu)}\left( Y_{\t_*(\nu)}\right),\qquad P-\hbox{a.s.}
  \eea
Here $\cS_{\nu,T}$ is the set of stopping times $\g$
satisfying $\nu \leq \g \leq T$, $P-$a.s., and the collection of
functionals $\big\{\rho_{\nu, \g}: \mathbb{L}^\infty(\cF_\g) \rightarrow
\mathbb{L}^\infty(\cF_\nu)\big\}_{  \nu  \in \cS_{0,T},\, \gamma  \in \cS_{\nu,T} \,   }\,$ is a ``dynamic convex risk
measure" in the sense of \cite{DPR_2009}. Our motivation is to solve the optimal stopping problem of a decision maker who evaluates future rewards/risks using dynamic convex risk measures rather than  statistical expectations. This question can also be cast  as a {\it robust optimal stopping} problem, in which the decision maker has to act in the presence of so-called ``Knightian uncertainty" regarding the underlying probability measure.

   When the filtration $\,\bF\,$ is generated by a  Brownian motion, the dynamic convex risk measure
admits the following representation: There exists a suitable nonnegative function  $f$, convex in its spatial argument, such that the representation
$$
 \rho_{\nu, \g}(\xi)\,=\,\esup{Q \in \cQ_\nu}\, E_Q \left[-\xi- \int_\nu^\g
  f \neg \left(s, \th^Q_s\right) \neg ds\,\Big|\, \cF_\nu
 \right], \q \pas
 $$
holds for all $\xi \in \mathbb{L}^{\infty}(\mathcal{F}_{\gamma})$. Here $\cQ_\nu$ is the collection
of  probability measures $Q$ which are equivalent to  $P$ on $\mathcal{F}$,   equal to $P$ on $\mathcal{F}_{\nu}$, and
satisfy  a certain integrability condition; whereas
$\theta^{Q}$ is the predictable process whose stochastic exponential gives the density of $Q$ with respect to $P$.
In this setting we  establish a minimax result, namely
\bea \label{eq:intro-value}
V(\nu) \triangleq  \underset{\g \in \cS_{\nu,T}}{\esssup}
 \left(  \underset{Q \in \cQ_\nu}{\essinf}\,E_Q \neg \left[Y_\g \neg
 +\neg \int_\nu^\g \neg f \neg \left(s, \th^Q_s\right) \neg ds \,\Big|\,\cF_\nu\right] \right)
 = \underset{Q \in \cQ_\nu}{\essinf} \neg \left(\, \underset{\g \in \cS_{\nu,T}}{\esssup}\,
  E_Q \neg \left[Y_\g \neg +\neg \int_\nu^\g \neg f \neg \left(s, \th^Q_s\right) \neg ds \,\Big|\,\cF_\nu\right]
  \right),
 \eea
and construct an optimal stopping time  $\tau(\nu)$
as the limit of stopping times that are optimal under expectation criteria ---
see Theorem~\ref{V_process}.
We show that  the  process
   $ \left\{ \b1_{\{t \ge \nu\}}  V\big(\t(\nu) \land t\big) \right\}_{t \in [0,T]}$ admits an RCLL
    modification $V^{0,\nu}$ such that for any $\gamma \in \cS_{0,T}$, we have  $\,V^{0,\nu}_\g =  \b1_{\{\g \ge \nu\}}  V\big(\t(\nu) \land \g \big), \, \,\pas$ We show that the stopping time $
 \, \t_V(\nu)  \dfnn \inf \big\{ t \in [\nu, T] :\, V^{0,\nu}_t=Y_{t} \big\}\, $
attains the infimum in \eqref{eq:defn-op-rm-intro}. Finally, we construct a saddle point of the stochastic game in \eqref{eq:intro-value}.

  The discrete-time optimal stopping problem for coherent risk measures was studied by \cite[Section 6.5]{Follmer_Schied_2004} and \cite[Sections 5.2 and 5.3]{CDK-2006}.  \cite{Delbaen_2006} and \cite{Kara_Zam_2006},
on the other hand, considered continuous-time optimal stopping problems in which the essential infimum over the stopping times in \eqref{eq:defn-op-rm-intro} is replaced by an essential supremum. The controller-and-stopper problem of \cite{Lepeltier_1985} and \cite{Kara_Zam_2008}, and the optimal stopping for non-linear expectations in \cite{OSNE}, are the closest  in spirit  to our work. However, since our assumptions    on the random function $f$ and the set $\mathcal{Q}_{\nu}$ are dictated by the representation theorem
for   dynamic convex risk measures, the results in these papers cannot be directly applied. In particular, because of the integrability assumption that appears in the definition of $\mathcal{Q}_{\nu}$ (subsection \ref{sec:notation}), this set may not be closed under \emph{pasting}; see Remark~\ref{rem:pasting}. Moreover, the extant results on controller-and-stopper games would require that $f$ and the $\theta^{Q}$'s be bounded.
We overcome these technical difficulties by using approximation arguments
which rely on \emph{truncation} and \emph{localization} techniques. On the other hand, in finding a saddle point \cite{Kara_Zam_2008} used the weak compactness of the collection of probability measures, in particular the boundedness of $\theta^{Q}$'s. We avoid making this assumption by using techniques from Reflected Backward Stochastic Differential Equations (RBSDEs). In particular, using a comparison theorem and the fact that $V$ can be approximated by solutions of BSDEs with Lipschitz generators, we show that $V$ solves a quadratic RBSDE (QRBSDE). The relationship between the solutions of QRBSDEs and the BMO martingales helps us construct a saddle point. We should point out that     the convexity of $f$ is not needed to derive our results; cf.$\,$Remark~\ref{rem:con-n-nec}.

The layout of the paper is simple. In Section~\ref{sec:dcrm}  we
recall the definition of the dynamic convex risk measures and a
representation theorem. We solve the optimal stopping problem in Section~\ref{sec:mainresults}. In Section~\ref{sec:saddle} we find a saddle point for the stochastic controller-and-stopper game in \eqref{eq:intro-value}. The proofs of our results are given in Section~\ref{sec:Proofs}.

\subsection{Notation and Preliminaries} \label{sec:notation}

 \bs
 Throughout this paper we let $B$ be a $d$-dimensional Brownian
 Motion defined on the probability space $(\O,\cF, P)$,
 and consider the augmented filtration generated by it, i.e.,
 \beas
 \bF= \big\{\cF_t \dfnn \si\big(B_s; s\in [0,t]\big)\vee \cN \big\}_{t
 \ge 0}, \hb{ where $\cN$ is the collection of all $P$-null sets in $\cF$.}
 \eeas
 We fix a finite time horizon $T>0$,  denote by $\sP$ (resp. $\widehat{\sP}$)  the
predictably (resp. progressively) measurable $\si$-field on $\O\times [0,T]$, and let
$\cS_{0,T}$ be the set of all $\bF$-stopping times $\nu$ such that
$0\leq\nu\leq T$, \pas~ From now on, when writing $\nu \le \g$, we
always mean two stopping times $ \nu, \g \in \cS_{0,T}$ such that
$\nu \le \g$, \pas~ For any $\nu \le \g$  we define $\cS_{\nu, \g}
\dfnn \{ \si \in \cS_{0,T}\, |\; \nu \le \si \le \g, ~\pas \}$
 and let $\cS^\star_{\nu, \g}$ denote all finite-valued stopping times in $\cS_{\nu, \g}\,$.

\ms  The following spaces of functions will be used in the sequel:

\ms \no $\bullet$   Let $\,\cG\,$ be a generic sub-$\si$-field of $\,\cF\,$.  $\mathbb{L}^0(\cG )$ denotes  the space of all real-valued,
$\,\cG-$measurable random variables.

\ms \no $\bullet$  $\mathbb{L}^\infty(\cG ) \dfnn \{\xi \in \mathbb{L}^0(\cG ):\,\|\xi\|_\infty \dfnn
\underset{\o \in \O}{\esssup}\,|\xi(\o)| <\infty\}$.

\ss \no $\bullet$ $\mathbb{L}^0_\bF[0,T]$ denotes the space of all real-valued,
$\bF$-adapted processes.

 \ss \no $\bullet$ $\mathbb{L}^\infty_\bF[0,T]  \dfnn  \big\{  X \in \mathbb{L}^0_\bF[0,T] :\,  \|X\|_\infty \dfnn \esup{(t,\o)
\in [0,T] \times \O}  |X_t(\o)| <\infty \big\}$.

  \no $\bullet$ $\hC^p_\bF[0,T]   \dfnn  \big\{  X \in \mathbb{L}^p_\bF[0,T] :\, \hb{$X$ has continuous
paths}\}$, $~~~p=0, \infty$.

 \ms \no $\bullet$ $\hC^2_\bF[0,T] \dfnn   \Big\{  X \in \hC^0_\bF[0,T] :\,    E\Big(\,\underset{t \in [0,T]}{\sup}|X_t|^2\Big) <\infty \Big\}$.

  \ss \no $\bullet$  $\hH^2_\bF([0,T];\hR^d) $\,\big(resp. $\widehat{\hH}^2_\bF([0,T];\hR^d)$\big) denotes the space of all $\hR^d-$valued,
 $\bF-$adapted predictably (resp. progressively) measurable processes $X$ with  $   E \neg\int_0^T \neg |X_t|^2  dt <\infty$.

 \ms \no $\bullet$  $\hH^\infty_\bF([0,T];\hR^d) $ denotes the space of all $\hR^d$-valued,
 $\bF$-adapted predictably measurable processes $X$ with  \\ $   \esup{(t,\o)
\in [0,T] \times \O}  |X_t(\o)| <\infty $.

 \ms \no $\bullet$  $  \hK_\bF[0,T] $ denotes the space of all real-valued, $\bF$-adapted continuous increasing processes $K$ with $K_0=0$.

\ms  Let us consider  the set $\cM^e$ of all probability measures on $(\O,
\cF)$ which are equivalent to $P$. For any $Q \in \cM^e $, it is well-known that there is
 an $\hR^d-$valued predictable process $\th^Q  $ with $\int_0^T |\th^Q_t|^2 dt < \infty$, \pas,
 such  that the density process $Z^Q$ of $Q$ with respect to $P$ is     the stochastic exponential of $\,\th^Q\,$, namely,
   \beas
   Z^Q_t= \sE\left(\th^Q \bullet B \right)_t = \exp\left\{ \int_0^t
\th^Q_s dB_s- \frac12  \int_0^t \big|\th^Q_s\big|^2 d s \right\}, \quad 0 \le t \le T\,.
 \eeas
 We denote
 $ \,Z^Q_{\nu,\g} \dfnn  Z^Q_\g / Z^Q_\nu  =  \exp\left\{ \int_\nu^\g
\th^Q_s dB_s- \frac12  \int_\nu^\g \big|\th^Q_s\big|^2 d s \right\}$
 for any $\nu \le \g$. Moreover, for any $\nu \in \cS_{0,T}$ and  with the notation  $\,\[0, \nu\[ \; \dfnn \{(t, \omega) \in [0,T] \times \Omega: 0
\leq t < \nu(\omega) \}\,$ for the   stochastic interval, we define
 \beas
 \cP_\nu &\dfnn& \big\{Q \in \cM^e :\, Q
 = P \hb{ on }\cF_\nu  \big\}=\big\{Q \in \cM^e :\,  \th^Q_t(\o) =0,
~ \dtp \hb{ on $\[0, \nu\[$} \big\}\,,  \\
    \q
   \cQ_\nu &\dfnn& \Big\{Q \in \cP_\nu :\, E_Q     \int_\nu^T f
\neg \left(s, \th^Q_s\right) \neg ds    < \infty \Big\}\,.
 \eeas

\section{Dynamic Convex Risk Measures}\label{sec:dcrm}


\begin{deff}  A {\rm dynamic convex risk measure} is a family of functionals
$\big\{\rho_{\nu, \g}: \mathbb{L}^\infty(\cF_\g) \rightarrow
\mathbb{L}^\infty(\cF_\nu)\big\}_{\nu \le \g}$ which satisfy the following properties:
For any stopping times $\,\nu \le \g\,$ and any $\, \mathbb{L}^\infty(\cF_\g)-$measurable random variables $\,\xi,\,\eta \, $, we have
 \bi

 \item  {\it ``Monotonicity":}\;  $ \rho_{\nu, \g}(\xi) \le  \rho_{\nu, \g}(\eta) $,
\pas~if $\,\xi \ge \eta\,$, \pas

\item  {\it ``Translation Invariance":}\;
$ \rho_{\nu, \g}(\xi+\eta)  = \rho_{\nu, \g}(\xi) - \eta$, \pas~ if
$\eta \in \mathbb{L}^\infty(\cF_\nu)$.

\item  {\it ``Convexity":}\;
$ \rho_{\nu, \g}\big(\l \xi+(1-\l)\eta\big) \le    \l \rho_{\nu, \g}
( \xi ) + (1-\l) \rho_{\nu, \g} ( \eta ) $, \pas~for any $\l \in
(0,1)$.

\item  {\it ``Normalization":}\;
$\rho_{\nu, \g}(0)=0$, \pas
 \ei
\end{deff}

\ss  \cite{DPR_2009} provide a representation result, Proposition \ref{prop_represent} below, for   dynamic convex risk measures $\big\{\rho_{\nu, \g} \big\}_{\nu \le \g}$ that satisfy  the following properties:
 \bi
 \item[ {\bf (A1)} ] {\it ``Continuity from above":}\; For any decreasing sequence
 $\{\xi_n\} \subset \mathbb{L}^\infty(\cF_\g) $
 with $\xi \dfnn \lmtd{n \to \infty} \xi_n \in \mathbb{L}^\infty(\cF_\g)$,
 it holds \pas~ that $ \lmtu{n \to \infty} \rho_{\nu, \g}(\xi_n) = \rho_{\nu, \g}(\xi)$.

 \item[ {\bf (A2)} ] {\it ``Time Consistency":}\; For any $\si \in \cS_{\nu, \g}$ we have: $\rho_{\nu, \si}\big(-\rho_{\si, \g}( \xi )\big)=\rho_{\nu, \g}( \xi
)$, \pas

\item[ {\bf (A3)} ] {\it ``Zero-One Law":}\; For any $A \in \cF_\nu$, we have: $\rho_{\nu, \g}(\b1_A \xi) = \b1_A \, \rho_{\nu, \g}(\xi)
$, \pas

\item[ {\bf (A4)} ]  $~~\einf{\xi \in \cA_t} E_P\,[\xi|\cF_t]=0$,
where $\cA_t \dfnn \{\xi \in \mathbb{L}^\infty(\cF_T):\, \rho_{t,T}(\xi) \le
0\}$.
 \ei

 \begin{prop} \label{prop_represent}
 Let $\big\{\rho_{\nu, \g}
\big\}_{\nu \le \g}$ be a dynamic convex risk measure satisfying
(A1)-(A4).  Then for any $\nu \le \g $ and $\xi \in
\mathbb{L}^\infty(\cF_\g)$, we have
 \bea \label{RM_represent}
 \rho_{\nu, \g}(\xi)\,=\,\esup{Q \in \cQ_\nu}\, E_Q \left[-\xi- \int_\nu^\g
  f \neg \left(s, \th^Q_s\right) \neg ds \,\Big| \, \cF_\nu
 \right]\,, \quad \pas
 \eea
Here $f: [0,T] \times \O \times \hR^d \rightarrow  [0, \infty]$ is a suitable
measurable  function, such that

 \ss \no  $(\mathfrak{f}\,1)$ $~f (\cd,\cd,z)$ is predictable for
 any $z \in \hR^d\,;$

 \ss \no  $(\mathfrak{f}\,2)$ $~f (t, \o, \cd)$ is proper
 convex, and
 lower semi-continuous for $\,$\dtp ~ $(t, \o) \in [0,T] \times \O\,;$ and

 \ss \no  $(\mathfrak{f}\,3)$  $~f (t, \o, 0)=0 $, ~\dtp

\end{prop}

 \smallskip We refer to \cite{rock}, p.$\,$24  for the notion of ``proper convex
 function", and  
 review  some basic properties of the essential extrema as in  
\cite[Proposition VI-\b1-1]{Neveu_1975} or \cite[Theorem
A.32]{Follmer_Schied_2004}.

 \begin{lemm}
 \label{lem_ess}
Let $\{\xi_i\}_{i \in \cI}$ and $\{\eta_i\}_{i \in \cI}$ be two
classes of $\cF$-measurable random variables with the same index set
$\cI$.

\ss \no  (1)  If $\xi_i \le (=)\; \eta_i$, \pas~ holds for all $i \in
\cI$,
 then $\underset{i \in \cI}{\esssup}\, \xi_i \le (=)\; \underset{i \in \cI}{\esssup}\, \eta_i$, \pas

\ss \no  (2)  For any $A \in \cF$, it holds \pas~that
 $\;\underset{i \in \cI}{\esssup}\, \big( \b1_A \xi_i + \b1_{A^c} \eta_i
\big) = \b1_A \,\underset{i \in \cI}{\esssup}\, \xi_i + \b1_{A^c}\,
\underset{i \in \cI}{\esssup}\, \eta_i$. In particular, $\underset{i
\in \cI}{\esssup}\, \big( \b1_A \xi_i  \big) =  \b1_A\, \underset{i
\in \cI}{\esssup}\, \xi_i $, \pas

\ss \no (3) For any $\cF$-measurable random variable $\g$ and any
 $\l>0$, we have $\underset{i \in \cI}{\esssup}\,  (\l \xi_i + \g )
 = \l \, \underset{i \in \cI}{\esssup}\, \xi_i + \g $, \pas

\ss \no Moreover, (1)-(3) hold when we replace $\,\underset{i \in
\cI}{\esssup}\,$  by $\,\underset{i \in \cI}{\essinf}\,$.
\end{lemm}

 \section{The Optimal Stopping Problem  }\label{sec:mainresults}


In this section we study the optimal stopping problem for dynamic
convex risk measures. More precisely, given $\nu \in \cS_{0,T}$, we
seek an optimal stopping time $\t_*(\nu) \in \mathcal{S}_{\nu, T}\,$ that
satisfies \eqref{eq:defn-op-rm-intro}. We shall assume throughout that the reward process $Y \in \mathbb{L}^\infty_\bF[0,T]$ is   right-continuous and {\it $\cQ_0-$quasi-left-continuous:} to wit, for any
  increasing sequence $\{\nu_n\}_{n \in \hN} $ in $\cS_{0,T}$ with
  $\,\nu \dfnn \lmtu{n \to \infty} \nu_n \in \cS_{0,T}\,$, and any
 $\,Q \in \cQ_0\,$, we have
   \beas
   \linf{n \to \infty} E_Q[Y_{\nu_n}|\cF_{\nu_1}] \le
   E_Q[Y_{\nu}|\cF_{\nu_1}], \q \pas
   \eeas
 In light of the representation \eqref{RM_represent}, we can
alternatively express \eqref{eq:defn-op-rm-intro} as a {\it robust
optimal stopping problem,} in the following sense:
 \bea \label{eq:defn-op-rm2}
  \underset{ \g  \in  \cS_{\nu,T}}{\esssup} \left( \underset{Q \in \cQ_\nu}{\essinf}\, E_Q \neg \,\Big[Y_\g  +\neg \int_\nu^\g \neg f \neg \left(s, \th^Q_s\right) \neg ds \,\Big|\,\cF_\nu\Big] \right)\,=\,\underset{Q \in \cQ_\nu}{\essinf}\, E_Q\Big[Y_{\t_*(\nu)}  + \int_\nu^{\t_*(\nu)} f \neg \left(s, \th^Q_s\right) \neg ds \, \Big|\, \cF_\nu  \Big].
 \eea

 \begin{rem} \label{rem:con-n-nec}
   We will  study the robust
optimal stopping problem \eqref{eq:defn-op-rm2} in a  setting more general than alluded to heretofore: From now on, we only assume that $f: [0,T] \times \O \times \hR^d \rightarrow  [0, \infty]$ is
 a  $ \sP \otimes \sB(\hR^d)/ \sB([0, \infty])$-measurable function which  satisfies \($\mathfrak{f}\,3$\); i.e., the convexity \($\mathfrak{f}\,2$\) is not necessary for solving \eqref{eq:defn-op-rm2}.
  \end{rem}

   In order to find a stopping time which is  optimal, i.e.,   attains the essential supremum in \eqref{eq:defn-op-rm2},
  we introduce the lower- and upper-value, respectively,  of the stochastic game suggested by \eqref{eq:defn-op-rm2}, to wit, for every
  $\nu \in \cS_{0,T}\,$:
 \beas
 \q \; \ul{V}(\nu)  \,\dfnn \, \underset{\g \in \cS_{\nu,T}}{\esssup}
 \left(  \underset{Q \in \cQ_\nu}{\essinf}\,E_Q \neg \Big[Y_\g
 +\neg \int_\nu^\g \neg f \neg \left(s, \th^Q_s\right) \neg ds \,\Big|\,\cF_\nu\Big] \right)
 ,~~\;\;  \ol{V}(\nu) \, \dfnn \, \underset{Q \in \cQ_\nu}{\essinf} \neg \left(\, \underset{\g \in \cS_{\nu,T}}{\esssup}\,
  E_Q \neg \,\Big[Y_\g  +\neg \int_\nu^\g \neg f \neg \left(s, \th^Q_s\right) \neg ds \,\Big|\,\cF_\nu\Big]
  \right).
 \eeas
  In Theorem \ref{V_process}  we shall show that  the quantities  $\ul{V}(\nu)$ and $\ol{V}(\nu)$ coincide    at any $\nu \in
\cS_{0,T}$, i.e., a min-max theorem holds; we shall also identify two optimal stopping times  in Theorems~\ref{V_process} and \ref{V_RC}, respectively.

  \smallskip    Given any   probability measure  $\,Q \in \cQ_0\,$, let us  introduce for each fixed  $\,\nu \in \cS_{0,T}$ the quantity
 \bea \label{def_R_Q}
 R^Q(\nu) \,\dfnn\,
   \underset{\z \in \cS_{\nu,T}}{\esssup}\,
  E_Q \Big[ \,Y_\z \neg +\int_\nu^\z f \neg \left(s, \th^Q_s\right)  \neg ds  \, \Big|\,\cF_\nu  \Big]
  = \underset{\si \in \cS_{0,T}}{\esssup}\,
  E_Q \Big[\,Y_{\si \vee \nu} \neg +\int_\nu^{\si \vee \nu}
  f \neg \left(s, \th^Q_s\right) \neg ds \, \Big|\,\cF_\nu \Big]\,
   \ge Y_\nu
 \eea
 and recall from the classical theory of optimal stopping (see \cite{El_Karoui_1981} or \cite[Appendix D]{Kara_Shr_MF})    the following result.

 \begin{prop} \label{Prop3.1}
Fix a probability measure  $Q \in \cQ_0$.

\noindent (1) The process $\big\{  R^Q( t)\big\}_{t \in [0,T]}$ admits an RCLL
 modification $R^{Q,0}$ such that, for any $\nu \in \cS_{0,T}\,,$  we have
 \bea
 R^{Q,0}_\nu = R^Q(\nu), \q \pas~\label{eqn-a03}
 \eea
 (2) For every $\nu \in \cS_{0,T}\,,$  the stopping time $ \t^Q(\nu) \dfnn \inf\{t \in [\nu, T]
:\, R^{Q,0}_t=Y_t\} \in \cS_{\nu, T}$ satisfies for any $\g \in \cS_{\nu, \t^Q(\nu)}:$
 \bea
   R^Q(\nu)  &=&   E_Q \Big[ Y_{\t^Q(\nu) }+   \int_\nu^{\t^Q(\nu)} f \neg \left(s, \th^Q_s\right) \neg ds
\,    \Big| \,\cF_\nu \Big]   \,= \,  E_Q \Big[  R^Q \big(\t^Q(\nu)\big)+
  \int_\nu^{\t^Q(\nu)} f \neg \left(s, \th^Q_s\right) \neg ds  \, \Big| \,\cF_\nu\Big] \nonumber \\
  & =&   E_Q \Big[  R^Q (\g) +\int_\nu^\g f \neg \left(s, \th^Q_s\right) \neg ds  \,\Big|\,\cF_\nu\Big], \q
   \pas~\label{eqn-a05}
 \eea
 Therefore, $\t^Q(\nu)$
  is an optimal stopping time for maximizing  $\,\dis E_Q \Big[Y_\z \neg +\int_\nu^\z f \neg \left(s, \th^Q_s\right)
  \neg ds \, \Big|\,\cF_\nu \Big]\,$ over $\,\z \in \cS_{\nu,T}$.
\end{prop}

\ms

For any $\nu \in \cS_{0,T}$ and $k \in \hN $, we introduce the collection of probability measures
   \beas
    \cQ^k_\nu \, \dfnn \, \Big\{ Q \in \cP_\nu : \,
   \big|\th^Q_t(\o)|\vee f \big(t, \o, \th^Q_t(\o)\big)
    \le k, ~ \dtp \hb{ on }\]\nu , T\]  \Big\}.
 \eeas
  \begin{rem} \label{rem_belong}
 It is clear that $\cQ^k_\nu \subset \cQ_\nu \,$; and  from
 $(\mathfrak{ f}\,3)$ one can deduce that for any $\nu \le \g$ we have
 \beas
 \cQ_\g \subset \cQ_\nu \q \hb{and} \q \cQ^k_\g \subset \cQ^k_\nu\,,
 \q \fa ~k \in \hN.
 \eeas
  \end{rem}

Given a $Q \in \cQ_\nu$ for some $\nu \in \cS_{0,T}$, we
\emph{truncate} it in the following way: The predictability of
process $\th^Q$ and Proposition \ref{prop_represent} imply that
$\big\{ f   \big(t,\th^Q_t\big)\big\}_{t \in [0,T]}$ is also
a predictable process. Therefore,  for any given $k \in \hN$, the set
  \bea
  \label{AQNK}
    A^Q_{\nu,k}\, \dfnn \, \Big\{ (t, \o) \in \, \]\nu , T\]:\,
    \big|\th^Q_t(\o)|\vee f\big(t, \o, \th^Q_t(\o)\big)
    \le k \Big\} \,\in \, \sP
  \eea
 is   predictable.   Then the predictable process $\,\th^{Q^{\nu, k}} \dfnn \b1_{A^Q_{\nu,k}}
\th^Q \,$ gives rise to a probability measure $Q^{\nu,k} \in
\cQ^k_\nu$ via the recipe $  \,dQ^{\nu,k} \dfnn \sE\big(\th^{Q^{\nu,
k}} \bullet B \big)_T \, dP\,$.  Let us define the stopping times
    \beas
     \si^Q_m\, \,\dfnn \, \,\inf\left\{ t \in [0,T]:
     \hb{$  \int_0^t \big|\th^Q_s\big|^2 ds  > m$}  \right\} \land T  \, , \q    m \in \hN.
   \eeas
     There exists a null set $N$ such that, for any $\,\o \in \Omega \setminus N\,$, we have  $ \si^Q_m (\omega) = T$ for some $m = m(\o) \in \hN$.
 Since $  \, E\int_0^{ \si^Q_m} \big|\th^Q_t\big|^2 dt   \le m \,$ holds  for each $m \in \hN\,$, we have
   $\big|\th^Q_t(\o)| <\infty,~
   \dtp$ on $\[0,  \si^Q_m\] $.

   \smallskip
   As $\,\Big(\underset{m \in \hN}{\cup}\[0,  \si^Q_m\] \,\Big) \bigcup  \big([0,T] \times   N \big) = [0,T] \times  \O \,   $,    it follows that  $\big|\th^Q_t(\o)| <\infty\,$ holds $\,    \dtp\,$ on $ [0,T] \times \O$.    On the other hand, since $Q \in
\cQ_\nu$ we have $E_Q   \int_\nu^T f \neg \left(s, \th^Q_s\right) \neg
ds   < \infty$, which implies    $   \b1_{
\]\nu , T\]  }   (t,\o)\, f\big(t, \o, \th^Q_t(\o)\big) <
\infty\,$ holds $~dt \otimes  dQ-$a.s., or equivalently \dtp~ Therefore, we see
 that
    \bea \label{eqn-a07}
     \lmtu{k \to \infty} \b1_{ A^Q_{\nu,k}}  =  1_{\]\nu , T\]},\q \dtp
    \eea

 \ms  For any $\nu \in \cS_{0,T}$, the upper value $\ol{V}(\nu)$
 can be approximated from above in two steps,   presented in the next two lemmas.

\begin{lemm} \label{lmt_V}
  Let $\nu \in \cS_{0,T}$.   (1) For any $\g \in \cS_{\nu,T} $   we have
  \bea \label{eqn-a09}
  \underset{Q \in \cQ_\nu}{\essinf}\,  E_Q \Big[Y_\g \neg
 +\int_\nu^\g f \neg \left(s, \th^Q_s\right) \neg ds \,\Big| \, \cF_\nu\Big]
 = \lmtd{k \to \infty}\, \underset{Q \in \cQ^k_\nu}{\essinf}\,
   E_Q \Big[Y_\g \neg +\int_\nu^\g f \neg \left(s, \th^Q_s\right) \neg ds\, \Big|\,\cF_\nu\Big],   \q  \pas
  \eea
   (2) It holds \pas~ that
     \bea \label{eqn-a11}
     \ol{V}(\nu)=\einf{Q \in \cQ_\nu} R^Q(\nu) = \lmtd{k \to \infty}\, \einf{Q \in \cQ^k_\nu}  R^Q(\nu).
     \eea
\end{lemm}

\begin{lemm} \label{lmt_R}
 ~Let $k \in \hN$ and $\,\nu \in \cS_{0,T}\,$.

  \noindent (1) For any $\g \in \cS_{\nu,T} $  there exists a sequence
 $\{Q^{\g,k}_n\}_{n \in \hN} \subset \cQ^k_\nu$ such that
  \bea \label{eqn-a13}
   \underset{Q \in \cQ^k_\nu}{\essinf}\,
   E_Q \Big[Y_\g \neg
 +\int_\nu^\g f \neg \left(s, \th^Q_s\right) \neg ds \,\Big|\,\cF_\nu\Big] \,=\, \underset{n \to \infty}{\lim}\, \dneg \da
  E_{Q^{\g,k}_n} \Big[Y_\g \neg
 +\int_\nu^\g f   \big(s, \th^{Q^{\g,k}_n}_s\big)   ds \,\Big|\,\cF_\nu\Big]\,, \q
 \pas
  \eea
  (2) There exists a sequence
 $\{Q^{(k)}_n\}_{n \in \hN} \subset \cQ^k_\nu$ such that
  \bea  \label{eqn-a15}
   \underset{Q \in \cQ^k_\nu}{\essinf}\,
  R^Q(\nu) = \underset{n \to \infty}{\lim} \,\dneg \da R^{Q^{(k)}_n}(\nu), \q \pas
  \eea
\end{lemm}

  Let us fix $\nu \in \cS_{0,T}$. For any $k \in \hN$, the infimum of the
family $\{\t^Q(\nu)\}_{Q \in \cQ^k_\nu}$ of optimal stopping times
can be approached by a decreasing sequence in this family. As a
result, the infimum is also a stopping time.

\begin{lemm} \label{lmt_tau_nu}
\bs Let $\nu \in \cS_{0,T}$ and $k \in \hN$. There exists a sequence
$ \{Q^{(k)}_n\}_{n \in \hN} \subseteq \cQ^k_\nu$ such that
  \beas
  \t_k (\nu) \,\dfnn \, \einf{Q \in \cQ^k_\nu} \, \t^Q(\nu)
  \, =\, \lmtd{n \to \infty} \t^{Q^{(k)}_n}(\nu), \q  \pas
  \eeas
in the notation of Proposition \ref{Prop3.1}, thus $ \,\t_k(\nu) \in \cS_{\nu,T}$.
  \end{lemm}
 Since $\left\{\cQ^k_\nu \right\}_{k \in \hN}$ is an
increasing sequence, $ \big\{\t_k(\nu)  \big\}_{k \in \hN} $  is in
turn a decreasing sequence. Hence
 \bea \label{defn_tau_nu}
  \t(\nu)\, \dfnn \,\lmtd{k \to \infty}
\t_k(\nu)
 \eea
  defines a stopping time  in  $\cS_{\nu,T}$.
The family of stopping times $\{\t(\nu)\}_{\nu \in \cS_{0,T}}$ will
play a crucial role in this section.

\ms The next lemma is concerned with the \emph{pasting} of two
probability measures.

\begin{lemm} \label{lem_paste}
Given $\nu \in \cS_{0,T}$, let $\,\wt{Q} \in \cQ^k_\nu\,$ for some $k
\in \hN$. For any $\,Q \in \cQ_\nu$ and $\g \in \cS_{\nu,T}$, the
predictable process
  \bea \label{eqn-a17}
  \th^{Q'}_t \,\dfnn \, \b1_{\{t \le \g \}} \th^Q_t  + \b1_{\{t > \g
  \}} \th^{\wt{Q}}_t \,, \q  t \in [0,T]
 \eea
 induces a probability measure $Q' \in \cQ_\nu$ by
$  \,  d Q'    \dfnn  \sE \big(\th^{Q'} \bullet B
\big)_T\, d P\,$. If $\,Q\,$ belongs to $\,   \cQ^k_\nu\,$, so does $Q'$.
Moreover, for any $\si
\in \cS_{\g,T}$,
 we have
   \bea
   R^{Q',0}_\si= R^{Q'}(\si)=R^{\wt{Q}}(\si)=R^{\wt{Q},0}_\si \,, \q
   \pas
   \eea
\end{lemm}

\begin{rem}\label{rem:pasting}
The probability measure $Q'$ in Lemma~\ref{lem_paste} is called the {\rm pasting} of $\,Q$ and $\,\wt{Q}$; see e.g. Section 6.7 of  \cite{Follmer_Schied_2004}. In general, $\cQ_\nu$ is not closed under such ``pasting".
\end{rem}

The proofs of the following results use  schemes similar to the ones  in \cite{Kara_Zam_2008}. The main technical difficulty in our case is   mentioned in Remark~\ref{rem:pasting}.
Moreover, in order to use the results of \cite{Kara_Zam_2008} directly, we would have to assume that $f$ and the $\,\theta^{Q}\,$'s are all bounded. We overcome these difficulties by using approximation arguments that rely on \emph{truncation} and \emph{localization} techniques.

 \ms First, we shall show that at any $\nu \in
\cS_{0,T}\,$ we have $\,\ul{V}(\nu)=\ol{V}(\nu)\,$, \pas

\begin{thm}
\label{V_process}
{\bf Existence of Value:}
For any $\nu \in \cS_{0,T}$, we have
 \bea \label{eqn-a19}
 \ul{V}(\nu ) \,=\, \einf{Q \in \cQ_\nu} E_Q \Big[   Y_{\t(\nu)}
    + \int_\nu^{\t(\nu) }  \neg   f \neg \left(s,  \th^Q_s\right) \neg ds
  \,  \Big| \,\cF_\nu\Big]\,=\, \ol{V}(\nu )  \ge Y_\nu\,,  \q \pas
  \eea
Therefore, the stopping time $ \t (\nu)$ of (\ref{defn_tau_nu}) is optimal for the robust optimal stopping problem \eqref{eq:defn-op-rm2} $($i.e., attains the essential infimum there$)$.    
\end{thm}

We shall denote the common value in (\ref{eqn-a19})  by $ \,\,V(\nu)~~\big( =\ul{V}(\nu)=
\ol{V}(\nu)\,\big)\,$.

\begin{prop} \label{V_Y_meet}
For any $\nu \in \cS_{0,T}$, we have $\,V( \t(\nu))=Y_{\t(\nu)}\,$,
\pas~
 \end{prop}

\smallskip
 Note that $\t(\nu)$ may not be the first time after $\nu$ when
the value process   coincides with the reward process. Actually, since the value process $\{V(t)\}_{t \in [0,T]}$ is not necessarily right-continuous, the random time
$\,\inf\{t \in [\nu, T] \;:\; V(t)=Y_t\} \,$ may not even be a stopping
time. We   address this issue in the next three results.

 \begin{prop} \label{V_sub}
Given $\nu \in \cS_{0,T}$, $ Q \in \cQ_\nu$, and $\g  \in
\cS_{\nu, \t(\nu)}$, we have
 \bea
   &&\q E_Q\Big[V (\g)+ \int_\nu^\g f \neg \left(s,\th^Q_s\right) \neg ds \,\Big|\,\cF_\nu\Big]
\,    \ge  \,V(\nu), \q \pas~\hspace{1cm}  \label{eqn-a21}
 \eea
\end{prop}

\begin{lemm} \label{lem_switch}
For any $\nu,\g,\si \in \cS_{0,T}$, we have the $\, \pas\,$ equalities
 \bea \label{eqn-a23}
  \b1_{\{ \nu =\g \}}\,
  \einf{Q \in \cQ_\nu}  E_Q \Big[  Y_{\si \vee \nu}
  + \int_\nu^{\si \vee \nu} f \neg \left(s,\th^Q_s\right) \neg ds \,\Big|\,\cF_{ \nu}\Big] \,=\,  \b1_{\{ \nu =\g \}}\,
  \einf{Q \in \cQ_\g}  E_Q \Big[  Y_{\si \vee \g}
  + \int_\g^{\si \vee \g} f \neg \left(s,\th^Q_s\right) \neg ds  \, \Big|\,\cF_{  \g}\Big]
 \eea
and 
  \bea \label{eqn-a25}
  \b1_{\{ \nu =\g \}} V(\nu)\,=\,\b1_{\{ \nu =\g \}} V(\g) .
  \eea
\end{lemm}

\smallskip Next, we show that for any given $\nu \in \cS_{0,T}$, the process
   $\left\{ \b1_{\{t \ge \nu\}}  V\big(\t(\nu) \land t\big) \right\}_{t \in [0,T]} $ admits an RCLL  modification $V^{0,\nu}$. As a consequence, the first time after $\nu$ when the process $V^{0,\nu}$ coincides with the process $Y$,  is an optimal stopping time for the robust optimal stopping problem \eqref{eq:defn-op-rm2}.   

    \begin{thm} \label{V_RC} {\bf Regularity of the Value:}
 \noindent Let us fix a stopping time $\nu \in \cS_{0,T}$.

 \noindent (1) The  process
   $ \left\{ \b1_{\{t \ge \nu\}}  V\big(\t(\nu) \land t\big) \right\}_{t \in [0,T]}
    $ admits an RCLL modification $V^{0,\nu}$ such that, for any $\g \in \cS_{0,T}\,:$  
    \bea  \label{V0}
    V^{0,\nu}_\g \,=\,  \b1_{\{\g \ge \nu\}}  V\big(\t(\nu) \land \g \big) \,, \q  \pas
    \eea
  (2) Consequently,
    \bea  \label{eqn-a27}
  \t_V\neg(\nu) \, \dfnn \, \inf\left\{ t \in [\nu, T ] :\, V^{0,\nu}_t=Y_t
  \right\}
   \eea
 is a stopping time which, in fact, attains the essential infimum in \eqref{eq:defn-op-rm2}.   
 \end{thm}

We  should point out that, in order to determine the optimal stopping time in \eqref{eq:defn-op-rm-intro},   knowledge of  the function $f$ in the representation \eqref{RM_represent} is not necessary. Indeed, let the $\rho-$Snell envelope be the RCLL modification of   $ \underset{ \g  \in  \cS_{\nu,T}}{\esssup}\left(-\rho_{\nu, \g}\left(Y_\g\right)\right)$, $\nu \in \cS_{0,T}$. From our results above, the first time after $\nu$ that the $\rho$-Snell envelope touches the reward process $Y$   is an optimal stopping time; this is consistent with    the classical theory of optimal stopping.

\section{The Saddle Point Problem}\label{sec:saddle}


   In this section  we will contruct a saddle point of the stochastic game in \eqref{eq:intro-value}. As in the previous section, we shall assume here that
  $f: [0,T] \times \O \times \hR^d \rightarrow  [0, \infty]$ is a $ \sP \otimes \sB(\hR^d)/ \sB([0, \infty])-$measurable function which satisfies \($\mathfrak{f}\,3$\).  For any given $Q \in \cQ_0  $ and $\nu \in  \cS_{0,T}\,$,  let us denote
    \beas
     Y^Q_\nu\, \dfnn \,Y_\nu+ \int_0^\nu f(s, \th^Q_s)ds   \,\, ~~\q \hb{and} \,\,~~\q   V^Q(\nu)\, \dfnn \,V(\nu)+ \int_0^\nu f(s, \th^Q_s)ds\,.
    \eeas

\begin{deff}
  A pair $(Q^*, \si_*) \in \cQ_0 \times \cS_{0,T}$ is called a {\rm saddle point}, if for every $Q \in \cQ_0$ and $\nu \in \cS_{0,T}$ we have
  \bea \label{defn_saddle}
  E_{Q^*} \big(   Y^{Q^*}_\nu  \big) \le E_{Q^*}\big(Y^{Q^*}_{\si_*}   \big) \le E_Q\big(Y^Q_{\si_*}  \big)\,.
  \eea
\end{deff}

\begin{thm} \label{Thm-saddle} {\bf Necessary Conditions for a Saddle Point:}
 A pair $(Q^*, \si_*) \in \cQ_0 \times \cS_{0,T}$ is a  saddle
 point,  if the following conditions are satisfied:

  \bi
 \item[\(i\)] $Y_{\si_*} = R^{Q^*}(\si_*)$, \pas;

 \item[\(ii\)]  for any $Q \in \cQ_0\,$, we have $\,V(0) \le E_Q \left[V^Q(\si_*) \right]\,$;

\item[\(iii\)]   for any $\nu \in \cS_{0, \si_*}\,$, we have $\, V^{Q^*}(\nu)=E_{Q^*}\left[V^{Q^*}(\si_*) \big|\cF_\nu\right]\,$, \pas
 \ei

 \end{thm}

   To construct a saddle point,  we need the following two notions.

 \begin{deff}
 We call  $\,\cZ \in \widehat{\hH}^2_\bF([0,T];\hR^d)$ a BMO (short for Bounded Mean Oscillation) process if
  \beas
    \|\cZ\|_{BMO} \, \, \dfnn \,\sup_{\t\in\cM_{0,T}}
\left\|E\Big[ \int_\t^T |\cZ_s|^2 ds \, \Big| \,\cF_\t \Big]^{1/2}\right\|_\infty < \infty.
  \eeas
 When $\cZ$ is a BMO process, $\cZ\bullet
B $ is a BMO martingale; see e.g. \cite{expM_BMO}.
  \end{deff}

 \begin{deff}  \label{defn_RBSDE} {\bf BSDE with Reflection:}
 Let $h:  [0,T] \times \O \times \hR \times \hR^d \rightarrow  \hR $ be a $\widehat{\sP} \times \sB(\hR) \times   \sB(\hR^d)/\sB(\hR)$-measurable function.     Given $S \in \hC^0_\bF[0,T]$ and   $\xi \in \hL^0(\cF_T)  $ with $\xi \ge S_T$, \pas,
  a triple $(\G, \cZ,K) \in \hC^0_\bF[0,T] \times \widehat{\hH}^2_\bF([0,T];\hR^d) \times \hK_\bF[0,T]$
  is called a solution to the {\rm reflected backward stochastic differential equation} with terminal condition $\xi$, generator $h$, and obstacle $S$
 \neg  \big(RBSDE $(\xi, h, S)$ for short\big), if   \pas, ~ we have the comparison
  \beas
   S_t \le \G_t= \xi  + \int_t^T h(s, \G_s, \cZ_s) \, ds +   K_T - K_t - \int_t^T  \cZ_s dB_s\,, \qquad \,\,\,~t \in [0, T]\, ,
  \eeas
  and the so-called {\rm flat-off  condition}
    \beas
       \int_0^T    \b1_{\{  \G_s   >  S_s  \}}   d K_s =0, \q \pas
    \eeas
     \end{deff}

   In the rest of this section  we shall assume that the reward process $Y \in \mathbb{L}^\infty_\bF[0,T]$ is continuous and that the function $f: [0,T] \times \O \times \hR^d \rightarrow  [0, \infty]$  satisfies the following additional conditions:

 \ss \no  {\bf (H1)} For every $\,(t,\o) \in [0, T] \times \O\,$, the mapping $\, z \mapsto f(t, \o, z)$ is continuous.

\ss \no   {\bf (H2)} It holds \dtp~ that
  \beas
     f(t, \o, z)  \ge \e \big|z-\U_t(\o)\big|^2 -\ell \,,  \q \fa z \in \hR^d\,.
  \eeas
  Here  $\e>0$ is a real contant,  $\U $ is an $\hR^d-$valued process
with $\|\U\|_\infty \dfnn   \esup{(t,\o)
\in [0,T] \times \O}  |\U_t(\o)| <\infty$, and  $\ell \ge \e \|\U\|^2_\infty$.

 \ss \no   {\bf (H3)} For any $(t,\o, u) \in [0, T] \times \O \times \hR^d$,
 the mapping $\,z \mapsto  f(t, \o, z)+\lan u, z \ran\,$ attains its infimum over $\hR^d$ at some $z^*=z^*(t,\o,u) \in \hR^d$, namely,
  \bea  \label{eqn-c020}
  \wt{f} (t, \o, u) \dfnn  \underset{z \in \hR^d}{\inf} \big( f(t, \o, z)+ \lan u, z \ran \big)= f(t, \o, z^*(t,\o,u))+\lan u, z^*(t,\o,u) \ran  .
  \eea
   Without loss of generality, we can assume that the mapping $z^*: [0, T]
\times \O \times \hR^d \rightarrow \hR^d$ is $\sP \otimes
\sB(\hR^d)/\sB(\hR^d)$-measurable thanks to the Measurable Selection Theorem (see e.g. Lemma 1 of  \cite{Benes_1970} or  Lemma 16.34 of
\cite{Elliott_1982}\big).  We further assume that there exist a non-negative BMO process $\p$ and a $M>0$ such that for \dtp~$(t, \o) \in [0, T] \times \O$
    \beas
       |z^*(t,\o,u)| \le  \p_t(\o) + M|u|, \q  \fa u \in \hR^d.
    \eeas

\begin{eg} {\rm
 Let $\l \ge 0$ and  let  $\L, \U \in \hH^\infty_\bF([0,T];\hR^d) $ with $\L_t (\o) \ge \e>0$, \dtp\,  Define
 \beas
 f (t,\o,z) \dfnn  \L_t(\o) \left(\big|z-\U_t(\o) \big|^{2+\l} - \big| \U_t(\o) \big|^{2+\l} \right) ,  \q  \fa  (t,\o,z) \in [0,T] \times \O \times \hR^d.
 \eeas
 Clearly, $f^+=f \vee 0$ is a $ \sP \otimes \sB(\hR^d)/ \sB([0, \infty])$-measurable function that satisfies  \($\mathfrak{f}\,3$\) and {\bf (H1)}. It turns out that   $f^+$ satisfies {\bf (H2)}, since  \dtp~ we have that
  \beas
 \q     f^+(t, \o, z) &\ge&  f(t, \o, z) \ge  \L_t(\o)\left(\big|z-\U_t(\o) \big|^2-1 \right) -\L_t(\o)\big| \U_t(\o) \big|^{2+\l} \\
    & \ge & \e \big|z-\U_t(\o) \big|^2  - \| \L \|_\infty\left(1+  \| \U \|_{\infty}^{2+\l} \right),  ~\fa z \in \hR^d\,.
  \eeas
  For any $(t,\o, u ) \in [0, T] \times \O \times \hR^d  $ the gradient
 \beas
 \nabla_z   \big( f(t, \o, z)+\lan u, z \ran \big)   =  (2+\l) \L_t(\o)\big|z-\U_t(\o) \big|^{\l} \big(z-\U_t(\o)\big) + u, \q \fa  z \in  \hR^d,
 \eeas
 is null only at $ \hat{z} (t,\o,u) =- \big[(2+\l)\L_t(\o) \big]^{-\frac{1}{1+\l}}|u|^{-\frac{\l}{1+\l}}u+\U_t(\o) $, where the mapping $z \rightarrow  f(t, \o, z)+\lan u, z \ran$ attains its infimum over $\hR^d$. When $|u|\ge r_t(\o) \dfnn (2+\l)\L_t(\o) |\U_t(\o)|^{1+\l}$,
  $\hat{z} (t,\o,u) \in A \dfnn \{z \in \hR^d: \,  |   z -\U_t(\o)  | \ge   | \U_t(\o)  | \} $. It  follows that
   \bea
       \underset{z \in \hR^d}{\inf}  \left( f^+(t,\o, z) + \lan u,  z\ran\right) & \le &  f^+(t,\o, \hat{z} (t,\o,u)) + \lan u,  \hat{z} (t,\o,u)\ran
       = f(t,\o, \hat{z} (t,\o,u)) + \lan u,  \hat{z} (t,\o,u)\ran \nonumber  \\
       &=& \underset{z \in \hR^d}{\inf}  \left( f(t,\o, z) + \lan u,  z\ran\right) \le  \underset{z \in \hR^d}{\inf}  \left( f^+(t,\o, z) + \lan u,  z\ran\right).
       \label{eqn-d01}
   \eea
   On the other hand, when $|u| < r_t(\o) $ or equivalently $\hat{z} (t,\o,u) \notin A $,
   the gradient $\nabla_z   \big( f(t, \o, z)+\lan u, z \ran \big) \ne 0$ for any $z \in A$,
   which implies that the mapping $z \rightarrow  f(t, \o, z)+\lan u, z \ran$
   can not attain its infimum over $A$ at an interior point of it. Thus
    \beas
     \underset{z \in A}{\inf}  \left( f(t,\o, z) + \lan u,  z\ran\right) = \underset{z \in \pa A}{\inf}  \left( f(t,\o, z) + \lan u,  z\ran\right)
     = \underset{z \in \pa A}{\inf}    \lan u,  z\ran .
    \eeas
     Then it follows that
    \beas
     \underset{z \in \hR^d}{\inf}  \left( f^+(t,\o, z) + \lan u,  z\ran\right)
     =\underset{z \in A^c}{\inf}   \lan u,  z\ran  \land    \underset{z \in A}{\inf}  \left( f(t,\o, z) + \lan u,  z\ran\right)
     = \underset{z \in \ol{A^c}}{\inf}     \lan u,  z\ran.
    \eeas
  The latter infimum is attained uniquely at some $\wt{z}(t,\o,u) \in \ol{A^c}$, which together with \eqref{eqn-d01} implies that
    \beas
    z^*(t,\o,u) = \b1_{\{|u|\ge r_t(\o) \}}   \hat{z} (t,\o,u) +   \b1_{\{|u| < r_t(\o) \}} \wt{z}(t,\o,u) .
    \eeas
     Therefore, $f^+$ satisfies  (H3), since for \dtp~$\,(t,\o) \in [0,T] \times \O\,$ we have
    \beas
     |z^*(t,\o,u)|  &\le&  |\hat{z} (t,\o,u) | +|\wt{z}(t,\o,u)|  \le  \big((2+\l)\e \big)^{-\frac{1}{1+\l}} |u|^{\frac{1}{1+\l}}     + 3\|\U\|_\infty \\
      &\le &  \big((2+\l)\e \big)^{-\frac{1}{1+\l}} |u| +\big((2+\l)\e \big)^{-\frac{1}{1+\l}}+ 3\|\U\|_\infty, \q \fa u \in \hR^d\,.
  \eeas}
   \end{eg}

\begin{rem}
{\rm The ``entropic" risk measure with  {\it risk tolerance coefficient} $ \,r > 0\,$, namely
\beas
 \rho^{\,r}_{\nu,\g}(\xi) \, \dfnn \,  r  \log \Big\{E\big[e^{-\frac{1}{r}\xi}\,\big|\,\cF_\nu\big]\Big\}, \q  \xi  \in
\hL^\infty (\cF_\g),
\eeas
 is a typical example of a dynamic convex risk measures satisfying (A1)-(A4). The corresponding $f$ in \eqref{RM_represent} is $f(z)=\frac{r}{2}|z|^2$, $z \in \hR^d$.}
\end{rem}

   \begin{eg}
  {\rm Let $b^1$, $b^2$ be two real-valued processes such that $ -\varpi \le b^1_t(\o) \le 0 \le b^2_t(\o) \le \varpi $, \dtp~ for some $\varpi >0$
  Let $\vf: [0,T] \times \O \times \hR \to \hR$ be a $\sP \otimes \sB(\hR)/\sB(\hR) $-measurable function that satisfies the following two assumptions:

   \ss  \no \,\,(i)  For any $(t, \o) \in [0,T] \times \O$,   $ \vf(t, \o, \cd)$ is a bijective locally-integrable function or a continuous surjective locally-integrable function on $\hR$.

  \ms  \no (ii) For some $ \e_1 , \e_2 >0$, it holds \dtp ~ that
   \beas
    \vf(t, \o, x) \left\{
     \begin{array}{ll}
     \ge \big(2 \,\e_1 x  +b^1_t(\o)\big)\vee 0, \qq & \hb{if } x > 0, \\
      \le   \big(2 \,\e_2  x +b^2_t(\o)\big) \land 0, \qq & \hb{if } x <  0.
     \end{array}
     \right.
   \eeas
 Then $f(t,\o,z) \dfnn \int_0^z \vf (t,\o, x) dx$, $z \in \hR$ defines
 a $\sP \otimes \sB(\hR)/\sB([0,\infty]) $-measurable non-negative function that satisfies $(\mathfrak{f}\,3)$ and {\bf (H1)}.
 Let $\e = \e_1 \land \e_2$. For \dtp~$(t,\o) \in [0,T] \times \O$,   if $z >0$, then
 \beas
        f(t,\o,z) &\ge& \int_0^z \big(2 \,\e_1 x  +b^1_t(\o)\big) dx = \e_1 z^2  +b^1_t(\o)z \ge \e z^2 - \varpi  z
  = \e\left(z- \frac{\varpi}{2\e}\right)^2 -\frac{\varpi^2}{4\e};
  \eeas
  on the other hand, if   $z < 0$, then
   \beas
        f(t,\o,z) &=& -\int_z^0  \vf (t,\o, x) dx \ge  -\int_z^0  \big( 2\e_2  x +b^2_t(\o)\big)  dx  = \e_2 z^2  +b^2_t(\o)z  \ge \e z^2  +\varpi z  \\
  & =& \frac12\e \left(z - \frac{\varpi}{2\e}\right)^2 +\frac12\e \left(z + \frac{3 \varpi}{2\e} \right)^2 -\frac{5\varpi^2}{4\e} .
  \eeas
Thus it holds \dtp~ that
 $  f(t,\o,z) \ge  \frac12\e \left(z - \frac{\varpi}{2\e}\right)^2   -\frac{5\varpi^2}{4\e}$, i.e.,  {\bf (H2)} is satisfied.

  \ms  For any $(t,\o, u) \in [0, T] \times \O \times \hR  $, since  $ \frac{d}{dz} \big( f(t, \o, z)+ u  z   \big)  =  \vf (t,\o, z)  + u$,  
  the mapping $z \mapsto  f(t, \o, z)+u  z$ attains its infimum over $\hR$ at each $z \in \{ z \in \hR: \, \vf(t,\o, z)=x \} $. Thus
        $ \vf^{-1}_-(t,\o, x) \le z^*(t,\o, u)  \le  \vf^{-1}_+(t,\o, x)$, where
        \beas
            \vf^{-1}_-(t,\o, x) \dfnn \inf \{ z \in \hR: \, \vf(t,\o, z)=x \}  \q \hb{and} \q \vf^{-1}_+(t,\o, x) \dfnn \sup \{ z \in \hR: \, \vf(t,\o, z)=x \}.
         \eeas
   It is clear that $\vf \left(t,\o, \vf^{-1}_-(t,\o, x) \right) =x $ and $\vf \left(t,\o, \vf^{-1}_+(t,\o, x) \right) =x $.
   For \dtp~$(t,\o) \in [0,T] \times \O$ and $u \in \hR$,   if $\vf^{-1}_-(t,\o, x) >0$, then
  \beas
 -u  =\vf \left(t,\o, \vf^{-1}(t,\o, -u) \right)     \ge 2 \,\e_1  \vf^{-1}_-(t,\o, x) +b^1_t(\o),
  \eeas
  which implies that  $     0< \vf^{-1}_-(t,\o, x)  \le \frac{1}{2 \e} \big( |u|+\varpi\big)$. On the other hand, if $\vf^{-1}_-(t,\o, x)  <0$, one can
   deduce that
    $     - \frac{1}{2 \e} \big( |u|+\varpi\big) \le  \vf^{-1}_-(t,\o, x)  < 0  $ by a similar argument. Hence $ \big|\vf^{-1}_-(t,\o, x) \big|  \le \frac{1}{2 \e} \big( |u|+\varpi\big)    $. Similarly, this inequality also holds for  $\vf^{-1}_+(t,\o, x)$, thus for  $z^*(t,\o, u)$. As a result, (H3) is also satisfied. \qed}
 \end{eg}

  One can easily deduce from {\bf (H2)} and $(\mathfrak{f}\,3)$ that \dtp
  \beas
      -    \frac{1+\e}{ 4\e}  |u|^2  -\|\U\|^2_\infty -\ell   \le \wt{f}(t, \o, u)  \le 0 ,  \q  \fa u \in \hR^d,
  \eeas
  which shows that $\wt{f}$ has quadratic growth in $u$. Thanks to Theorems 1 and 3 of \cite{KLQT_RBSDE}, the RBSDE $(Y_T,\wt{f}, Y)$
 admits a solution $(\wt{\G}, \wt{\cZ},  \wt{K}) \in \hC^\infty_\bF[0,T] \times \hH^2_\bF([0,T];\hR^d) \times  \hK_\bF[0,T]$.

In fact, $\wt{\cZ}$ is a BMO process. To see this, we set $\k \dfnn \frac{1+\e}{ 4\e}  \vee \big( \|\U\|^2 +\ell \big) $.
 For any $\nu \in \cS_{0,T}$, applying It\^o's formula to $\dis  e^{-4\k \wt{\G}_t}$ we get
  \beas
  e^{-4\k \wt{\G}_\nu}+8\k^2 \neg \int_\nu^T e^{-4\k \wt{\G}_s}|\wt{\cZ}_s|^2ds
  & \tneg= & \tneg e^{-4\k Y_T}   \neg - 4\k \neg \int_\nu^T e^{-4\k \wt{\G}_s} \wt{f}(s, \wt{\cZ}_s)ds
  - 4\k \neg \int_\nu^T e^{-4\k \wt{\G}_s} d \wt{K}_s+ 4\k \neg \int_\nu^T e^{-4\k \wt{\G}_s}\wt{\cZ}_sdB_s  \nonumber \\
 & \tneg\leq& \tneg  e^{-4\k Y_T} \neg +4\k^2 \neg \int_\nu^Te^{-4\k \wt{\G}_s}\big(1+|\wt{\cZ}_s|^2\big)ds
 + 4\k \neg \int_\nu^T e^{-4\k \wt{\G}_s}\wt{\cZ}_sdB_s .
  \eeas
 Taking   conditional expectations in the above expression, we obtain
  \beas
 e^{-4\k \|  \wt{\G} \|_\infty}    E\Big[  \int_\nu^T |\wt{\cZ}_s|^2 ds \,\Big|\,\cF_\nu\Big]
 \le E\Big[  \int_\nu^Te^{-4\k \wt{\G}_s} |\wt{\cZ}_s|^2 ds\,\Big|\,\cF_\nu\Big]  \le  \frac{1}{4\k^2}
 E\Big[e^{-4\k Y_T} \,\big|\,\cF_\nu \Big]+e^{4\k \|  \wt{\G} \|_\infty}T\,.
 \eeas
 which implies that  $\|\wt{\cZ}\|_{BMO} \le  e^{4 \k \|  \wt{\G} \|_\infty} \big(\frac{1}{4\k^2} + T\big)^{1/2}\,$.

  \ms  Since the mapping $z^*: [0, T]
\times \O \times \hR^d \rightarrow \hR^d$ is $\sP \otimes
\sB(\hR^d)/\sB(\hR^d)$-measurable (see (H3)),
 \bea   \label{eqn-c030}
 \th^*_t(\o)  \dfnn z^*(t,  \o, \wt{\cZ}_t (\o) ), \q (t , \o)  \in [0, T] \times \O
 \eea
 is a predictable process. It follows from {\bf (H3)} that for any $\nu \in [0, T]$
  \beas
  E\Big[\int_\nu^T |\th^*_s|^2 ds\,\Big|\,\cF_\nu\Big]        \le 2  E\Big[\int_\nu^T \p^2_s  ds\,\Big|\,\cF_\nu\Big]
   + 2M^2  E\Big[\int_\nu^T |\wt{\cZ}_s|^2 ds\,\Big|\,\cF_\nu\Big], \q \pas,
  \eeas
 which implies that $\th^*$ is a BMO process.

  \ms Fix $\nu \in \cS_{0,T}$.
 Since $\th^{*,\nu}_t \dfnn \b1_{\{t > \nu\}} \th^*_t $, $t \in [0,T]$ is also a BMO process,
 we know from Theorem 2.3 of \cite{expM_BMO} that the stochastic exponential
 $  \big\{\sE\left(\th^{*,\nu} \bullet B \right)_t \big\}_{t \in [0, T]} $ is a uniformly integrable martingale.
 Therefore,  $  \,d Q^{*,\nu}   \dfnn  \sE\left(\th^{*,\nu} \bullet B \right)_T  \, d P\,$  defines a probability measure $\,Q^{*,\nu} \in   \cP_\nu \,$.
  As
   \beas    
       \wt{f}(s,   \wt{\cZ}_s ) =      f \big(s,    z^*(s,   \wt{\cZ}_s )    \big)+\lan \wt{\cZ}_s , z^*(s,   \wt{\cZ}_s )  \ran
       =      f(s,    \th^*_s )+\lan \wt{\cZ}_s , \th^*_s  \ran ,  \q \dtp
  \eeas
 by \eqref{eqn-c020} and \eqref{eqn-c030} and the Girsanov Theorem, we can deduce
    \bea   \label{eqn-c034}
   \wt{\G}_{\nu \vee t}
   &=&   Y_T  + \int_{\nu \vee t}^T  \left[ f(s,    \th^{*, \nu}_s )+\lan \wt{\cZ}_s , \th^{*, \nu}_s  \ran \right]  ds
   +   \wt{K}_T - \wt{K}_{\nu \vee t} - \int_{\nu \vee t}^T \wt{\cZ}_s dB_s    \nonumber    \\
     &=&    Y_T  + \int_{\nu \vee t}^T f(s,    \th^{*, \nu}_s )   ds
   +   \wt{K}_T - \wt{K}_{\nu \vee t} - \int_{\nu \vee t}^T \wt{\cZ}_s dB^{Q^{*, \nu}}_s    ,\q  t \in [0, T] ,
  \eea
  where $B^{Q^{*, \nu}}$ is a Brownian Motion under $Q^{*, \nu}$. Letting $t=0$ and  taking the expectation $E_{Q^{*, \nu}}$ yield that
    \beas
        E_{Q^{*, \nu}}  \int_\nu^T f \big(s,    \th^{*, \nu}_s \big)  ds   \le   E_{Q^{*, \nu}} \big(\,\wt{\G}_\nu   -   Y_T \big) \le 2\| \wt{\G} \|_\infty\,,
    \eeas
 thus $Q^{*, \nu} \in \cQ_\nu $.  The   lemma below shows that
  $\wt{\G}$  is indistinguishable from $ R^{{Q^{*, \nu}},0}  $ on the stochastic interval $ \[\nu,T\] $.

 \begin{lemm}  \label{G_R_same} Given $\nu \in \cS_{0,T}$,  it holds \pas ~ that
         \bea  \label{eqn-c060}
         \wt{\G}_t  =  R^{{Q^{*, \nu}},0}_t , \q \fa  t \in [\nu, T] .
         \eea
 \end{lemm}

 \ms    Let $k \in \hN$ and $Q \in \cQ^k_\nu$. It is easy to see that the function
 $h_Q(s , \o, z) \dfnn f(s , \o, \th^Q_s(\o))  + \lan z, \th^Q_s(\o)  \ran$ is Lipschitz continuous in $z$: to wit, for \dtp $(t, \o) \in [0, T] \times \O$
  \beas
   \big| h_Q(s , \o, z) -h_Q(s , \o, z')   \big| = \big|\lan z -z',  \th^Q_s \ran \big| \le \big| \th^Q_s\big|\cd|  z -z' |  \le k |  z -z' | , \q \fa z, z' \in \hR^d.
  \eeas
 Moreover, we have
  \beas
   E\int_0^T | h_Q(s , 0)|^2 ds = E\int_0^T | f(s , \th^Q_s ) |^2 ds = E\int_\nu^T | f(s , \th^Q_s ) |^2 ds \le k^2 T.
  \eeas
  Theorem 5.2 of \cite{EKPPQ-1997} assures now that there exists a unique solution $(\G^Q,\cZ^Q,K^Q) \in \hC^2_\bF[0,T] \times   \hH^2_\bF([0,T];\hR^d) \times  \hK_\bF[0,T]$ to the RBSDE$(Y_T,h_Q, Y)$.  Fix $t \in [0, T]$. For any $\g \in \cS_{t,T}$, Girsanov Theorem implies
    \beas
    \G^Q_t &=& Y_T  + \int_t^T   h_Q(s,   \cZ^Q_s  ) \,   ds +   K^Q_T - K^Q_t - \int_t^T \cZ^Q_s dB_s   \\
                & =& \G^Q_\g  + \int_t^\g   f(s, \th^Q_s)   ds +   K^Q_\g - K^Q_t - \int_t^\g \cZ^Q_s dB^Q_s , \q \pas,
 \eeas
 where $B^Q$ is a Brownian Motion under $Q$. By analogy with Lemma \ref{G_R_same}, it holds \pas ~ that
   \bea  \label{eqn-c080}
     \G^Q_t  =  R^{Q,0}_t , \q \fa   t \in [0, T].
      \eea
      In particular, we see that $R^{Q,0}$ is, in fact, a continuous process.

 \bs Next, we recall a comparison theorem of RBSDEs; see  Theorem 4.1 of \cite{EKPPQ-1997}. (We restate it in a more general form.)
 \begin{prop} \label{prop_comp_RBSDE}
    Let $(\G, \cZ,K )$ \(resp. $(\G',\cZ',K')$\) be a solution of RBSDE $(\xi, h, S)$ \(resp.  RBSDE $(\xi', h', S')$\) in the sense of Definition \ref{defn_RBSDE}. Additionally, assume that

     \ss    \no \hb{}\,\,\,\,(i) either $h$ or $h'$ is Lipschitz in $(y,z)$;

     \ss    \no \hb{}\,\,(ii) it holds \pas ~ that     $ \xi \le \xi'  $ and $ S_t \le S'_t$ for any $  t \in [0, T]$;

     \ss    \no (iii)  it holds \dtp ~ that      $    h(t,\o,y,z) \le h'(t,\o,y,z)$ for any  $ (y,z)  \in \hR \times \hR^d$.

  \ms  \no    Then it holds \pas~ that $      \G_t \le \G'_t$ for any  $ t \in [0, T] $.
       \end{prop}

         Since it holds \dtp~ that
  \beas
      \wt{f} (t, \o, u) \dfnn  \underset{z \in \hR^d}{\inf} \big( f(t, \o, z)+ \lan u, z \ran \big) \le  f(s , \o, \th^Q_s(\o))  + \lan u,\th^Q_s(\o)  \ran
      =  h_Q(s , \o, u)  , \q \fa u \in \hR^d.
  \eeas
    we see from Proposition \ref{prop_comp_RBSDE} and \eqref{eqn-c080} that we have
    \pas
     \bea   \label{eqn-c098}
               \wt{\G}_t \le  \G^Q_t  =        R^{Q, 0}_t, \q \fa t \in [0, T].
     \eea
 Letting $t =\nu$, taking the essential infimum of right-hand-side over $Q \in \cQ^k_\nu$, and then letting $k \to \infty$,
 we can deduce from Lemma \ref{G_R_same},  \eqref{eqn-a11}, and \eqref{eqn-a03} 
 that
   \beas
  R^{Q^{*,\nu},0}_\nu    = \wt{\G}_\nu   \le      \lmtd{k \to \infty}\, \einf{Q \in \cQ^k_\nu}  R^{Q,0}_\nu
  = \lmtd{k \to \infty}\, \einf{Q \in \cQ^k_\nu}   R^Q(\nu)= \ol{V}(\nu)= V(\nu) \le R^{Q^{*,\nu}}(\nu)= R^{Q^{*,\nu},0}_\nu, \q \pas
   \eeas
   which implies that $V(\nu) = \wt{\G}_\nu  $, \pas~ Applying Lemma \ref{G_R_same}  and \eqref{eqn-a03} once again yields that
    \bea   \label{eqn-c100}
          V(\nu) = \wt{\G}_\nu =R^{Q^*,0}_\nu =R^{Q^*}(\nu)     , \q \pas
    \eea
     where $Q^* \dfnn Q^{*,0}  \in \cQ_0$. It is clear that
     $      \, d Q^* =  d Q^{*,0}      = \sE\left(\th^{*,0} \bullet B \right)_T\, dP=\sE\left(\th^* \bullet B \right)_T\, dP\,$.\\

We are now ready to state the main result of this section.

 \begin{thm}
 \label{Thm-saddle-constr}
 {\bf Existence  of a Saddle Point:}
 The pair  $(Q^*, \t^{Q^*}(0))$ is a saddle point as in (4.1).
 \end{thm}

\section{Proofs}\label{sec:Proofs}

\subsection{Proof of the Results in Sections~\ref{sec:dcrm} and \ref{sec:mainresults}}

 \ms \no {\bf Proof of Proposition
\ref{prop_represent}: }
 \cite[Proposition 1]{Bion_2009} shows that
  \bea \label{eqn-b001}
 \rho_{\nu, \g}(\xi)\,=\,\esup{Q \in \cQ_{\nu, \g}}
 \Big( E_Q\left[-\xi \big| \cF_\nu
 \right]- \a_{\nu, \g}(Q)\Big)  , \q \pas
 \eea
Here we have set $\,\cQ_{\nu, \g}
\dfnn    \big\{Q \in \cP_\nu :    E_Q \big[ \a_{\nu, \g}(Q) \big]<\infty \big\}\,$, and the quantity $$\a_{\nu, \g}(Q) \,\dfnn \,\esup{\eta \in \mathbb{L}^\infty( \cF_\g )}
\dneg \,\,\Big( E_Q [-\eta|\cF_\nu]- \rho_{\nu, \g}(\eta)\Big)$$ is
known as the ``minimal penalty" of $\rho_{\nu, \g}\,$. (The representation \eqref{eqn-b001}
 was shown for $Q <\dneg< P$ rather than $Q \sim P$ in \cite{Bion_2009}.
  However, our assumption {\bf (A4)} assures that \eqref{eqn-b001} also holds.  For a proof, see \cite[Lemma 3.5]{Follmer_Penner_2006} and \cite[Theorem 3.1]{Kloppel_Schweizer_2007}. )

Thanks to \cite[Theorem 5(i) and the proof of
Proposition 9(v)]{DPR_2009}, there exists a nonnegative function
$f: [0,T] \times \O \times \hR^d \rightarrow  [0, \infty]$ satisfying
$(\mathfrak{f}\,1)$-$(\mathfrak{f}\,3)$, such that for each $Q \in \cQ_{\nu, \g}$ we have
 \beas
 \a_{\nu, \g}(Q)\,=\,E_Q\Big(  \int_\nu^\g
  f \neg \left(s, \th^Q_s\right) \neg ds\,\Big|\, \cF_\nu
 \Big), \q  \pas
 \eeas
 Hence we can rewrite $ \cQ_{\nu, \g}   = \left\{Q \in \cP_\nu :\, E_Q
 \int_\nu^\g f \neg \left(s, \th^Q_s\right) \neg ds  <
\infty \right\},$ and  \eqref{eqn-b001} becomes
 \bea \label{eqn-b002}
  \rho_{\nu, \g}(\xi)\,=\,\esup{Q \in \cQ_{\nu, \g}}\, E_Q\left[-\xi- \int_\nu^\g
  f \neg \left(s, \th^Q_s\right) \neg ds\,\Big|\, \cF_\nu
 \right], \q \pas
 \eea
  Since $\,\cQ_\nu \equiv \cQ_{\nu,T} \subset \cQ_{\nu, \g}\, $, it follows readily that
 \bea \label{eqn-b003}
  \underset{Q \in \cQ_\nu }{\essinf}\,E_Q \Big[Y_\g \neg
 +\int_\nu^\g f \neg \big(s, \th^Q_s\Big)   ds
\, \Big|\,\cF_\nu\big]\, \ge\, \underset{Q \in \cQ_{\nu, \g} }{\essinf}\,E_Q \Big[Y_\g \neg
 +\int_\nu^\g f   \big(s, \th^Q_s\big)   ds \,\Big|\,\cF_\nu \Big], \q \pas
 \eea
On the other hand, for any given $Q \in \cQ_{\nu, \g}$, the predictable
process $ \th^{\wt{Q}}_t \dfnn \b1_{\{t \le \g\}} \th^Q_t$, $t \in
[0,T]$ induces a probability measure $\wt{Q} \in \cP_\nu$ via
  $ \,   d \wt{Q}  \dfnn \sE \big(\th^{\wt{Q}}\bullet B \big)_T\, dP\,$.
 Since  $\,f  \big(t, \th^{\wt{Q}}_t\big)=
 \b1_{\{t \le \g\}}  f \big(t, \th^Q_t \big) $, \dtp\, from $(\mathfrak{f}\,3)$, it follows 
 \beas
   E_{\wt{Q}}  \int_\nu^T f   \big(s, \th^{\wt{Q}}_s\big)
\neg \,ds \, = \,E_{\wt{Q}}   \int_\nu^\g f   \big(s,
\th^Q_s \big) \, ds \, =\,E_Q  \int_\nu^\g f
\big(s, \th^Q_s\big) \, ds \,<\, \infty\,,
 \eeas
 thus $\wt{Q} \in \cQ_\nu\,$. Then we can deduce
 \beas
  \underset{Q \in \cQ_\nu }{\essinf}\,E_Q \left[Y_\g \neg
 +\int_\nu^\g f \neg \left(s, \th^Q_s\right) \neg ds
 \bigg|\cF_\nu\right]
 &\le& E_{\wt{Q}} \left[Y_\g \neg +\int_\nu^\g f \neg \left(s, \th^{\wt{Q}}_s\right) \neg ds
 \bigg|\cF_\nu\right]
 =E_{\wt{Q}} \left[Y_\g \neg +\int_\nu^\g f \neg \left(s, \th^{Q}_s\right) \neg ds
 \bigg|\cF_\nu\right] \\
 &=&  E_Q \left[Y_\g \neg +\int_\nu^\g f \neg \left(s, \th^{Q}_s\right) \neg ds
 \bigg|\cF_\nu\right], \q \pas
 \eeas
Taking the essential infimum of the right-hand-side over $Q \in
\cQ_{\nu, \g}$ yields
 \beas
  \underset{Q \in \cQ_\nu }{\essinf}\,E_Q \left[Y_\g \neg
 +\int_\nu^\g f \neg \left(s, \th^Q_s\right) \neg ds
 \bigg|\cF_\nu\right]
  \le  \underset{Q \in \cQ_{\nu, \g}}{\essinf}\,  E_Q \left[Y_\g \neg +\int_\nu^\g f \neg
\left(s, \th^{Q}_s\right) \neg ds
 \bigg|\cF_\nu\right], \q  \pas;
 \eeas
 this, together with \eqref{eqn-b003} and \eqref{eqn-b002},
  proves \eqref{RM_represent}. \qed

\ms \no {\bf Proof of Lemma \ref{lmt_V}:} {\bf (1)}  Since
$\left\{\cQ^k_\nu \right\}_{k \in \hN} $
 is an increasing sequence of sets contained in $\cQ_\nu$, it follows that
   \bea \label{eqn-b004}
      \einf{Q \in \cQ_\nu} E_Q \left[Y_\g \neg
  +\int_\nu^\g f \neg \left(s, \th^Q_s\right) \neg ds \bigg|\cF_\nu\right]
  \le \lmtd{k \to \infty}\, \einf{Q \in \cQ^k_\nu}
     E_Q \left[Y_\g \neg
  +\int_\nu^\g f \neg \left(s, \th^Q_s\right) \neg ds \bigg|\cF_\nu\right], \q  \pas
     \eea
Now let us fix a probability measure $\,Q \in \cQ_\nu\,$, and  define the stopping times
    \beas
     \d^Q_m \, \dfnn \, \inf\left\{ t \in [\nu,T]:
     \hb{$\int_\nu^t \big[ f \neg \left(s, \th^Q_s\right)   +  \big|\th^Q_s\big|^2 \big] ds  > m$}  \right\} \land T  , \q m \in \hN.
   \eeas
   It is easy to see that $\lmtu{m \to \infty} \d^Q_m = T$, \pas~
 For any $(m, k) \in \hN^2$,  the predictable process
   $   \,  \th^{Q^{m,k}}_t \dfnn \b1_{\{t \le \d^Q_m\}} \b1_{A^Q_{\nu,k}} \th^Q_t$, $ t \in  [0,T]  $
 induces a probability measure $Q^{m,k} \in \cQ^k_\nu$ by
  \bea \label{eqn-b005}
 d Q^{m,k} \,\, \dfnn \,\, \sE\big(\th^{Q^{m,k}}\bullet B \big)_T \cdot d P
   \eea
(recall the notation of (\ref{AQNK})).    It follows from $(\mathfrak{f}\,3)$ that
    \bea \label{eqn-b007}
    f   \big(t, \th^{Q^{m,k}}_t\big)\,=\,\b1_{\{t \le \d^Q_m\}} \b1_{A^Q_{\nu,k}} f   \big(t,    \th^Q_t  \big)\,, \q \dtp
    \eea
Then we can deduce from Bayes' Rule (see, e.g., \cite[Lemma
3.5.3]{Kara_Shr_BMSC}) that
   \bea  \label{eqn-b009}
  \qq \qq && \hspace{-2cm} \einf{Q \in \cQ^k_\nu}  E_{Q} \left[Y_\g \neg +\int_\nu^\g
   f \neg \left(s, \th^{Q}_s\right) \neg ds
  \bigg|\cF_\nu\right]  \le E_{Q^{m,k}} \left[Y_\g \neg +\int_\nu^\g f \neg \left(s, \th^{Q^{m,k}}_s\right) \neg ds
  \bigg|\cF_\nu\right] \nonumber \\
  &=& E \left[   Z^{Q^{m,k}}_{\nu, T}   \left( Y_{\g}+ \int_\nu^{\g \land \d^Q_m}
  \neg  \b1_{A^Q_{\nu,k}} f \neg \left(s,  \th^Q_s\right) \neg ds\right) \Bigg|\cF_\nu\right]
   \le  E \left[   Z^{Q^{m,k}}_{\nu, T}   \left( Y_{\g}+ \int_\nu^{\g \land \d^Q_m}
  \neg  f \neg \left(s,  \th^Q_s\right) \neg ds\right) \Bigg|\cF_\nu\right] \nonumber \\
   &=& E \left[ \Big( Z^{Q^{m,k}}_{\nu, T}-Z^Q_{\nu, \d^Q_m} \Big)  \left( Y_{\g}+ \int_\nu^{\g \land \d^Q_m}
  \neg    f \neg \left(s,  \th^Q_s\right) \neg ds\right) \Bigg|\cF_\nu\right]
  +  E \left[ \big( Z^Q_{\nu, \d^Q_m} -Z^Q_{\nu, T} \big) \cd Y_{\g} \Big|\cF_\nu\right] \nonumber\\
  && + E \left[ Z^Q_{\nu, T} Y_{\g} \Big|\cF_\nu\right]
   +  E  \left[  Z^Q_{\nu, \d^Q_m} \int_\nu^{\g \land \d^Q_m}
  \neg   f \neg \left(s,  \th^Q_s\right) \neg ds \bigg|\cF_\nu\right] \nonumber \\
   &\le & \big(
   \|Y\|_\infty +m\big)\cd E \left[ \left| Z^{Q^{m,k}}_{\nu, T}
  -Z^Q_{\nu,  \d^Q_m}  \right| \bigg|\cF_\nu\right]
  +  \|Y\|_\infty \cd E \left[ \left| Z^Q_{\nu,   \d^Q_m} -Z^Q_{\nu, T} \right|   \Big|\cF_\nu\right]
  + E_Q \left[   Y_{\g} \Big|\cF_\nu\right] \nonumber\\
  && +  E_Q  \left[   \int_\nu^{\g \land \d^Q_m}
  \neg   f \neg \left(s,  \th^Q_s\right) \neg ds \bigg|\cF_\nu\right]  \nonumber\\
   &\le & \big(
   \|Y\|_\infty +m\big)\cd E \left[ \left| Z^{Q^{m,k}}_{\nu, T}
  -Z^Q_{\nu,  \d^Q_m}  \right| \bigg|\cF_\nu\right]
  +  \|Y\|_\infty \cd E \left[ \left| Z^Q_{\nu,   \d^Q_m} -Z^Q_{\nu, T} \right|   \Big|\cF_\nu\right] \nonumber\\
  && + E_Q \left[   Y_{\g} +   \int_\nu^\g \neg
   f \neg \left(s,  \th^Q_s\right) \neg ds \bigg|\cF_\nu\right] , \qq  \pas
   \eea
From the equation   \eqref{eqn-a07} and the Dominated Convergence Theorem, we observe 
 \beas
 \lmt{k \to \infty } E \left( \int_\nu^{  \d^Q_m} \big(\b1_{A^Q_{\nu,k}}-1
\big) \th^Q_s dB_s \right)^2  \, = \, \lmt{k \to \infty } E
\int_\nu^{  \d^Q_m} \big(1-\b1_{A^Q_{\nu,k}} \big)
\big|\th^Q_s\big|^2 d s\, =\,0\,, \q \pas
 \eeas
 Thus we can find a subsequence of $\Big\{A^Q_{\nu,k}\Big\}_{k \in
 \hN}$ \Big(we still denote it by $\Big\{A^Q_{\nu,k}\Big\}_{k \in
 \hN}$\Big) such that
  \beas
  \lmt{k \to \infty} \int_\nu^{\d^Q_m} \b1_{A^Q_{\nu,k}} \th^Q_s dB_s
   = \int_\nu^{\d^Q_m} \th^Q_s dB_s \q \hb{ and } \q  \lmt{k \to \infty}
  \int_\nu^{\d^Q_m}  \b1_{A^Q_{\nu,k}}
 \big|\th^Q_s\big|^2 d s  =\int_\nu^{\d^Q_m}
  \big|\th^Q_s\big|^2 d s, \q \pas
  \eeas
and consequently, $\, \pas\,$:
 \bea
    \lmt{k \to \infty} Z^{Q^{m,k}}_{\nu, T}
    =\, \lmt{k \to
  \infty} \exp\left\{ \int_\nu^{\d^Q_m} \neg \b1_{A^Q_{\nu,k}}
\Big(\th^Q_s dB_s   -  \frac12 \big|\th^Q_s\big|^2 d s \Big)\neg
\right\} \nonumber 
  \,= \,\exp\left\{ \int_\nu^{\d^Q_m}
\neg \Big(\th^Q_s dB_s  -  \frac12 \big|\th^Q_s\big|^2 d s
\Big)\neg\right\}
\, =  \,Z^Q_{\nu, \d^Q_m}\,. 
 \label{eqn-b010}
 \eea
Since $ E\big(   Z^{Q^{m,k}}_{\nu, T} \big|\,\cF_\nu \big) =E\big(
Z^Q_{\nu, \d^Q_m} \big|\,\cF_\nu \big)=1 $, \pas~for any $k \in
\hN$, it follows from Scheff\'e's Lemma (see e.g. \cite[Section
5.10]{DW_PwM_91}) that
  \bea \label{eqn-b011}
 \lmt{k \to \infty} E \left[ \left| Z^{Q^{m,k}}_{\nu, T}
  -Z^Q_{\nu, \d^Q_m}  \right| \bigg|\cF_\nu\right]  = 0,  \q \pas
  \eea
  Hence, letting $k \to \infty$ in \eqref{eqn-b009}, we obtain
 \bea \label{eqn-b013}
 && \hspace{-2cm} \lmtd{k \to \infty}
   \einf{Q \in \cQ^k_\nu}  E_{Q} \left[Y_\g \neg +\int_\nu^\g f \neg \left(s, \th^{Q}_s\right) \neg ds
  \bigg|\cF_\nu\right] \nonumber \\
  &\le&   E_Q \left[   Y_{\g}
  +   \int_\nu^\g \neg  f \neg \left(s,  \th^Q_s\right) \neg ds \bigg|\cF_\nu\right]
   + \|Y\|_\infty \cd E \left[ \big| Z^Q_{\nu,   \d^Q_m} -Z^Q_{\nu, T} \big|   \Big|\cF_\nu\right]   , \q  \pas
 \eea
It is easy to see that $\lmtu{m \to \infty} \d^Q_m = T$, \pas\, The
right-continuity of the process $Z^Q$ then implies that $\lmt{m \to
\infty} Z^Q_{\nu, \d^Q_m} = Z^Q_{\nu, T} $, \pas~ Since $ E\left[
Z^Q_{\nu, \d^Q_m} \Big|\cF_\nu \right] =E\left[ Z^Q_{\nu, T}
\Big|\cF_\nu \right]=1$ , \pas~for any $m \in \hN$, using Scheff\'e's
Lemma once again we obtain
   \bea \label{eqn-b015}
 \lmt{m \to \infty} E \left[ \left| Z^Q_{\nu, \d^Q_m}
    -Z^Q_{\nu, T} \right|   \Big|\cF_\nu\right]  = 0,  \q \pas
  \eea
Therefore, letting $m \to \infty$ in \eqref{eqn-b013} we obtain
 \beas
  \lmtd{k \to \infty}
   \einf{Q \in \cQ^k_\nu}  E_{Q} \left[Y_\g \neg +\int_\nu^\g f \neg \left(s, \th^{Q}_s\right) \neg ds
  \bigg|\cF_\nu\right]  \le   E_Q \left[   Y_{\g}
  +   \int_\nu^\g \neg  f \neg \left(s,  \th^Q_s\right) \neg ds \bigg|\cF_\nu\right]    , \q  \pas
 \eeas
Taking the essential infimum of right-hand-side over $Q \in \cQ_\nu$ gives
 \beas
  \lmtd{k \to \infty} \einf{Q \in \cQ^k_\nu}   E_{Q} \left[Y_\g \neg +\int_\nu^\g
  f \neg \left(s, \th^{Q}_s\right) \neg ds
  \bigg|\cF_\nu\right]
   \le \einf{Q \in \cQ_\nu} E_Q \left[   Y_{\g}
  +   \int_\nu^\g \neg  f \neg \left(s,  \th^Q_s\right) \neg ds \bigg|\cF_\nu\right]   , \q
  \pas
 \eeas
 which, together with \eqref{eqn-b004}, proves \eqref{eqn-a09}.

 \ms \no {\bf (2)} By analogy with  \eqref{eqn-b004}, we have
 \bea \label{eqn-b017}
    \einf{Q \in \cQ_\nu}
 R^Q(\nu) \le \lmtd{k \to \infty}\, \einf{Q \in \cQ^k_\nu}
     R^Q(\nu), \q \pas
 \eea
 Taking the essential supremum in \eqref{eqn-b009} over $\g \in  \cS_{\nu,T}$ we get
  \bea
  \label{eqn-b019}
  \einf{Q \in \cQ^k_\nu}   R^Q(\nu) \le R^{Q^{m,k}} (\nu)  & \le&  R^Q (\nu) +\big(
   \|Y\|_\infty +m\big)\cd E \left[ \left| Z^{Q^{m,k}}_{\nu, T}
  -Z^Q_{\nu,  \d^Q_m}  \right| \bigg|\cF_\nu\right] \nonumber \\
  && +  \|Y\|_\infty \cd E \left[ \big| Z^Q_{\nu,   \d^Q_m} -Z^Q_{\nu, T} \big|   \Big|\cF_\nu\right] , \q  \pas
  \eea
 In light of \eqref{eqn-b011} and \eqref{eqn-b015}, letting $k \to \infty$
 and subsequently letting $m \to \infty$ in \eqref{eqn-b019}, we obtain
 \beas
  \lmtd{k \to \infty} \einf{Q \in \cQ^k_\nu}   R^Q(\nu) \le R^Q (\nu)   , \q  \pas
 \eeas
Taking the essential infimum of right-hand-side over $Q \in \cQ_\nu$
yields $\,
  \lmtd{k \to \infty} \einf{Q \in \cQ^k_\nu}   R^Q(\nu) \le \einf{Q \in \cQ_\nu}  R^Q (\nu)\,$, $
  \pas$ 
  which, together with \eqref{eqn-b017}, proves \eqref{eqn-a11}. \qed

\ms \no {\bf Proof of Lemma \ref{lmt_R}:} \textbf{(1)} We first show
that the family $ \left\{   E_Q  \left[ Y_\g \neg +\int_\nu^\g f
\neg \left(s, \th^{Q}_s\right) \neg ds \Big|\cF_\nu\right]
\right\}_{Q \in \cQ^k_\nu }$ is directed downwards, i.e., for any
$Q_1, Q_2 \in \cQ^k_\nu$, there exists a $Q_3 \in \cQ^k_\nu$ such
that
 \bea \label{eqn-b021}
    E_{Q_3}\neg \left[ Y_\g \neg +\int_\nu^\g f \neg \left(s, \th^{Q_3}_s\right) \neg ds \Big|\cF_\nu\right] \le
E_{Q_1}\neg \left[ Y_\g \neg +\int_\nu^\g f \neg \left(s,
\th^{Q_1}_s\right) \neg ds \Big|\cF_\nu\right] \land E_{Q_2}\neg
\left[ Y_\g \neg +\int_\nu^\g f \neg \left(s, \th^{Q_2}_s\right)
\neg ds \Big|\cF_\nu\right] \q  \pas \q
 \eea
To see this, we let $Q_1, Q_2 \in \cQ^k_\nu$ and let $A \in
\cF_\nu$. It is clear that
 \bea \label{defn_Q3}
  \th^{Q_3}_t \,\dfnn \,\b1_{\{t > \nu \}} \left(\b1_A \,\th^{Q_1}_t + \b1_{A^c}\,
 \th^{Q_2}_t\right), \q t \in [0,T]
 \eea
  forms a predictable process, thus we can define a probability measure $Q_3 \in \cM^e$
  via $  \,d Q_3   \dfnn  \sE\left(\th^{Q_3}\bullet B \right)_T  \, d P\,$.
   It follows from ($\mathfrak{f}\,3$) that
   \bea \label{eqn-b023}
    f   \big( t, \th^{Q_3}_t\big) =  \b1_{\{t > \nu \}} \left(\b1_A f \big(t, \th^{Q_1}_t\big) + \b1_{A^c} f \big(t, \th^{Q_2}_t\big)\right), \q \dtp,
   \eea
which together with \eqref{defn_Q3} implies that
   $\,\th^{Q_3} =0 $ $\dtp$ on
  $\,\[0,\nu\]\,$ and $\,\big|\th^{Q_3}_t(\o)|\vee f  \big(t, \o, \th^{Q_3}_t(\o)\big)
  = \b1_A(\o) \big|\th^{Q_1}_t(\o)|\vee f \big(t, \o, \th^{Q_1}_t(\o) \big)
  + \b1_{A^c}(\o) \big|\th^{Q_2}_t(\o)|\vee f \big(t, \o, \th^{Q_2}_t(\o) \big)  \le k$, \dtp~on $\]\nu , T\]$.
  Hence $Q_3 \in \cQ^k_\nu$. For any $\g \in \cS_{\nu, T}$, we have
   \bea
 Z^{Q_3}_{\nu,\g} &=&   \exp\left\{ \int_\nu^\g
 \big(\b1_A \th^{Q_1}_s + \b1_{A^c}
 \th^{Q_2}_s\big) dB_s- \frac12  \int_\nu^\g \big(\b1_A  |\th^{Q_1}_s |^2 + \b1_{A^c}
 |\th^{Q_2}_s |^2 \big)  ds \right\} \nonumber \\
 &=&   \exp\left\{\b1_A \Big( \int_\nu^\g
 \th^{Q_1}_s dB_s- \frac12  \int_\nu^\g \big|\th^{Q_1}_s\big|^2 d s
  \Big) + \b1_{A^c} \Big(\int_\nu^\g \th^{Q_2}_s dB_s- \frac12
 \int_\nu^\g \big|\th^{Q_2}_s\big|^2 d s \Big) \right\} \nonumber\\
 &=&  \b1_A  \exp\left\{ \int_\nu^\g
 \th^{Q_1}_s dB_s- \frac12  \int_\nu^\g \big|\th^{Q_1}_s\big|^2 d s
 \right\} + \b1_{A^c} \exp\left\{ \int_\nu^\g \th^{Q_2}_s dB_s- \frac12
 \int_\nu^\g \big|\th^{Q_2}_s\big|^2 d s \right\} \label{eqn-b025}  \\
 &=&   \b1_A   Z^{Q_1}_{\nu,\g}   + \b1_{A^c}  Z^{Q_2}_{\nu,\g} , \q \pas \nonumber
 \eea
Then Bayes' Rule implies that
 \bea  \label{eqn-b027}
 && \hspace{-2cm}  E_{Q_3}\neg \left[ Y_\g \neg +\int_\nu^\g f \neg \left(s, \th^{Q_3}_s\right) \neg ds
 \Big|\cF_\nu\right]
    =   E \left[ Z^{Q_3}_{\nu,T} \left( Y_\g \neg +\int_\nu^\g f \neg \left(s, \th^{Q_3}_s\right) \neg ds \right)
 \Big|\cF_\nu\right]  \nonumber  \\
 &  =&  E \left[ \b1_A Z^{Q_1}_{\nu,T} \left( Y_\g \neg +\int_\nu^\g f \neg \left(s, \th^{Q_1}_s\right) \neg ds
 \right)\neg +\b1_{A^c} Z^{Q_2}_{\nu,T} \left( Y_\g \neg +\int_\nu^\g f \neg \left(s, \th^{Q_2}_s\right) \neg ds \right)
 \Big|\cF_\nu\right]  \nonumber \\
 &  =&  \b1_A E_{Q_1} \neg \left[  Y_\g \neg +\int_\nu^\g f \neg \left(s, \th^{Q_1}_s\right) \neg ds
 \Big|\cF_\nu\right]
  \neg +\b1_{A^c} E_{Q_2} \neg \left[  Y_\g \neg +\int_\nu^\g f \neg \left(s, \th^{Q_2}_s\right) \neg ds
  \Big|\cF_\nu\right] , \q
  \pas
 \eea
 Letting $A = \left\{E_{Q_1}   \left[  Y_\g \neg +\int_\nu^\g f \neg \left(s, \th^{Q_1}_s\right) \neg ds
 \Big|\cF_\nu\right] \le
E_{Q_2} \left[  Y_\g \neg +\int_\nu^\g f \neg \left(s,
\th^{Q_2}_s\right) \neg ds\Big|\cF_\nu\right]\right\} \in \cF_\nu$
above, one obtains that
 \beas
 \qq  E_{Q_3}\neg \left[ Y_\g \neg +\int_\nu^\g f \neg \left(s, \th^{Q_3}_s\right) \neg ds \Big|\cF_\nu\right]
 = E_{Q_1}\neg \left[ Y_\g \neg +\int_\nu^\g f \neg \left(s, \th^{Q_1}_s\right) \neg ds
 \Big|\cF_\nu\right] \land E_{Q_2}\neg \left[ Y_\g \neg +\int_\nu^\g
 f \neg \left(s, \th^{Q_2}_s\right) \neg ds \Big|\cF_\nu\right] \q  \pas
 \eeas
proving \eqref{eqn-b021}.  Appealing to the basic properties of the
essential infimum (e.g., \cite[Proposition VI-\b1-1]{Neveu_1975}),
  we can find a sequence
 $\left\{Q^{\g,k}_n \right\}_{n \in \hN}$ in $ \cQ^k_\nu$ such that
 \eqref{eqn-a13} holds.

 \ms \no \textbf{(2)} Taking essential suprema over $\g
\in \cS_{\nu, T}$ on both sides of \eqref{eqn-b027}, we can deduce
from Lemma \ref{lem_ess} that
 \beas
  R^{Q_3}(\nu) &=& \esup{\g \in \cS_{\nu, T}}E_{Q_3}\left[ Y_\g \neg +\int_\nu^\g f \neg
  \left(s, \th^{Q_3}_s\right) \neg ds
  \Big|\cF_\nu\right]\\
  &=& \b1_A \, \esup{\g \in \cS_{\nu, T}}
   E_{Q_1}   \left[  Y_\g \neg +\int_\nu^\g f \neg \left(s, \th^{Q_1}_s\right) \neg ds\Big|\cF_\nu\right]
  +\b1_{A^c} \, \esup{\g \in \cS_{\nu, T}}
   E_{Q_2} \left[  Y_\g \neg +\int_\nu^\g f \neg \left(s, \th^{Q_2}_s\right) \neg ds\Big|\cF_\nu\right]\\
  &=&  \b1_A R^{Q_1}(\nu)  +\b1_{A^c} R^{Q_2}(\nu) , \q \pas
 \eeas
Taking $A = \big\{R^{Q_1}(\nu) \le R^{Q_2}(\nu)\big\} \in \cF_\nu$
yields that $ R^{Q_3}(\nu)= R^{Q_1}(\nu) \land  R^{Q_2}(\nu)$,
\pas, thus the family $  \{ R^Q(\nu)  \}_{Q \in \cQ^k_\nu }$ is
directed downwards. Applying Proposition VI-\b1-1 of
\cite{Neveu_1975} once again, one can find a sequence
$\{Q^{(k)}_n\}_{n \in \hN}$ in $ \cQ^k_\nu$ such that
\eqref{eqn-a15} holds.    \qed

  \ms \no {\bf Proof of Lemma \ref{lmt_tau_nu}:} Let $Q_1, Q_2 \in \cQ^k_\nu$.
  We define the stopping time $\g \dfnn \t^{Q_1}(\nu) \land \t^{Q_2}(\nu)
\in \cS_{\nu, T}$ and the event $A \dfnn \{R^{Q_1,0}_\g \le
R^{Q_2,0}_\g \} \in \cF_\g$. It is clear that
 \bea \label{defn_Q3_2}
  \th^{Q_3}_t \dfnn \b1_{\{t > \g \}} \left(\b1_A \th^{Q_1}_t + \b1_{A^c}
 \th^{Q_2}_t\right), \q t \in [0,T]
 \eea
  forms a predictable process, thus we can define a probability measure $Q_3 \in \cM^e$ by
 $  \,(  d Q_3/ d P )  \dfnn  \sE\left(\th^{Q_3}\bullet B \right)_T$.  By analogy with  \eqref{eqn-b023}, we have
   \bea \label{eqn-b029}
    f  \big(t, \th^{Q_3}_t \big) =  \b1_{\{t > \g \}} \left(\b1_A f  \big(t, \th^{Q_1}_t \big) + \b1_{A^c}
 f  \big( t, \th^{Q_2}_t \big) \right), \q \dtp
   \eea
   which together with \eqref{defn_Q3_2} implies that    $\th^{Q_3} =0, \, $ $\dtp$ on
  $\[0,\g\]$ and $\big|\th^{Q_3}_t(\o)|\vee f \big( t, \o, \th^{Q_3}_t(\o)\big)
   \le k$, \dtp~on $\]\g , T\]$.
 Hence $Q_3 \in \cQ^k_\g \subset \cQ^k_\nu$, thanks to Remark \ref{rem_belong}.
   Moreover, by analogy with \eqref{eqn-b025}, we can deduce that
  for any $\,\z \in \cS_{\g, T}\,$ we have
  \bea \label{eqn-b031}
 Z^{Q_3}_{\g, \z}\, = \, \b1_A   Z^{Q_1}_{\g, \z}   + \b1_{A^c}  Z^{Q_2}_{\g, \z}\,\,, \q \pas
 \eea
 Now fix $t \in [0,T]$. For any $\si \in \cS_{\g \vee t, T}$,
 \eqref{eqn-b031} shows that
  \beas
  Z^{Q_3}_{\g \vee t, \si}= \frac{Z^{Q_3}_{\g, \si}}{Z^{Q_3}_{\g, \g \vee
  t}}= \b1_A \frac{Z^{Q_1}_{\g, \si}}{Z^{Q_1}_{\g, \g \vee
  t}} + \b1_{A^c} \frac{Z^{Q_2}_{\g, \si}}{Z^{Q_2}_{\g, \g \vee
  t}}= \b1_A  Z^{Q_1}_{\g \vee t, \si} + \b1_{A^c} Z^{Q_2}_{\g \vee t,
  \si}, \q \pas,
  \eeas
and Bayes' Rule together with \eqref{eqn-b029} imply then
 \beas
 && \hspace{-2cm} E_{Q_3}\left[ Y_\si \neg +\int_{\g \vee t}^\si f \neg \left(s, \th^{Q_3}_s\right) \neg ds
 \bigg|\cF_{\g \vee t}\right]
 = E \left[ Z^{Q_3}_{\g \vee t, \si} \left( Y_\si \neg +\int_{\g \vee t}^\si f \neg
 \left(s, \th^{Q_3}_s\right) \neg ds \right)
 \Bigg|\cF_{\g \vee t}\right]\\
 &=& E \left[ \b1_A \,Z^{Q_1}_{\g \vee t, \si} \left( Y_\si \neg
 +\int_{\g \vee t}^\si f \neg \left(s, \th^{Q_1}_s\right) \neg ds
 \right)+ \b1_{A^c}\, Z^{Q_2}_{\g \vee t, \si} \left( Y_\si \neg
 +\int_{\g \vee t}^\si f \neg \left(s, \th^{Q_2}_s\right) \neg ds \right)
 \Bigg|\cF_{\g \vee t}\right]\\
 &=& \b1_A \,E_{Q_1}   \left[  Y_\si \neg +\int_{\g \vee t}^\si f \neg \left(s, \th^{Q_1}_s\right) \neg ds
 \bigg|\cF_{\g \vee t}\right]
  +\b1_{A^c} \,E_{Q_2} \left[  Y_\si \neg +\int_{\g \vee t}^\si f \neg \left(s, \th^{Q_2}_s\right) \neg ds
  \bigg|\cF_{\g \vee t}\right] , \q  \pas
 \eeas
Taking essential suprema over $\si \in \cS_{\g \vee t, T}$ on both
sides above, we can deduce from Lemma \ref{lem_ess} as well as
\eqref{eqn-a03} that
 \beas
 R^{Q_3,0}_{\g \vee t}  =  R^{Q_3}(\g \vee t) =  \b1_A R^{Q_1}(\g \vee t)  +\b1_{A^c} R^{Q_2}(\g \vee t) =
   \b1_A R^{Q_1,0}_{\g \vee t}  +\b1_{A^c} R^{Q_2,0}_{\g \vee t}, \q \pas
 \eeas
 Since $R^{Q_i,0}$, $i=1,2,3$ are all RCLL processes, we have $\,
  R^{Q_3,0}_{\g \vee t} = \b1_A R^{Q_1,0}_{\g \vee t}  +\b1_{A^c} R^{Q_2,0}_{\g \vee t}, \q
 \fa t \in [0,T]\,$
outside   a null set $N$, and this      implies
 \bea \label{eqn-b033}
 \t^{Q_3}(\nu) &\dneg =& \dneg \inf\left\{t \in [\nu,T] :\, R^{Q_3,0}_t=Y_t\right\}  \le
   \inf\left\{t \in  [\g, T] :\,  R^{Q_3,0}_t=Y_t \right\}  \nonumber \\
 &\dneg =& \dneg \b1_A \inf\left\{t \in [\g,T] :\,  R^{Q_1,0}_t=Y_t  \right\}
 +\b1_{A^c} \inf\left\{t \in [\g,T] :\, R^{Q_2,0}_t=Y_t \right\}, \q \pas~
 \eea
Since $R^{Q_j,0}_{\t^{Q_j}(\nu)}=Y_{\t^{Q_j}(\nu)}$, \pas~ for
$j=1,2\,$, and since $\g = \t^{Q_1}(\nu) \land \t^{Q_2}(\nu)$, it holds
\pas~that $ Y_\g$ is equal either to $R^{Q_1,0}_\g$ or to $
R^{Q_2,0}_\g$. Then the definition of the set $A$ shows that
$R^{Q_1,0}_\g =Y_\g$ holds \pas~on $A$, and that $R^{Q_2,0}_\g =Y_\g$ holds \pas~on $A^c$, both of which further imply that
 \beas
 \b1_A  \inf\left\{t \in [\g,T] :\, R^{Q_1,0}_t =Y_t\right\} =\g\b1_A   \q \hb{and}
 \q \b1_{A^c} \inf\left\{t \in [\g, T] :\, R^{Q_2,0}_t =Y_t\right\}=\g\b1_{A^c} , \q
 \pas~
 \eeas
We conclude from (\ref{eqn-b033}) that $\,\t^{Q_3}(\nu) \le   \g =
\t^{Q_1}(\nu) \land \t^{Q_2}(\nu)\,$ holds \pas, hence the
family $ \{\t^Q(\nu) \}_{Q \in \cQ^k_\nu}$ is directed downwards.
Thanks to \cite[page 121]{Neveu_1975}, we can find a sequence
$\left\{Q^{(k)}_n\right\}_{n \in \hN}$ in $\, \cQ^k_\nu\,$, such that $\,
 \t_k (\nu) = \einf{Q \in \cQ^k_\nu} \, \t^Q(\nu)
   = \lmtd{n \to \infty} \t^{Q^{(k)}_n}(\nu)\,$, $\,\pas$
The limit $\underset{n \to \infty}{\lim} \dneg \da
 \t^{Q^{(k)}_n}(\nu)$ is also a stopping time in $\cS_{\nu,T}$. \qed

\ms \no {\bf Proof of Lemma \ref{lem_paste}:}
  It is easy to see from \eqref{eqn-a17} and ($\mathfrak{ f}\,$3) that
   \bea \label{eqn-b037}
      \th^{Q'}= \th^Q =0,  ~\dtp \hb{ on } \[0,\nu\],
   \eea
   and that
   \bea \label{eqn-b039}
    f \big(t, \th^{Q'}_t\big) =  \b1_{\{t \le \g \}} f \big(t, \th^Q_t\big)
    + \b1_{\{t > \g  \}} f \big(t, \th^{\wt{Q}}_t\big) , \q \dtp
   \eea
As a result
  \beas
  E_{Q'}  \int_\nu^T f \big(s, \th^{Q'}_s\big)  ds
  &=& E_{Q'}   \int_\nu^{\g} f   \left(s, \th^Q_s\right)   ds
  +E_{Q'}  \int_{\g}^T f   \big(s, \th^{\wt{Q}}_s \big)   ds
 \\
 &\le& E_Q   \int_\nu^\g f \neg \left(s, \th^Q_s\right) \neg ds
      +E_{Q'}   \hb{$\int_{\g}^T k
  ds$}      \le  E_Q   \int_\nu^T f \neg \left(s, \th^Q_s\right) \neg ds
    +  kT  < \infty,
  \eeas
  thus $Q' \in \cQ_\nu$. If $Q \in \cQ^k_\nu$, we see from
  \eqref{eqn-a17} and \eqref{eqn-b039} that
 \beas
  \big|\th^{Q'}_t(\o)|\vee f\left(t, \o, \th^{Q'}_t(\o)\right)= \left\{
 \begin{array}{ll}
  \big|\th^{Q}_t(\o)|\vee f\left(t, \o, \th^{Q}_t(\o)\right) \le k \q & \dtp \hb{ on }\]\nu, \g\],\vspace{0.1cm} \\
   \big|\th^{\wt{Q}}_t(\o)|\vee f\left(t, \o, \th^{\wt{Q}}_t(\o)\right)
 \le k \q & \dtp \hb{ on }\] \g, T \],
 \end{array} \right.
 \eeas
 which, together with \eqref{eqn-b037}, shows that $Q' \in \cQ^k_\nu$.

 \ms Now we fix $\si \in \cS_{\g,T}$. For any $\d \in \cS_{\si,T}$,   Bayes' Rule shows
 \beas
   E_{Q'}\left[ Y_\d +  \int_\si^\d f \neg \left(s, \th^{Q'}_s\right) \neg ds \bigg|\cF_{\si}\right]
 = E_{Q'}\left[ Y_\d +  \int_\si^\d f \neg \left(s, \th^{\wt{Q}}_s\right) \neg ds \bigg|\cF_{\si}\right]
= E_{\wt{Q}}\left[ Y_\d +  \int_\si^\d f \neg \left(s,
\th^{\wt{Q}}_s\right) \neg ds\,
 \Big|\,\cF_{\si}\right], \q \pas,
 \eeas
 and  \eqref{eqn-a03} implies
 \beas
 \hspace{2cm}  R^{Q',0}_{\si}= R^{Q'} (\si )
  &= & \underset{\d \in \cS_{\si, T}}{\esssup}\,E_{Q'}\left[ Y_\d
 +  \int_\si^\d f \neg \left(s, \th^{Q'}_s\right) \neg ds \bigg|\cF_{\si}\right]  \\
 & = & \underset{\d \in \cS_{\si,T}}{\esssup}\,
  E_{\wt{Q}}\left[ Y_\d +  \int_\si^\d f \neg \left(s, \th^{\wt{Q}}_s\right) \neg ds
 \bigg|\cF_{\si}\right]
 = R^{\wt{Q}}(\si)=R^{\wt{Q},0}_{\si}, \q \pas
 \hspace{2.2cm}\hb{\qed}
\eeas

  \no {\bf Proof of Theorem \ref{V_process}: }
 Fix $Q \in \cQ_\nu$. For any $m, k \in \hN$, we consider the probability
measure $Q^{m,k} \in \cQ^k_\nu$ as defined in \eqref{eqn-b005}.
  In light of Lemma \ref{lmt_tau_nu}, for any $l \in \hN$ there exists a
 sequence $\,\big\{Q^{(l)}_n\big\}_{n \in \hN}\,$ in $\cQ^l_\nu$ such that $\,
  \t_l (\nu)   = \lmtd{n \to \infty} \t^{Q^{(l)}_n}(\nu)\,$, $    \pas$
  Now let $k,l,m,n \in \hN$ with $k \le l$. Lemma \ref{lem_paste} implies that
   the predictable process
 \beas
  \th^{Q^{m,k,l}_n}_t \, \dfnn \,\, \b1_{\{t \le \t_l(\nu) \}}  \th^{Q^{m,k}}_t  + \b1_{\{t > \t_l(\nu)
  \}} \th^{Q^{(l)}_n}_t , \q  t \in [0,T]
 \eeas
  induces a probability measure $Q^{m,k,l}_n \in \cQ^l_\nu $ via
 $   d Q^{m,k,l}_n   = \sE \big( \th^{Q^{m,k,l}_n} \bullet B \big)_T\, d P\,$,
 such that for any $t \in [0,T]$, we have $\,
   R^{Q^{m,k,l}_n,0}_{\t_l (\nu) \vee t} =R^{Q^{(l)}_n,0}_{\t_l (\nu) \vee t}\,$,  $\, \pas$
  Since $R^{Q^{m,k,l}_n,0} $ and $R^{Q^{(l)}_n,0} $ are both RCLL processes,
 outside  a null set $N$ we have
 \beas
 R^{Q^{m,k,l}_n,0}_{\t_l (\nu) \vee t} = R^{Q^{(l)}_n,0}_{\t_l (\nu) \vee t}, \q  \fa t \in [0,T]
 \eeas
and this, together with the fact that $\t_l (\nu) \le
 \t^{Q^{m,k,l}_n}(\nu) \land \t^{Q^{(l)}_n}(\nu) $, \pas~implies
 \bea  \label{eqn-b041}
 \t^{Q^{m,k,l}_n}(\nu)& \dneg  =& \dneg \inf\left\{t \in \big[\nu,T\big] :\, R^{Q^{m,k,l}_n,0}_t=Y_t\right\}
 = \inf\left\{t \in \big[\t_l (\nu),T\big] :\,
 R^{Q^{m,k,l}_n,0}_t=Y_t\right\} \nonumber \\
 & \dneg=&\dneg \inf\left\{t \in \big[\t_l (\nu),T\big] :\, R^{Q^{(l)}_n,0}_t=Y_t\right\}
 = \inf\left\{t \in \big[\nu,T\big] :\,
 R^{Q^{(l)}_n,0}_t=Y_t\right\}=\t^{Q^{(l)}_n}(\nu), \q \pas  \qq
 \eea
 Similar to \eqref{eqn-b007}, we have
  \bea \label{eqn-b043}
  f \neg \left(t,\th^{Q^{m,k,l}_n}_t\right) = \b1_{\{t \le \t_l(\nu) \}} f\left(t, \th^{Q^{m,k}}_t\right)
   + \b1_{\{t > \t_l(\nu)  \}} f\left(t,\th^{Q^{(l)}_n}_t\right),\q \dtp
  \eea
Then one can deduce from \eqref{eqn-b041} and \eqref{eqn-b043} that
 \bea   \label{eqn-b045}
 \ol{V}(\nu) & \dneg =& \underset{Q \in \cQ_\nu}{\essinf}\,
  R^Q(\nu) \le  R^{Q^{m,k,l}_n}(\nu) \,= \, \dneg E_{Q^{m,k,l}_n}\bigg( Y_{\t^{Q^{m,k,l}_n}(\nu)}
  +  \int_\nu^{\t^{Q^{m,k,l}_n}(\nu)}   f   \big(s, \th^{Q^{m,k,l}_n}_s\big)   ds
 \, \Big|\,\cF_\nu\bigg)
  \nonumber  \\
  & \dneg =& \dneg   E_{Q^{m,k,l}_n}\left[Y_{\t^{Q^{(l)}_n}(\nu)}
  +  \int_{\t_l (\nu)}^{\t^{Q^{(l)}_n}(\nu)} \neg f \neg \left(s, \th^{Q^{m,k,l}_n}_s\right) \neg ds
  \Bigg|\cF_\nu\right] +   E_{Q^{m,k}} \left[  \int_\nu^{\t_l (\nu)}
  \neg f \neg \left(s,  \th^{Q^{m,k,l}_n}_s \right) \neg ds\bigg|\cF_\nu\right] \nonumber \\
 &\dneg  =&\dneg  E  \left[\left( Z^{Q^{m,k,l}_n}_{\nu, \t^{Q^{(l)}_n}(\nu)} -Z^{Q^{m,k}}_{\nu,\t_l (\nu)}
  \right)\cd \left( Y_{\t^{Q^{(l)}_n}(\nu)}
  +  \int_{\t_l (\nu)}^{\t^{Q^{(l)}_n}(\nu)} \neg f \neg \left(s, \th^{Q^{(l)}_n}_s\right) \neg ds
  \right)\Bigg|\cF_\nu \right]  \nonumber \\
  && \dneg + \, E  \left[ Z^{Q^{m,k}}_{\nu,\t_l (\nu)}
    \left( Y_{\t^{Q^{(l)}_n}(\nu)}
  +  \int_{\t_l (\nu)}^{\t^{Q^{(l)}_n}(\nu)} \neg f \neg \left(s, \th^{Q^{(l)}_n}_s\right) \neg ds
  \right)\Bigg|\cF_\nu \right] +   E_{Q^{m,k}} \left[  \int_\nu^{\t_l (\nu)}
  \neg f \neg \left(s, \th^{Q^{m,k}}_s\right) \neg ds\bigg|\cF_\nu\right] \nonumber \\
  &\dneg  \le  &\dneg  \big( \|Y\|_\infty+ lT\big) \cd
   E\left[ \left| Z^{Q^{m,k,l}_n}_{\nu, \t^{Q^{(l)}_n}(\nu)} -Z^{Q^{m,k}}_{\nu,\t_l (\nu)}\right|
   \bigg|\cF_\nu  \right]  +E \left[Z^{Q^{m,k}}_{\nu,\t_l (\nu)}
    \left( Y_{\t^{Q^{(l)}_n}(\nu)}
  + k \big(\t^{Q^{(l)}_n}(\nu)  - \t_l (\nu)\big)\right)
  \bigg|\cF_\nu \right]    \nonumber \\
   && \dneg + E_{Q^{m,k}} \left[  \int_\nu^{\t_l (\nu)}
  \neg f \neg \left(s, \th^{Q^{m,k}}_s\right) \neg ds\bigg|\cF_\nu\right], \q  \pas
  \eea
Because $\dis E \Big(\int_{\t_l (\nu)}^{\t^{Q^{(l)}_n}(\nu)}
\th^{Q^{(l)}_n}_s dB_s \Big)^2   =E \int_{\t_l
(\nu)}^{\t^{Q^{(l)}_n}(\nu)} \big|\th^{Q^{(l)}_n}_s\big|^2 d s  \le
l^2 E\left[ \t^{Q^{(l)}_n}(\nu)- \t_l (\nu) \right]$, which goes to
zero as $ n \to \infty$, using similar arguments to those that lead
to \eqref{eqn-b010}, we can find a subsequence of
$\left\{Q^{(l)}_n\right\}_{n \in
 \hN}$ \Big(we still denote it by $\left\{Q^{(l)}_n\right\}_{n \in
 \hN}$\Big) such that
 $  \lmt{n \to \infty} Z^{Q^{m,k,l}_n}_{\nu, \t^{Q^{(l)}_n}(\nu)}
      =  Z^{Q^{m,k}}_{\nu,\t_l (\nu)} $, \pas\,
 Since $ E\left[   Z^{Q^{m,k,l}_n}_{\nu, \t^{Q^{(l)}_n}(\nu)}
\bigg|\cF_\nu \right] =E\left[   Z^{Q^{m,k}}_{\nu,\t_l (\nu)}
\Big|\cF_\nu \right]=1    $, \pas ~for any $n \in \hN$, Scheff\'e's
Lemma implies
  \bea \label{eqn-b047}
 \lmt{n \to \infty}  E\Big( \,\Big| Z^{Q^{m,k,l}_n}_{\nu,
 \t^{Q^{(l)}_n}(\nu)} \neg -Z^{Q^{m,k}}_{\nu,\t_l (\nu)}\Big|
 \,  \Big|\,\cF_\nu  \Big)\, =\, 0,  \q \pas
  \eea
  On the other hand, since
  \beas
  Z^{Q^{m,k}}_{\nu,\t_l (\nu)}
    \left| Y_{\t^{Q^{(l)}_n}(\nu)}
  + k \big(\t^{Q^{(l)}_n}(\nu)  - \t_l (\nu)\big)\right|
  \le Z^{Q^{m,k}}_{\nu,\t_l (\nu)}
    \left(  \|Y \|_\infty  + kT \right), \q \pas,
  \eeas
 and since $Y$ is right-continuous, the Dominated Convergence Theorem gives
   \bea  \label{eqn-b049}
  \lmt{n \to \infty}  E \left[Z^{Q^{m,k}}_{\nu,\t_l (\nu)}
    \left( Y_{\t^{Q^{(l)}_n}(\nu)}
  + k \big(\t^{Q^{(l)}_n}(\nu)  - \t_l (\nu)\big)\right)
  \Big|\cF_\nu \right]= E \left[Z^{Q^{m,k}}_{\nu,\t_l (\nu)}
     Y_{\t_l(\nu)}  \Big|\cF_\nu \right]= E_{Q^{m,k}} \big[
     Y_{\t_l(\nu)}  \big|\cF_\nu \big] , \q \pas \q \;\;
  \eea
  Therefore, letting $n \to \infty$ in \eqref{eqn-b045}, we can deduce from
  \eqref{eqn-b047} and \eqref{eqn-b049} that
   \beas
   \ol{V}(\nu) \le     E_{Q^{m,k}} \left[ Y_{\t_l(\nu)}+ \int_\nu^{\t_l (\nu)}
  \neg f \neg \left(s, \th^{Q^{m,k}}_s\right) \neg ds\, \bigg|\,\cF_\nu\right], \q \pas
   \eeas
As $l \to \infty$, the Bounded Convergence Theorem gives
 \beas
   \ol{V}(\nu)   \le     E_{Q^{m,k}} \left[ Y_{\t(\nu)}+ \int_\nu^{\t(\nu)}
  \neg f \neg \left(s, \th^{Q^{m,k}}_s\right) \neg
  ds\bigg|\cF_\nu\right], \q \pas
 \eeas
whence, just as in \eqref{eqn-b009}, we   deduce
   \bea \label{eqn-b051}
   \ol{V}(\nu) &\dneg \le& \dneg    E_{Q^{m,k}} \left[ Y_{\t(\nu)}+ \int_\nu^{\t(\nu)}
  \neg f \neg \left(s, \th^{Q^{m,k}}_s\right) \neg
  ds\bigg|\cF_\nu\right] \nonumber \\
  &\dneg \le & \dneg \big(
   \|Y\|_\infty +m\big)\cd E \left[ \left| Z^{Q^{m,k}}_{\nu, \t(\nu)}
  -Z^Q_{\nu, \t(\nu)\land \d^Q_m}  \right| \bigg|\cF_\nu\right]
  +  \|Y\|_\infty \cd E \left[ \big| Z^Q_{\nu, \t(\nu)\land \d^Q_m}
   -Z^Q_{\nu, \t(\nu)} \big|   \Big|\cF_\nu\right] \nonumber \\
  &&\dneg + E_Q \left[   Y_{\t(\nu)} +   \int_\nu^{\t(\nu) }
  \neg   f \neg \left(s,  \th^Q_s\right) \neg ds \bigg|\cF_\nu\right]   , \q \pas
   \eea
 By analogy with  \eqref{eqn-b011} and \eqref{eqn-b015}, one can show that for any $m \in \hN$  we have   $ \lmt{k \to \infty} E \left[ \left| Z^{Q^{m,k}}_{\nu, \t(\nu)}
  -Z^Q_{\nu, \t(\nu)\land \d^Q_m}  \right| \bigg|\cF_\nu\right]  = 0$,
  \pas~ and that  $ \lmt{m \to \infty} E \left[ \left| Z^Q_{\nu, \t(\nu)\land \d^Q_m}
    -Z^Q_{\nu, \t(\nu)} \right| \,  \Big|\cF_\nu\right]  = 0$, \pas\,
 Therefore, letting $k \to \infty$ and subsequently letting $m \to \infty$ in \eqref{eqn-b051}, we obtain
  \beas
   \ol{V}(\nu ) \le  E_Q \left[   Y_{\t(\nu)}
    + \int_\nu^{\t(\nu) }  \neg   f \neg \left(s,  \th^Q_s\right) \neg ds \bigg|\cF_\nu\right]   , \q \pas
   \eeas
Taking the essential infimum of the right-hand-side over $Q \in
\cQ_\nu$ yields
 \beas
 \ol{V}(\nu ) &\le&  \einf{Q \in \cQ_\nu} E_Q \left[   Y_{\t(\nu)}
    + \int_\nu^{\t(\nu) }  \neg   f \neg \left(s,  \th^Q_s\right) \neg ds
    \bigg|\cF_\nu\right]\\
    &\le& \esup{\g \in \cS_{\nu, T}}\,\einf{Q \in \cQ_\nu} E_Q \left[
    Y_\g   + \int_\nu^\g  \neg   f \neg \left(s,  \th^Q_s\right)
    \neg ds \bigg|\cF_\nu\right]=\ul{V}(\nu )\le \ol{V}(\nu ), \q \pas
 \eeas
 and the result follows.
\qed

\ms \no {\bf Proof of Proposition \ref{V_Y_meet}:}
 \ms For each fixed $k \in \hN$, there exists on the strength of Lemma \ref{lmt_tau_nu}  a
 sequence $\big\{Q^{(k)}_n\big\}_{n \in \hN}$ in $\cQ^k_\nu$ such that $\,
    \t_k (\nu)   = \lmtd{n \to \infty} \t^{Q^{(k)}_n}(\nu) \,$, $\,  \pas$

  For any $ n \in \hN$, the predictable process
  $ \th^{\wt{Q}^{(k)}_n}_t \dfnn \b1_{\{t > \t_k(\nu)\}} \th^{Q^{(k)}_n}_t
  $, $t \in [0,T]$ induces a probability measure
  $\wt{Q}^{(k)}_n$ by  $  \,d \wt{Q}^{(k)}_n  \dfnn  \sE\big(\wt{Q}^{(k)}_n\bullet B \big)_T\, dP
  = Z^{Q^{(k)}_n}_{\t_k(\nu),T}\, dP\,
  $. Since $\nu \le \si  \dfnn  \t(\nu) \le \t_k(\nu) \le \t^{\wt{Q}^{(k)}_n}(\nu)$, \pas, we have
   $\,\wt{Q}^{(k)}_n \in  \cQ^k_{\t_k(\nu)} \subset \cQ^k_\si \subset \cQ^k_\nu\,  $ and
 \bea \label{eqn-b055}
  \t^{\wt{Q}^{(k)}_n}(\nu)=\inf \big\{t \in [\nu,T] :\, R^{\wt{Q}^{(k)}_n,0}_t = Y_t \big\}   =\inf \big\{t \in [\si,T] :\, R^{\wt{Q}^{(k)}_n,0}_t = Y_t  \big\}=\t^{\wt{Q}^{(k)}_n}(\si), \q \pas~
 \eea
We also know from Lemma \ref{lem_paste} that for any $t \in [0,T]\,$:  $\,
  R^{\wt{Q}^{(k)}_n,0}_{\t_k (\nu) \vee t} =R^{Q^{(k)}_n,0}_{\t_k
(\nu) \vee t} \,$, $\,  \pas$

  Since $R^{\wt{Q}^{(k)}_n,0} $ and $R^{Q^{(k)}_n,0} $ are both RCLL processes,
there exists  a null set $N$ outside which we have $\,
 R^{\wt{Q}^{(k)}_n,0}_{\t_k (\nu) \vee t} = R^{Q^{(k)}_n,0}_{\t_k (\nu) \vee t}\,, ~~  \fa \, t \in
 [0,T]\,$.
By analogy  with \eqref{eqn-b041} and \eqref{eqn-b007}, respectively, we have
 \bea  \label{eqn-b057}
 \t^{\wt{Q}^{(k)}_n}(\nu) =\t^{Q^{(k)}_n}(\nu), \q \pas
 \eea
and 
 $  f  \big(t, \th^{\wt{Q}^{(k)}_n}_t\big)
 =\b1_{\{t > \t_k(\nu)\}}  f \big(t, \th^{Q^{(k)}_n}_t\big), ~~\dtp $
 Then we can deduce from \eqref{eqn-b055}, \eqref{eqn-b057} that
 \bea   \label{eqn-b059}
 V(\si) & \dneg =& \dneg \ol{V}(\si) \, =\,  \underset{Q \in \cQ_\si}{\essinf}\,
  R^Q(\si)\, \le \, R^{\wt{Q}^{(k)}_n}(\si) \,=\, E_{\wt{Q}^{(k)}_n}\Bigg(Y_{\t^{Q^{(k)}_n}(\nu)}
  +  \int_\si^{\t^{Q^{(k)}_n}(\nu)} \neg \b1_{\{s > \t_k(\nu)\}}\, f  \Big(s,  \th^{Q^{(k)}_n}_s\Big)   ds\,\bigg|\, \cF_\si\Bigg)
  \nonumber \\
   & \dneg =& \dneg   E_{\wt{Q}^{(k)}_n}\left[Y_{\t^{Q^{(k)}_n}(\nu)}
  +  \int_\si^{\t^{Q^{(k)}_n}(\nu)} \neg \b1_{\{s > \t_k(\nu)\}} \,f \neg \left(s,  \th^{Q^{(k)}_n}_s\right) \neg ds
  \Bigg|\cF_\si\right] \nonumber \\
 &\dneg  =&\dneg  E  \left[\left( Z^{\wt{Q}^{(k)}_n}_{\si, \t^{Q^{(k)}_n}(\nu)}
 -1 \right)\cd \left(  Y_{\t^{Q^{(k)}_n}(\nu)}
  +  \int_{\t_k(\nu)}^{\t^{Q^{(k)}_n}(\nu)} \neg f \neg \left(s, \th^{Q^{(k)}_n}_s\right) \neg ds
  \right)\Bigg|\cF_\si \right]  \nonumber \\
  && \dneg +E  \left[ Y_{\t^{Q^{(k)}_n}(\nu)}
  +  \int_{\t_k(\nu)}^{\t^{Q^{(k)}_n}(\nu)} \neg f \neg \left(s, \th^{Q^{(k)}_n}_s\right) \neg ds\Bigg|\cF_\si \right]
   \nonumber \\
  &\dneg  \le  &\dneg  \big( \|Y\|_\infty+ kT\big) \neg \cd \neg
   E\left[ \left| Z^{Q^{(k)}_n}_{\t_k(\nu), \t^{Q^{(k)}_n}(\nu)} -1 \right|
   \bigg|\cF_\si  \right] \neg  + E  \left[ Y_{\t^{Q^{(k)}_n}(\nu)}
  \neg +  k\big( \t^{Q^{(k)}_n}(\nu)-\t_k(\nu) \big)  \bigg|\cF_\si \right] ,
  ~~  \pas \q
  \eea
  Just as in \eqref{eqn-b047}, it can shown that $\,
 \lmt{n \to \infty} E\Big( \,\Big|  Z^{Q^{(k)}_n}_{\t_k(\nu), \t^{Q^{(k)}_n}(\nu)} -1 \Big|
\,   \Big| \,\cF_\si  \Big)\, =\, 0\,$,  $\, \pas$;
  on the other hand, the Bounded Convergence Theorem implies
   \beas
   \lmt{n \to \infty}E  \Big( Y_{\t^{Q^{(k)}_n}(\nu)}
  +  k\big( \t^{Q^{(k)}_n}(\nu)-\t_k(\nu) \big) \, \Big|\,\cF_\si
  \Big)\,=E\, \left[ Y_{\t_k(\nu)} \Big|\cF_\si \right], \q \pas
   \eeas
Letting $n \to \infty$ in \eqref{eqn-b059} yields
 $  \,V(\si) \le E \left[ Y_{\t_k(\nu)} \Big|\cF_\si \right], ~$ \pas, and applying the   Bounded Convergence Theorem  we obtain $\,
 V(\si) \le \lmt{k \to \infty} E \left[ Y_{\t_k(\nu)} \big|\cF_\si \right]
 =E \left[ Y_\si \big|\cF_\si \right]=Y_\si \,$, $  \pas$
 The reverse inequality is rather obvious. \qed

\ms \no {\bf Proof of Proposition \ref{V_sub}:}
 Fix $k \in \hN$. In light of \eqref{eqn-a15}, we can find a sequence
  $\{Q^{(k)}_n\}_{n \in \hN} \subset \cQ^k_\g$ such that
  \bea \label{eqn-b061}
   \underset{Q \in \cQ^k_\g}{\essinf}\,
  R^Q(\g) = \underset{n \to \infty}{\lim} \dneg \da R^{Q^{(k)}_n}(\g), \q \pas~
  \eea
 For any $n \in \hN$, Lemma \ref{lem_paste} implies that
   the predictable process $\,
  \th^{\wt{Q}^{(k)}_n}_t \dfnn \b1_{\{t \le \g \}}  \th^Q_t  + \b1_{\{t >
  \g  \}} \th^{Q^{(k)}_n}_t \,$, $\,  t \in [0,T]\,$
  induces a probability measure $\, \wt{Q}^{(k)}_n \in \cQ_\g $ via
 $  \,  d \wt{Q}^{(k)}_n    \dfnn \sE\big(\th^{\wt{Q}^{(k)}_n} \bullet B \big)_T\,d P\,$,
 such that for any $t \in [0,T]$, $ R^{\wt{Q}^{(k)}_n}(\g) =R^{Q^{(k)}_n}(\g)$,
 \pas~ Since $\g \le \t(\nu) \le \t^{\wt{Q}^{(k)}_n}(\nu) $, \pas,
applying \eqref{eqn-a05} yields
 \bea \label{eqn-b063}
  V(\nu)&\le & \dneg  R^{\wt{Q}^{(k)}_n}(\nu)   =  E_{\wt{Q}^{(k)}_n}\left[
R^{\wt{Q}^{(k)}_n}(\g)+ \int_\nu^\g f \neg \left(s,
\th^{\wt{Q}^{(k)}_n}_s\right) \neg ds
\bigg|\cF_\nu\right]=E_{\wt{Q}^{(k)}_n}\left[ R^{Q^{(k)}_n}(\g)+
\int_\nu^\g f \neg \left(s, \th^Q_s\right) \neg ds
\bigg|\cF_\nu\right] \nonumber  \\
&\dneg =& \dneg E_Q\left[ R^{Q^{(k)}_n}(\g)+ \int_\nu^\g f \neg
\left(s, \th^Q_s\right) \neg ds \bigg|\cF_\nu\right] , \q \pas~
 \eea
 It follows from \eqref{def_R_Q} that
  \bea \label{eqn-b065}
  - \|Y\|_\infty \le  Y_\g \le  R^{Q^{(k)}_n}(\g) \le  \|Y\|_\infty + kT, \q \pas
  \eea
Letting $n \to \infty$ in \eqref{eqn-b063}, we can deduce from the
Bounded Convergence Theorem that
  \beas
  V(\nu) \le E_Q\left[ \lmtd{n \to \infty} R^{Q^{(k)}_n}(\g) \bigg|\cF_\nu\right]+E_Q\left[
\int_\nu^\g f \neg \left(s, \th^Q_s\right) \neg ds
\bigg|\cF_\nu\right]= E_Q\left[\underset{Q \in \cQ^k_\g}{\essinf}\,
  R^Q(\g) +
\int_\nu^\g f \neg \left(s, \th^Q_s\right) \neg ds
\bigg|\cF_\nu\right], \q \pas
  \eeas
 Letting $ n \to \infty$ in \eqref{eqn-b065}, one sees from \eqref{eqn-b061} that $\,
   - \|Y\|_\infty \le  \underset{Q \in \cQ^k_\g}{\essinf}\,  R^Q(\g)  \le  \|Y\|_\infty + kT\,$ holds $   \pas$,
 and this leads to
  \beas
   - \|Y\|_\infty \le  \underset{Q \in \cQ^k_\g}{\essinf}\,  R^Q(\g)
   \le  \underset{Q \in \cQ^1_\g}{\essinf}\,  R^Q(\g)  \le  \|Y\|_\infty + T, \q \pas
  \eeas
  From the Bounded Convergence Theorem and
  Lemma \ref{lmt_V} we obtain now
 \beas
  V(\nu) \le E_Q\left[ \lmtd{k \to \infty} \, \underset{Q \in \cQ^k_\g}{\essinf}\,
  R^Q(\g) \bigg|\cF_\nu\right]+E_Q\left[
\int_\nu^\g f \neg \left(s, \th^Q_s\right) \neg ds
\bigg|\cF_\nu\right]= E_Q\left[V (\g) + \int_\nu^\g f \neg \left(s,
\th^Q_s\right) \neg ds \bigg|\cF_\nu\right], \q \pas \q \hb{\qed}
 \eeas

\ms \no {\bf Proof of Lemma \ref{lem_switch}:} Fix $k \in \hN$. For any $Q \in \cQ^k_\nu$, the predictable process  $ \th^{\wt{Q}}_t \dfnn \b1_{\{t > \nu \vee \g  \}} \th^Q_t  $, $t \in [0,T]$ induces a probability measure  $\wt{Q}$ by  $   \,(d \wt{Q}/ d P) \, \dfnn \sE\left(\wt{Q}\bullet B \right)_T = Z^{Q}_{\nu \vee \g,T}  $. Remark  \ref{rem_belong} shows that
  $\wt{Q} \in  \cQ^k_{\nu \vee \g} \subset \cQ^k_\nu \cap \cQ^k_\g
  $. By analogy with \eqref{eqn-b007}, we have  $  f \neg \left(t, \th^{\wt{Q}}_t \right)
 =\b1_{\{t > \nu \vee \g  \}}  f\left(t, \th^Q_t \right) $,  \dtp~
 Then one can deduce that
   \bea \label{eqn-b067}
 && \hspace{-1.2cm} \b1_{\{ \nu =\g \}} E_{\wt{Q}} \left[  Y_{\si \vee \g}
  + \int_\g^{\si \vee \g} f \neg \left(s,\th^{\wt{Q}}_s\right) \neg ds \bigg|\cF_{
  \g}\right]
   = \b1_{\{ \nu =\g \}} E_{\wt{Q}} \left[  Y_{\si \vee \g}
  + \int_\g^{\si \vee \g} \b1_{\{s > \nu \vee \g  \}} f \neg \left(s,\th^Q_s\right) \neg ds \bigg|\cF_{ \nu}\right]
 \nonumber \\
  &\tneg =& \tneg  E_{\wt{Q}} \left[ \b1_{\{ \nu =\g \}} \left( Y_{\si \vee \nu}
  + \int_\nu^{\si \vee \nu} f \neg \left(s,\th^{Q}_s\right) \neg ds \right) \bigg|\cF_{
  \nu}\right]=   E\left[ E_Q \left[  \b1_{\{ \nu =\g \}} \left( Y_{\si \vee \nu}
  + \int_\nu^{\si \vee \nu} f \neg \left(s,\th^{Q}_s\right) \neg ds \right) \bigg|\cF_{\nu \vee \g}\right]
  \Bigg|\cF_{ \nu}\right] \nonumber \\
  &\tneg=&\tneg  E\left[ \b1_{\{ \nu =\g \}}  E_Q \left[   Y_{\si \vee \nu}
  +\neg \int_\nu^{\si \vee \nu} \neg f \neg \left(s,\th^{Q}_s\right) \neg ds   \bigg|\cF_{\nu }\right]
  \Bigg|\cF_{ \nu}\right] \neg = \b1_{\{ \nu =\g \}}  E_Q \neg \left[   Y_{\si \vee \nu}
  + \neg \int_\nu^{\si \vee \nu}\neg f \neg \left(s,\th^{Q}_s\right) \neg ds   \bigg|\cF_{\nu
  }\right], ~~ \pas,
   \eea
   which implies
    \beas
    \b1_{\{ \nu =\g \}}  E_Q \left[   Y_{\si \vee \nu}
  + \int_\nu^{\si \vee \nu} f \neg \left(s,\th^{Q}_s\right) \neg ds   \bigg|\cF_{\nu
  }\right] \ge \b1_{\{ \nu =\g \}}
  \einf{Q \in \cQ^k_\g}  E_Q \left[  Y_{\si \vee \g}
  + \int_\g^{\si \vee \g} f \neg \left(s,\th^Q_s\right) \neg ds   \bigg|\cF_{  \g}\right],
  \q \pas
    \eeas
Taking the essential infimum of the left-hand-side over $Q \in
\cQ^k_\nu$, one can deduce from Lemma \ref{lem_ess} that
 \beas
 \b1_{\{ \nu =\g \}} \einf{Q \in \cQ^k_\nu}  E_Q \left[   Y_{\si \vee \nu}
  + \int_\nu^{\si \vee \nu} f \neg \left(s,\th^{Q}_s\right) \neg ds   \bigg|\cF_{\nu
  }\right] & \neg=& \neg \einf{Q \in \cQ^k_\nu} \b1_{\{ \nu =\g \}} E_Q \left[   Y_{\si \vee \nu}
  + \int_\nu^{\si \vee \nu} f \neg \left(s,\th^{Q}_s\right) \neg ds   \bigg|\cF_{\nu
  }\right] \\
  & \neg \ge & \neg \b1_{\{ \nu =\g \}}
  \einf{Q \in \cQ^k_\g}  E_Q \left[  Y_{\si \vee \g}
  + \int_\g^{\si \vee \g} f \neg \left(s,\th^Q_s\right) \neg ds   \bigg|\cF_{  \g}\right] , \q \pas
   \eeas
 Letting $k \to \infty$, we see from Lemma \ref{lmt_V} (1) that
  \beas
  \b1_{\{ \nu =\g \}} \einf{Q \in \cQ_\nu}  E_Q \left[   Y_{\si \vee \nu}
  + \int_\nu^{\si \vee \nu} f \neg \left(s,\th^{Q}_s\right) \neg ds   \bigg|\cF_{\nu
  }\right] \ge \b1_{\{ \nu =\g \}}
  \einf{Q \in \cQ_\g}  E_Q \left[  Y_{\si \vee \g}
  + \int_\g^{\si \vee \g} f \neg \left(s,\th^Q_s\right) \neg ds   \bigg|\cF_{  \g}\right]  , \q \pas
   \eeas
Reversing the roles of $\nu$ and $\g$, we obtain \eqref{eqn-a23}.

 \ms On the other hand, taking essential supremum over $\si \in \cS_{0,T}$ on both sides of
 \eqref{eqn-b067}, we can deduce from Lemma \ref{lem_ess} that
  \beas
   \b1_{\{ \nu =\g \}} R^{\wt{Q}}(\g)
   &=&   \esup{\si \in \cS_{0,T}}
  \b1_{\{ \nu =\g \}} E_{\wt{Q}} \left[  Y_{\si \vee \g}
  + \int_\g^{\si \vee \g} f \neg \left(s,\th^{\wt{Q}}_s\right) \neg ds   \bigg|\cF_{
  \g}\right] \\
  &=&  \esup{\si \in \cS_{0,T}}
  \b1_{\{ \nu =\g \}}  E_Q \left[   Y_{\si \vee \nu}
  + \int_\nu^{\si \vee \nu} f \neg \left(s,\th^{Q}_s\right) \neg ds   \bigg|\cF_{\nu
  }\right] = \b1_{\{ \nu =\g \}} R^{Q}(\nu), \q \pas
  \eeas
  which implies that
   $  \b1_{\{ \nu =\g \}} R^{Q}(\nu) \ge   \b1_{\{ \nu =\g \}} \einf{Q \in \cQ^k_\g}
  R^{Q}(\g)$ , \pas~ Taking the essential infimum of the left-hand-side over $Q \in
\cQ^k_\nu$, one can deduce from Lemma \ref{lem_ess} that
 \beas
   \b1_{\{ \nu =\g \}} \einf{Q \in
\cQ^k_\nu} R^{Q}(\nu) = \einf{Q \in \cQ^k_\nu} \b1_{\{ \nu =\g \}}
R^{Q}(\nu)  \ge \b1_{\{ \nu =\g \}} \einf{Q \in \cQ^k_\g}
  R^{Q}(\g) , \q  \pas
 \eeas
 Letting $k \to \infty$, we see from Lemma \ref{lmt_V} (2) that
  \beas
  \b1_{\{ \nu =\g \}} V(\nu) = \b1_{\{ \nu =\g \}} \einf{Q \in
\cQ_\nu} R^{Q}(\nu)    \ge \b1_{\{ \nu =\g \}} \einf{Q \in \cQ_\g}
  R^{Q}(\g) = \b1_{\{ \nu =\g \}} V(\g) , \q \pas
   \eeas
 Reversing the roles of $\nu$ and $\g$, we obtain \eqref{eqn-a25}. \qed

  \ms \no {\bf Proof of Theorem \ref{V_RC}: Proof of (1).}

\noindent {\bf Step 1:} For any $\si, \nu \in \cS_{0,T}$, we define
 \beas
  \P^\si(\nu)   \dfnn     \b1_{\left\{\si \le \nu \right\}}
  Y_{\si  }\neg + \neg \b1_{\left\{\si > \nu \right\}}\,\underset{Q \in \cQ_{\nu}}{\essinf}\,
  E_Q \left[  Y_{\si \vee \nu }
  \neg+ \neg \int_{ \nu}^{\si \vee \nu } f \neg \left(s, \th^Q_s\right) \neg ds \bigg|\cF_{  \nu}\right].
  \eeas
 We see from \eqref{eqn-a09} that
 \bea \label{eqn-b069}
 \underset{Q \in \cQ_\nu}{\essinf}\, E_Q \left[Y_{\si \vee \nu}
  +\int_{  \nu}^{\si \vee \nu} f \neg \left(s, \th^Q_s\right) \neg ds \bigg|\cF_{  \nu}\right]
  = \lmtd{k \to \infty}  \underset{Q \in \cQ^k_\nu}{\essinf}\, E_Q \left[Y_{\si \vee \nu}
  +\int_{  \nu}^{\si \vee \nu} f \neg \left(s, \th^Q_s\right) \neg ds \bigg|\cF_{  \nu}\right],\q  \pas~
 \eea
 Fix $k \in \hN$. In light of \eqref{eqn-a13}, we can find a sequence
$\left\{ Q^{(k)}_n  \right\}_{n \in \hN}$ in $ \cQ^k_\nu$ such that
 \bea  \label{eqn-b071}
 \underset{Q \in \cQ^k_\nu}{\essinf}\, E_Q \left[Y_{\si \vee \nu}
  +\int_{  \nu}^{\si \vee \nu} f \neg \left(s, \th^Q_s\right) \neg ds \bigg|\cF_{  \nu}\right]
  = \lmtd{n \to \infty}   E_{Q^{(k)}_n} \left[Y_{\si \vee \nu}
  +\int_{ \nu}^{\si \vee \nu} f \neg \left(s, \th^{Q^{(k)}_n}_s\right) \neg ds \bigg|\cF_{\nu}\right],\q  \pas~
 \eea
 By analogy with \eqref{eqn-b065}, we have
  \bea \label{eqn-b073}
   - \|Y\|_\infty \le E_{Q^{(k)}_n} \left[  Y_{\si \vee \nu}
  \neg+ \neg \int_{ \nu}^{\si \vee \nu} f \neg \left(s, \th^{Q^{(k)}_n}_s\right) \neg ds
   \bigg|\cF_{ \nu}\right]  \le  \|Y\|_\infty + kT
  \eea
\pas; letting $n \to \infty$, we see from \eqref{eqn-b071} that
  \beas
   - \|Y\|_\infty \le  \einf{Q \in \cQ^k_\nu}  E_Q \left[  Y_{\si \vee \nu}
  \neg+ \neg \int_{ \nu}^{\si \vee \nu}  f \neg \left(s, \th^Q_s\right) \neg ds
   \bigg|\cF_{ \nu}\right]  \le  \|Y\|_\infty + kT, \q
  \pas
  \eeas
Therefore,
  \bea \label{eqn-b075}
   - \|Y\|_\infty &\le&  \einf{Q \in \cQ^k_\nu}  E_Q \left[  Y_{\si \vee \nu}
  \neg+ \neg \int_{ \nu}^{\si \vee \nu}  f \neg \left(s, \th^Q_s\right) \neg ds   \bigg|\cF_{
  \nu}\right] \nonumber \\
 & \le & \einf{Q \in \cQ^1_\nu}  E_Q \left[  Y_{\si \vee \nu}
  \neg+ \neg \int_{ \nu}^{\si \vee \nu}  f \neg \left(s, \th^Q_s\right) \neg ds   \bigg|\cF_{ \nu}\right]
   \le  \|Y\|_\infty + T, \q
  \pas
  \eea
  Letting $k \to \infty$, we see from \eqref{eqn-b069} that
  \beas
   - \|Y\|_\infty \le  \einf{Q \in \cQ_\nu}  E_Q \left[  Y_{\si \vee \nu}
  \neg+ \neg \int_{ \nu}^{\si \vee \nu}  f \neg \left(s, \th^Q_s\right) \neg ds
   \bigg|\cF_{ \nu}\right]  \le  \|Y\|_\infty + T, \q \pas
  \eeas
 which implies that
  \bea \label{eqn-b077}
  - \|Y\|_\infty \le  \P^\si(\nu)   \le  \|Y\|_\infty + T, \q \pas
  \eea
Let $\g \in \cS_{0,T}$. It follows from \eqref{eqn-a23} that
  \bea \label{eqn-b079}
   \b1_{\{ \nu =\g\}}\P^\si(\nu) & \tneg =& \tneg \b1_{\{\si \le \nu=\g \}} Y_{\si  } \neg + \b1_{\{\si > \nu =\g \}}
  \einf{Q \in \cQ_\nu}  E_Q \left[  Y_{\si \vee \nu}
  + \neg \int_\nu^{\si \vee \nu} \neg f \neg \left(s,\th^Q_s\right) \neg ds \bigg|\cF_{ \nu}\right] \nonumber \\
  & \tneg=&\tneg \b1_{\{\si \le \g = \nu  \}} Y_{\si  } \neg + \b1_{\{\si > \g = \nu \}}
  \einf{Q \in \cQ_\g}  E_Q \left[  Y_{\si \vee \g}
  + \neg \int_\g^{\si \vee \g} \neg f \neg \left(s,\th^Q_s\right) \neg ds   \bigg|\cF_{  \g}\right]
  \neg = \b1_{\{ \nu =\g\}}\P^\si(\g), ~\;\; \pas \qq \;\;
  \eea

 \ss \no {\bf Step 2:} Fix $\si \in \cS_{0,T}$.  For any $\z \in \cS_{0,T}$
 , $\nu  \in \cS_{\z,T}$ and $k \in \hN$, we let $\left\{ Q^{(k)}_n  \right\}_{n \in \hN} \subset \cQ^k_\nu$
 be the sequence described in \eqref{eqn-b071}. Then we can deduce that
 \bea  \label{eqn-b081}
    \P^\si(\z) &\le& \b1_{\left\{\si \le \z \right\}}
  Y_{\si  }\neg + \neg \b1_{\left\{\si > \z \right\}}\, E_{Q^{(k)}_n} \left[  Y_{\si \vee \z}
  \neg+ \neg \int_{ \z}^{\si \vee \z} f \neg \left(s, \th^{Q^{(k)}_n}_s\right) \neg ds
  \bigg|\cF_{ \z}\right] \nonumber \\
  &=& \b1_{\left\{\si \le \z \right\}}
  Y_{\si \land \z  }\neg + \neg \b1_{\left\{\si > \z \right\}}\,
  E  \left[ E_{Q^{(k)}_n} \left[  Y_{\si \vee \z}
  \neg+ \neg \int_{ \z}^{\si \vee \z} f \neg \left(s, \th^{Q^{(k)}_n}_s\right) \neg ds
  \bigg|\cF_{ \nu}\right] \Bigg|\cF_{
  \z}\right] \nonumber \\
  &=&  E  \left[ \b1_{\left\{\si \le \z \right\}}
  Y_{\si \land \z  }\neg + \neg \b1_{\left\{\si > \z \right\}}\, E_{Q^{(k)}_n} \left[  Y_{\si \vee \z}
  \neg+ \neg \int_{ \z}^{\si \vee \z} f \neg \left(s, \th^{Q^{(k)}_n}_s\right) \neg ds
  \bigg|\cF_{ \nu}\right] \Bigg|\cF_{
  \z}\right], \q \pas
  \eea
On the other hand, it holds \pas~ that
  \beas
   && \hspace{-2cm} \b1_{\left\{ \si > \nu \right\}}\, E_{Q^{(k)}_n} \left[  Y_{\si \vee \z}
  \neg+ \neg \int_{ \z}^{\si \vee \z} f \neg \left(s, \th^{Q^{(k)}_n}_s\right) \neg ds \bigg|\cF_{ \nu}\right]
  =   E_{Q^{(k)}_n} \left[ \b1_{\left\{\si > \nu \right\}} \left( Y_{\si  }
  \neg+ \neg \int_{ \z}^{\si  } f \neg \left(s, \th^{Q^{(k)}_n}_s\right) \neg ds \right) \bigg|\cF_{ \nu}\right]\\
  &=& E_{Q^{(k)}_n} \left[ \b1_{\left\{\si > \nu \right\}} \left( Y_{\si \vee \nu}
  \neg+ \neg \int_{ \nu}^{\si \vee \nu} f \neg \left(s, \th^{Q^{(k)}_n}_s\right) \neg ds \right) \bigg|\cF_{ \nu}\right]
  =\b1_{\left\{\si > \nu \right\}} \, E_{Q^{(k)}_n} \left[  Y_{\si \vee \nu}
  \neg+ \neg \int_{ \nu}^{\si \vee \nu} f \neg \left(s, \th^{Q^{(k)}_n}_s\right) \neg ds   \bigg|\cF_{
  \nu}\right]
  \eeas
  and that
 \beas
   && \hspace{-2cm}  \b1_{\left\{ \z < \si \le \nu \right\}}\, E_{Q^{(k)}_n} \left[  Y_{\si \vee \z}
  \neg+ \neg \int_{ \z}^{\si \vee \z} f \neg \left(s, \th^{Q^{(k)}_n}_s\right) \neg ds \bigg|\cF_{ \nu}\right]
  =   E_{Q^{(k)}_n} \left[ \b1_{\left\{ \z < \si \le \nu \right\}} \left( Y_{\si }
  \neg+ \neg \int_{ \z}^{\si} f \neg \left(s, \th^{Q^{(k)}_n}_s\right) \neg ds \right) \bigg|\cF_{ \nu}\right]\\
  &=& E_{Q^{(k)}_n} \left[ \b1_{\left\{ \z < \si \le \nu \right\}}  Y_{\si \land \nu}
    \bigg|\cF_{ \nu}\right]= \b1_{\left\{ \z < \si \le \nu \right\}}  Y_{\si \land
    \nu}= \b1_{\left\{ \z < \si \le \nu \right\}}  Y_{\si }\,;
  \eeas
 recall the definitions of the classes $\, \mathcal{P}_\nu\,$,  $\, \mathcal{Q}_\nu\,$ from subsection 1.1. Therefore, \eqref{eqn-b081} reduces to
   \beas
  \P^\si(\z)  \le E  \left[ \b1_{\left\{\si \le \nu \right\}}
  Y_{\si  }\neg + \neg  \b1_{\left\{\si > \nu \right\}} \, E_{Q^{(k)}_n} \left[  Y_{\si \vee \nu}
  \neg+ \neg \int_{ \nu}^{\si \vee \nu} f \neg \left(s, \th^{Q^{(k)}_n}_s\right) \neg ds
   \bigg|\cF_{ \nu}\right]  \Bigg|\cF_{
  \z}\right], \q \pas
   \eeas
 We obtain then   from \eqref{eqn-b071}, \eqref{eqn-b073}  and the Bounded
 Convergence Theorem, that
  \beas
  \P^\si(\z)  & \dneg \le& \dneg \lmtd{n \to \infty} E  \left[ \b1_{\left\{\si \le \nu \right\}}
  Y_{\si  }\neg + \neg  \b1_{\left\{\si > \nu \right\}}    E_{Q^{(k)}_n} \left[  Y_{\si \vee \nu}
  \neg+ \neg \int_{ \nu}^{\si \vee \nu} f \neg \left(s, \th^{Q^{(k)}_n}_s\right) \neg ds
    \bigg|\cF_{ \nu}\right]  \Bigg|\cF_{
  \z}\right]\\
 & \dneg =& \dneg E  \left[ \b1_{\left\{\si \le \nu \right\}}
  Y_{\si  }\neg + \neg  \b1_{\left\{\si > \nu \right\}} \,\einf{Q \in \cQ^k_\nu}  E_Q \left[  Y_{\si \vee \nu}
  \neg+ \neg \int_{ \nu}^{\si \vee \nu} f \neg \left(s, \th^{Q^{(k)}_n}_s\right) \neg ds
   \bigg|\cF_{ \nu}\right]  \Bigg|\cF_{
  \z}\right] , \q \pas
  \eeas
 On the other hand, we can deduce from \eqref{eqn-b069}, \eqref{eqn-b075} and the Bounded
Convergence Theorem once again that
  \bea \label{eqn-b083}
  \P^\si(\z)  & \dneg \le& \dneg \lmtd{k \to \infty} E  \left[ \b1_{\left\{\si \le \nu \right\}}
  Y_{\si  }\neg + \neg  \b1_{\left\{\si > \nu \right\}}
   \einf{Q \in \cQ^k_\nu}  E_Q \left[  Y_{\si \vee \nu}
  \neg+ \neg \int_{ \nu}^{\si \vee \nu} \neg f \neg \left(s, \th^{Q^{(k)}_n}_s\right) \neg ds
    \bigg|\cF_{ \nu}\right]  \Bigg|\cF_{
  \z}\right] \nonumber \\
  & \dneg =& \dneg E  \left[ \b1_{\left\{\si \le \nu \right\}}
  Y_{\si  }\neg + \neg  \b1_{\left\{\si > \nu \right\}} \, \einf{Q \in \cQ_\nu}  E_Q \left[  Y_{\si \vee \nu}
  \neg+ \neg \int_{ \nu}^{\si \vee \nu} \neg f \neg \left(s, \th^Q_s\right) \neg ds
   \bigg|\cF_{ \nu}\right]  \Bigg|\cF_{
  \z}\right] =E \left[ \P^\si(\nu) \big|\cF_{ \z}\right], ~ \;\;
  \pas,\qq \;
  \eea
 which implies that $\{\P^\si(t)\}_{t \in [0,T]}$ is a
 submartingale. Therefore \cite[Proposition 1.3.14]{Kara_Shr_BMSC} shows that
  \bea \label{eqn-b085}
 P\left( \hb{the limit } \P^{\si,+}_t \dfnn \lmt{n \to \infty}
   \P^\si\big(q_n(t)\big)    \hb{ exists for any }  t \in [0,T] \right)=1
 \eea
\big(where $q_n(t) \dfnn \frac{\lceil 2^nt \rceil}{2^n}\land
T$\big), and that $ \P^{\si,+} $ is an RCLL process.

 \ms \no {\bf Step 3:} For any $\nu \in \cS_{0, T}$ and $n \in \hN$,
 $ q_n(\nu)$ takes values
 in a finite set $\cD^n_T \dfnn \big([0,T)\cap \{k2^{-n}\}_{k \in \hZ}
\big) \cup \{T\}  $. Given an $\l \in \cD^n_T$, it holds for any $m
\ge n$ that $q_m (\l)=\l $ since $ \cD^n_T \subset \cD^m_T$. It
follows from \eqref{eqn-b085} that
 \beas
 \P^{\si,+}_\l
 = \underset{m \to \infty}{\lim} \P^\si\big(q_m(\l)\big)=\P^\si ( \l
 ) , \q \pas
 \eeas
 Then one can deduce from \eqref{eqn-b079} that
  \beas
  \qq \P^{\si,+}_{q_n(\nu)} =   \sum_{\l \in \cD^n_T} \b1_{\{q_n(\nu) = \l
  \}}\P^{\si,+}_{\l} = \sum_{\l \in \cD^n_T} \b1_{\{q_n(\nu) = \l
  \}}\P^\si(\l)
   =  \sum_{\l \in \cD^n_T} \b1_{\{q_n(\nu) = \l  \}}\P^\si\big(q_n(\nu)\big)
  = \P^\si\big(q_n(\nu)\big), \q \pas
  \eeas
  Thus the right-continuity of the process $\P^{\si,+}$ implies that
   \bea \label{eqn-b087}
   \P^{\si,+}_{\nu}= \lmt{n \to \infty} \P^{\si,+}_{q_n(\nu) }= \lmt{n \to \infty}
   \P^\si\big(q_n(\nu)\big), \q \pas
   \eea
  Hence \eqref{eqn-b083}, \eqref{eqn-b077} and the Bounded Convergence
Theorem imply
   \bea \label{eqn-b089}
    \P^\si(\nu) \le  \underset{n \to \infty }{\lim}
    E\big[    \P^\si(q_n(\nu)) \big|\cF_\nu\big]
   =E\big[   \P^{\si,+}_\nu \big|\cF_\nu \big]
   =    \P^{\si,+}_\nu , \q  \pas
 \eea
In the last equality we used the fact that $ \P^{\si,+}_\nu = \underset{n
 \to \infty}{\lim} \P^\si\big( q_n(\nu) \big)
 \in \cF_\nu $, thanks to the right-continuity of the Brownian filtration $\bF$.

 \ms \no {\bf Step 4:} Set $\nu, \g \in \cS_{0,T}$ and $\,
 \z \dfnn \t(\nu) \land \g, \q \z_n \dfnn \t(\nu) \land q_n(\g)\,$, $\, \fa n \in \hN\,$.
Now, let $\si \in \cS_{\z,T}$. Since $\lmtu{n \to \infty}
\b1_{\{\t(\nu) > q_n (\g)\}}
   = \b1_{\{\t(\nu) >  \g\}} $ and
    \beas
    \{\t(\nu)   > \g  \} \subset
 \left\{ q_n (\g) = q_n \big( \t(\nu) \land \g \big)
 \right\},\q \{\t(\nu) > q_n (\g)\} \subset \left\{ q_n(\g) = \t(\nu) \land
   q_n(\g) \right\} , \q \fa n \in \hN,
 \eeas
 one can deduce
 from \eqref{eqn-b089}, \eqref{eqn-b087}, and \eqref{eqn-b079}  that
 \bea \label{eqn-b091}
    \b1_{\{\t(\nu)   > \g\}} \P^\si(\z) &\le & \b1_{\{\t(\nu)   > \g\}} \P^{\si,+}_\z
    =\b1_{\{\t(\nu)   > \g\}}   \lmt{n \to \infty}  \P^\si \big( q_n (\z) \big)
    = \lmt{n \to \infty} \b1_{\{\t(\nu)   > \g\}}
    \P^\si \left(q_n \big( \t(\nu) \land \g \big) \right) \nonumber\\
    &=& \lmt{n \to \infty} \b1_{\{\t(\nu)   > \g\}}  \P^\si \big(
    q_n(\g)\big) =  \lmt{n \to \infty}\b1_{\{\t(\nu) > q_n (\g)\}}
   \P^\si \big(  q_n(\g)\big) \nonumber\\
   &=&  \lmt{n \to \infty}\b1_{\{\t(\nu) > q_n (\g)\}} \P^\si \big( \t(\nu) \land
   q_n(\g) \big) = \b1_{\{\t(\nu)   > \g\}} \lmt{n \to \infty} \P^\si \big(
     \z_n\big)
    , \q  \pas
 \eea
For any $n \in \hN$, we see from \eqref{eqn-a19} and Lemma
\ref{lem_ess} that
 \beas
 V(\z_n)&  =& \tneg \ul{V}(\z_n) = \underset{\beta \in \cS_{\z_n, T}}{\esssup}\,
 \left(\, \underset{Q \in \cQ_{\z_n}}{\essinf}\,E_Q
 \left[Y_\beta\neg +\neg \int_{\z_n}^\beta f \neg \left(s, \th^Q_s\right) \neg ds
 \bigg|\cF_{\z_n}\right]\right)   \\
 &\ge&  \underset{Q \in \cQ_{\z_n}}{\essinf}\, E_Q
 \left[Y_{\si \vee \z_n}\neg +\neg \int_{\z_n}^{\si \vee \z_n} f \neg \left(s, \th^Q_s\right) \neg ds
 \bigg|\cF_{\z_n}\right] \\
 &=&  \underset{Q \in \cQ_{\z_n}}{\essinf}\,  E_Q  \left[\b1_{\left\{\si \le \z_n \right\}} Y_{\z_n}
 \dneg +\neg  \b1_{\left\{\si > \z_n \right\}}
 \left( Y_{\si \vee \z_n}\neg +\neg \int_{\z_n}^{ \si \vee \z_n} f \neg \left(s, \th^Q_s\right) \neg ds
 \right)  \bigg|\cF_{\z_n}\right] \\
 &=& \underset{Q \in \cQ_{\z_n}}{\essinf}\, \left( \b1_{\left\{\si \le \z_n \right\}} Y_{\z_n}
 \dneg +\neg  \b1_{\left\{\si > \z_n \right\}} E_Q  \left[
   Y_{\si \vee \z_n}\neg +\neg \int_{\z_n}^{ \si \vee \z_n} f \neg \left(s, \th^Q_s\right) \neg ds
    \bigg|\cF_{\z_n}\right] \right) \\
 &=&  \b1_{\left\{\si \le \z_n \right\}} Y_{\z_n}
 \dneg +\neg  \b1_{\left\{\si > \z_n  \right\}} \underset{Q \in \cQ_{\z_n}}{\essinf}\, E_Q  \left[
   Y_{\si \vee \z_n}\neg +\neg \int_{\z_n}^{ \si \vee \z_n} f \neg \left(s, \th^Q_s\right) \neg ds
    \bigg|\cF_{\z_n}\right]  , \q  \pas
  \eeas
 Since $ \{ \t(\nu) \le \g\} \subset \{ \z_n=\z=\t(\nu)\}
 $ and $\{\si > \z_n\} \subset \{\si > \z\} $, it follows from \eqref{eqn-a23} that
 \beas
 V(\z_n)&\tneg \ge & \tneg       \b1_{\left\{\si \le \z_n \right\}}
  Y_{\z_n}\dneg + \neg \b1_{\left\{\si> \z_n, \t(\nu) > \g \right\}}\underset{Q \in
 \cQ_{\z_n}}{\essinf}\, E_Q  \left[
   Y_{\si \vee \z_n}\neg +\neg \int_{\z_n}^{ \si \vee \z_n} f \neg \left(s, \th^Q_s\right) \neg ds
    \bigg|\cF_{\z_n}\right]  \nonumber \\
 &&  +   \b1_{\left\{\si> \z_n, \t(\nu) \le \g \right\}}\underset{Q \in
 \cQ_{\z}}{\essinf}\, E_Q  \left[
   Y_{\si \vee \z}\neg +\neg \int_{\z}^{ \si \vee \z } f \neg \left(s, \th^Q_s\right) \neg ds
    \bigg|\cF_{\z}\right]   \nonumber \\
 &\tneg =& \tneg   \b1_{\left\{\si \le \z_n \right\}}  Y_{\z_n}
  \neg + \neg \b1_{\left\{\si> \z_n, \t(\nu) > \g \right\}}
    \P^\si(\z_n)  +   \b1_{\left\{\si> \z_n, \t(\nu) \le \g \right\}}    \P^\si(\z)   , \q  \pas
 \eeas
As $n \to \infty$,  the right-continuity of processes $Y$,
\eqref{eqn-b091} as well as Lemma \ref{lem_ess} show that
 \beas
   \underset{n \to \infty}{\liminf}V(\z_n)
    &\ge& \b1_{\left\{\si=\z\right\}}  Y_{\z}
  \neg + \neg \b1_{\left\{\si> \z, \t(\nu) > \g \right\}}
    \lmt{n \to \infty} \P^\si(\z_n)
  +   \b1_{\left\{\si> \z, \t(\nu) \le \g \right\}}   \P^\si(\z)   \\
    &\ge& \b1_{\left\{\si=\z\right\}}  Y_\z
  +\b1_{\left\{\si> \z\right\}}    \P^\si(\z)
    =   \b1_{\left\{\si = \z   \right\}}Y_\z + \b1_{\left\{\si > \z \right\}}
    \underset{Q \in \cQ_{\z}}{\essinf}\, E_Q \left[  Y_{\si \vee \z }
  \neg+ \neg \int_{ \z}^{\si \vee \z } f \neg \left(s, \th^Q_s\right) \neg ds \bigg|\cF_{ \z}\right] \\
  &=&  \underset{Q \in \cQ_{\z}}{\essinf}\,\left(\b1_{\left\{\si = \z   \right\}}Y_\z+\b1_{\left\{\si> \z
  \right\}}E_Q \left[  Y_{\si \vee \z }
  \neg+ \neg \int_{ \z}^{\si \vee \z } f \neg \left(s, \th^Q_s\right) \neg ds \bigg|\cF_{ \z}\right]
  \right)\\
  &=& \underset{Q \in \cQ_{\z}}{\essinf}\, E_Q \left[ \b1_{\left\{\si = \z   \right\}}Y_\z+\b1_{\left\{\si> \z
  \right\}} \left( Y_{\si  }
  \neg+ \neg \int_{ \z}^{\si } f \neg \left(s, \th^Q_s\right) \neg ds \right) \bigg|\cF_{ \z}\right] \\
   &=& \underset{Q \in \cQ_{\z}}{\essinf}\,
  E_Q \left[   Y_{\si  }
  \neg+ \neg \int_{ \z}^{\si } f \neg \left(s, \th^Q_s\right) \neg ds   \bigg|\cF_{ \z}\right], \q \pas
  \eeas
 Taking the essential supremum of the right-hand-side over $\si \in
\cS_{\z, T}$, we obtain
  \bea  \label{eqn-b093}
 \underset{n \to \infty}{\liminf}V(\z_n) \ge \underset{\si \in \cS_{\z, T}}{\esssup}\, \left(
 \underset{Q \in \cQ_\z}{\essinf}\,E_Q\left[Y_\si+\int_\z^\si f \neg \left(s, \th^Q_s\right) \neg ds
 \bigg|\cF_\z\right] \right)
 = \ul{V}(\z)= V(\z) , \q \pas~
  \eea

\ms  Let us show the reverse inequality. Fix $Q \in \cQ_\z$ and $n
\in \hN$. For any $k, m \in \hN$, the predictable process
  \beas
  \th^{Q^{m,k}_n}_t  \, \dfnn \, \b1_{\{\z_n < t \le \d^{Q,n}_m\}} \b1_{A^Q_{\z,k}} \th^Q_t, \q t \in [0,T]
  \eeas
 induces a probability measure $Q^{m,k}_n \in \cQ^k_{\z_n}$ by
 $ \,  d Q^{m,k}_n   \dfnn  \sE\left(\th^{Q^{m,k}_n} \bullet B \right)_T\, dP \,$,
 where $\d^{Q,n}_m$ is defined by
 \beas
  \d^{Q,n}_m \, \dfnn \, \inf\left\{ t \in [\z_n,T]:
  \hb{$\int_{\z_n}^t f \neg \left(s, \th^Q_s\right) \neg ds > m$}  \right\} \land T , \q m \in \hN.
 \eeas
 For any $\beta \in \cS_{\z_n, T}$, using  arguments similar to those that lead to
\eqref{eqn-b009}, we obtain
 \beas
 \qq \qq && \hspace{-2.5cm} E_{Q^{m,k}_n} \left[ Y_\beta + \int_{\z_n}^\beta f \neg
 \left(s, \th^{Q^{m,k}_n}_s\right) \neg ds \Big|\cF_{\z_n}\right]
   \le    \big( \|Y\|_\infty+ m \big) \cd  E \left[ \left|  Z^{Q^{m,k}_n}_{\z_n,T}
  -Z^Q_{\z_n,\d^{Q,n}_m}  \right|  \bigg|\cF_{\z_n}\right] \\
 && + \|Y\|_\infty \cd E \left[   \left| Z^Q_{\z_n, \d^{Q,n}_m}-Z^Q_{\z_n,T}
  \right| \bigg|\cF_{\z_n}\right] + E_Q \left[ Y_\beta + \int_{\z_n}^{\beta  }  f \neg \left(s, \th^Q_s\right) \neg ds
   \bigg|\cF_{\z_n}\right], \q \pas
 \eeas
    Then taking the essential supremum of both sides over $\beta \in \cS_{\z_n, T}$ yields that
    \bea \label{eqn-b099}
  \einf{Q \in \cQ^k_{\z_n}} R^Q (\z_n)  & \le &
     R^{Q^{m,k}_n}(\z_n)  \le    \big( \|Y\|_\infty+ m \big)  E \left[ \left|  Z^{Q^{m,k}_n}_{\z_n,T}
  -Z^Q_{\z_n,\d^{Q,n}_m}  \right|  \bigg|\cF_{\z_n}\right] \nonumber \\
  &&+ \|Y\|_\infty \cd E \left[   \left| Z^Q_{\z_n, \d^{Q,n}_m}-Z^Q_{\z_n,T}
  \right| \bigg|\cF_{\z_n}\right] +R^Q (\z_n)   , \q \pas \qq \qq
    \eea
  Just as in  \eqref{eqn-b011}, we can show that
  \beas
 \lmt{k \to \infty} E \left[  \left| Z^{Q^{m,k}_n}_{\z_n,T}
  -Z^Q_{\z_n,\d^{Q,n}_m}  \right|
  \bigg|\cF_{\z_n} \right]  = 0,  \q \pas
  \eeas
Therefore, letting $k \to \infty$ in \eqref{eqn-b099}, we know from
Lemma \ref{lmt_V} (2) that
 \bea \label{eqn-b101}
 V(\z_n) =  \lmtd{k \to \infty} \einf{Q \in \cQ^k_{\z_n}} R^Q (\z_n)    \le
    \|Y\|_\infty \cd E \left[   \left| Z^Q_{\z_n, \d^{Q,n}_m}-Z^Q_{\z_n,T}
  \right| \bigg|\cF_{\z_n}\right]   +R^Q (\z_n)   , \q \pas
    \eea
 Next, by analogy with \eqref{eqn-b015}, we have $\,
 \lmt{m \to \infty} E \Big( \Big| Z^Q_{\z_n, \d^{Q,n}_m}-Z^Q_{\z_n,T}
  \Big| \, \big| \, \cF_{\z_n}\Big) = 0,  ~ \pas
 $ 
Letting $m \to \infty$ in \eqref{eqn-b101}, we obtain  $\,
   V(\z_n)    \le  R^Q (\z_n)= R^{Q,0}_{ \z_n }, ~  \pas
$ from
\eqref{eqn-a03}. Then the right-continuity of the process $R^{Q,0}$, as well as (\ref{eqn-a03}), imply that
  \beas
      \underset{n \to \infty}{\limsup} V(\z_n) \le
      \underset{n \to \infty}{\lim}R^{Q,0}_{\z_n}=R^{Q,0}_\z=R^Q(\z), \q \pas~
  \eeas
    Taking the essential infimum of $R^Q(\z)$ over $Q \in \cQ_\z$ yields $\,
     \underset{n \to \infty}{\limsup} V(\z_n) \le  \underset{Q \in \cQ_\z}{\essinf}\,
     R^Q(\z) = \ol{V}(\z) =V(\z)\,$, $\,\pas$
    This inequality, together with \eqref{eqn-b093}, shows that
 \bea \label{eqn-b103}
 \underset{n \to \infty}{\lim} V \big( \t(\nu) \land q_n(\g) \big)  = V \big(\t(\nu) \land \g \big), \q \pas
 \eea

 \ms \no {\bf Step 5:} \ms  Now fix $\nu \in \cS_{0,T}$. It is clear that $ P
\in \cQ_\nu$ and that $\th^P_\cd \equiv 0$.  For any $t \in [0,T]$,
 \eqref{eqn-a25} implies that
 \beas
  \b1_{\{t \ge \nu\}}  V\big(\t(\nu) \land t\big) =
    \b1_{\{t \ge \nu\}} V \big(\t(\nu) \land (t \vee \nu) \big), \q
    \pas,
 \eeas
  since $\{t \ge \nu\} \subset \{\t(\nu) \land t = \t(\nu) \land (t \vee
  \nu)\}$. Then we can deduce from \eqref{eqn-a21},
  ($\mathfrak{f}\,3$), and \eqref{eqn-a19} that for any $s \in [0,t)$
 \beas
  \b1_{\{s \ge \nu\}}  V\big(\t(\nu) \land s\big) &\tneg=& \tneg
    \b1_{\{s \ge \nu\}} V \big(\t(\nu) \land (s \vee \nu) \big) \\
     &\tneg \le & \tneg  \b1_{\{s \ge \nu\}} E \left[V\big(\t(\nu) \land (t \vee \nu)\big)\neg
     + \neg \int_{\t(\nu) \land (s \vee \nu)}^{\t(\nu) \land (t \vee \nu)} \neg
f \neg \left(r, \th^P_r\right)dr \bigg|\cF_{\t(\nu) \land (s \vee
\nu)}\right]\\
&\tneg=& \tneg \b1_{\{s \ge \nu\}} E \left[ V\big(\t(\nu) \land (t
\vee \nu)\big) \bigg|\cF_{\t(\nu) \land s}\right] = E \left[
\b1_{\{s \ge \nu\}} V\big(\t(\nu) \land (t \vee \nu)\big)
\bigg|\cF_{\t(\nu)
\land s}\right]\\
 &\tneg \le&\tneg  E \left[\b1_{\{t \ge \nu\}} V\big(\t(\nu)
\land (t \vee \nu)\big)+ \b1_{\{t \ge \nu > s \}} \|Y\|_\infty
\bigg|\cF_{\t(\nu) \land s}\right]\\
 &\tneg =&\tneg E \left[ E \left[\b1_{\{t \ge \nu\}} \left( V\big(\t(\nu) \land t\big)
 +   \|Y\|_\infty \right) \bigg|\cF_{\t(\nu) }\right] \Bigg|\cF_s\right]
 - \b1_{\{s \ge \nu   \}} \|Y\|_\infty\\
 &\tneg =&\tneg  E \left[ \b1_{\{t \ge \nu\}}  \left( V\big(\t(\nu) \land t\big)
 +   \|Y\|_\infty \right)  \Bigg|\cF_s\right]- \b1_{\{s \ge \nu   \}} \|Y\|_\infty  , \q \pas,
 \eeas
 which shows that $ \Big\{ \b1_{\{t \ge \nu\}} \left( V\big(\t(\nu) \land t\big)
 +   \|Y\|_\infty \right)\Big\}_{t \in [0,T]}
    $ is a submartingale.  Hence it follows from \cite[Proposition 1.3.14]{Kara_Shr_BMSC} that
  \beas
 P\left( \hb{the limit } V^{0,\nu}_t \dfnn \lmt{n \to \infty}
  \b1_{\{q_n(t) \ge \nu\}} V\big(\t(\nu) \land q_n(t)\big)     \hb{ exists for any }  t \in [0,T] \right)=1,
 \eeas
  and that $ V^{0,\nu} $ is an RCLL process.

  Let $\z \in \cS^\star_{0,T}$ take values
 in a finite set $\{t_1 < \cds< t_m\}$. For any $\l \in \{1 \cds m\}$ and $n \in \hN$,
 since  $  \{\z=t_\l  \} \subset \{\t(\nu) \land
q_n(\z)=\t(\nu) \land q_n(t_\l)  \}$,
 one can deduce from \eqref{eqn-a25} that
 \beas
 \b1_{\{\z=t_\l  \}} V\left(\t(\nu) \land q_n(\z)\right) = \b1_{\{\z=t_\l \}}
V\left(\t(\nu) \land q_n(t_\l)\right), \q \pas
 \eeas
As $n \to \infty$,  \eqref{eqn-b103} shows
  \beas
   \b1_{ \{\z=t_\l  \}}V^{0,\nu}_\z &=& \b1_{ \{\z=t_\l   \}}V^{0,\nu}_{t_\l}
   =\b1_{\{t_\l  \ge \nu\}} \underset{n \to \infty}{\lim} \b1_{\{\z=t_\l\}} V\left(\t(\nu) \land
   q_n(t_\l)\right)\\
 &=& \b1_{\{t_\l  \ge \nu\}} \underset{n \to \infty}{\lim} \b1_{\{\z=t_\l\}}
 V\left(\t(\nu) \land q_n(\z)\right)
   = \b1_{\{\z  \ge \nu\}} \b1_{\{\z=t_\l\}}   V  \left( \t(\nu) \land \z \right)  , \q \pas
  \eeas
 Summing the above expression over $\l$, we obtain
  $  V^{0,\nu}_\z =\b1_{\{\z  \ge \nu\}} V  \left( \t(\nu) \land \z \right)$, \pas~
 Then for any $\g \in \cS_{0,T}$, the right-continuity of the
process $V^{0,\nu}$ and \eqref{eqn-b103} imply
 \beas
  V^{0,\nu}_\g= \lmt{n \to \infty}V^{0,\nu}_{q_n(\g)}
  = \lmt{n \to \infty}  \b1_{\{q_n(\g)  \ge \nu\}} V \left( \t(\nu) \land q_n(\g) \right)
   =  \b1_{\{ \g  \ge \nu\}}  V \left( \t(\nu) \land \g \right) , \q  \pas,
 \eeas
 proving (\ref{V0}). In particular, $ V^{0,\nu} $ is an RCLL modification of
  the process $\left\{ \b1_{\{ t  \ge \nu\}}  V \left( \t(\nu) \land t \right) \right\}_{t \in [0,T]} $.

\ms \no {\bf Proof of Theorem \ref{V_RC}: Proof of (2).}  Proposition
\ref{V_Y_meet} and \eqref{V0} imply that $ V^{0,\nu}_{\t(\nu)}=   V\big(\t(\nu)
\big) = Y_{\t(\nu) }$, \pas~ Hence we can deduce from the
right-continuity of processes $V^{0,\nu}$ and $Y$ that
$\t_V\neg(\nu)$ in (\ref{eqn-a27}) is a stopping time belonging to
$\cS_{\nu, \t(\nu)}$ and that
 \beas
  Y_{\t_V\neg(\nu)}= V^{0,\nu}_{\t_V\neg(\nu)}   = V(\t_V\neg(\nu)), \q
  \pas,
 \eeas
where the second equality is due to \eqref{V0}. Then it follows from
\eqref{eqn-a21} that for any $Q \in \cQ_\nu$
 \beas
 V(\nu)   \le  E_Q \left[V(\t_V\neg(\nu))+ \int_\nu^{\t_V\neg(\nu)}
 f \neg \left(s,\th^Q_s\right) \neg ds \bigg|\cF_\nu\right]
  = E_Q \left[Y_{\t_V\neg(\nu)} + \int_\nu^{\t_V\neg(\nu)}
 f \neg \left(s,\th^Q_s\right) \neg ds \bigg|\cF_\nu\right], \q \pas
 \eeas
 Taking the essential infimum of the right-hand-side over $Q \in \cQ_\nu$ yields that
 \beas
 \q  V(\nu) &\le&  \einf{Q \in \cQ_\nu}   E_Q \neg \left[Y_{\t_V\neg(\nu)} + \neg \int_\nu^{\t_V\neg(\nu)}
\neg f \neg \left(s,\th^Q_s\right) \neg ds \bigg|\cF_\nu\right]  \\
   &\le&  \esup{\g \in \cS_{\nu,T}}  \left( \einf{Q \in
  \cQ_\nu}  E_Q \neg \left[Y_{\g} + \neg \int_\nu^{\g}
 f \neg \left(s,\th^Q_s\right) \neg ds \bigg|\cF_\nu\right]
  \right) \neg =\ul{V}(\nu)=V(\nu), \q \pas,
 \eeas
 from which the claim follows. \qed

\subsection{Proof of the Results in Section~\ref{sec:saddle}}\label{sec:pf2}
 {\bf Proof of Theorem \ref{Thm-saddle}:}  It is easy to see from (i) that
  \bea   \label{eqn-c010}
   Y_{\si_*} = V(\si_*) = R^{Q^*}(\si_*)  , \q  \pas
   \eea
      which together with  (ii) and (iii) shows that for any $Q  \in \cQ_0$
  \beas
   E_{Q^*}\left[Y^{Q^*}_{\si_*}\right]
   =   E_{Q^*}\big[V^{Q^*}(\si_*)\big] =V^{Q^*}(0)= V(0)
   \le E_Q \left[V^Q( \si_*) \right]=E_Q \left[Y^Q_{\si_*} \right].
 \eeas
Thus the second inequality in  \eqref{defn_saddle}
 holds for  $(Q^*, \si_*) $. Now we show that $(Q^*, \si_*)$ satisfies the first
inequality in \eqref{defn_saddle}  in three steps:

 \ss \no $\bullet$ When $ \,\nu \in \cS_{0,\si_*}\,$,   property (iii) and \eqref{eqn-c010}
 imply that
 \bea \label{eqn-k52}
  Y^{Q^*}_\nu  \le  V^{Q^*}(\nu) =E_{Q^*} \neg \left[ V^{Q^*}(\si_*)\big|\cF_\nu\right]
  =E_{Q^*} \neg \left[Y^{Q^*}_{\si_*}\big|\cF_\nu\right],  \q \pas
 \eea
Taking the expectation $E_{Q^*} $ on both sides yields that $E_{Q^*}\big[Y^{Q^*}_\nu \big]
  \le     E_{Q^*}\big[Y^{Q^*}_{\si_*}\big]$.

 \ss \no $\bullet$ When $\, \nu \in \cS_{\si_*,T}\,$, it follows from \eqref{eqn-c010} that
 \beas
  E_{Q^*} \neg \left[Y^{Q^*}_\nu \right]
     &=& E_{Q^*} \neg \left[ E_{Q^*} \neg \left[Y_\nu +\int_{\si_*}^\nu  f\big(s, \th^{Q^*}_s\big) ds   \Big|\cF_{\si_*}
   \right]+\int_0^{\si_*}   f\big(s, \th^{Q^*}_s\big) ds    \right]  \\
  & \le&  E_{Q^*} \neg \left[ R^{Q^*} (\si_*)+\int_0^{\si_*}   f\big(s, \th^{Q^*}_s\big) ds    \right]
 = E_{Q^*} \neg \left[ Y^{Q^*}_{\si_*} \right].
 \eeas

\no $\bullet$ For a general stopping time $\nu \in \cS_{0,T}$, let us
define $\nu_1= \nu \land \si_* \in \cS_{0,\si_*}$ and $\nu_2= \nu
\vee \si_* \in \cS_{ \si_*,T}$. Since $\left\{\nu \le  \si_*\right\}
\in \cF_{\nu \land \si_*}=\cF_{\nu_1}$,  one can deduce from
(\ref{eqn-k52}) that
 \beas
  \qq  E_{Q^*} \neg\big[Y^{Q^*}_\nu \big]
  &=&E_{Q^*} \neg\left[E_{Q^*} \neg \left[ \b1_{\left\{\nu \le \si_*\right\}}Y^{Q^*}_{\nu_1}
  \neg +\b1_{\left\{\nu > \si_*\right\}}Y^{Q^*}_{\nu_2} \big|\cF_{\si_*}\right]\right]
  \neg =E_{Q^*} \neg \left[ \b1_{\left\{\nu \le \si_*\right\}}Y^{Q^*}_{\nu_1}
  \neg +\b1_{\left\{\nu > \si_*\right\}}E_{Q^*} \neg \left[Y^{Q^*}_{\nu_2} \big|\cF_{\si_*}\right]\right]\\
  &\le & E_{Q^*}  \neg\left[ \b1_{\left\{\nu \le \si_*\right\}}Y^{Q^*}_{\nu_1}
   \neg +\b1_{\left\{\nu > \si_*\right\}} \neg \left(R^{Q^*}(\si_*)   + \int_0^{\si_*}  f(s, \th^{Q^*}_s)ds \right) \right] \\
  &=& E_{Q^*}  \neg\left[\b1_{\left\{\nu \le  \si_*\right\}}Y^{Q^*}_{\nu_1} \neg +\b1_{\left\{\nu >
  \si_*\right\}}Y^{Q^*}_{\si_*}\right]=   E_{Q^*}  \neg\left[ \b1_{\left\{\nu \le  \si_*\right\}}Y^{Q^*}_{\nu_1} \neg+\b1_{\left\{\nu >  \si_*\right\}}
   E_{Q^*}  \neg\left[ Y^{Q^*}_{\si_*}\big|\cF_{\nu_1}\right]\right] \\
   & \le &   E_{Q^*}  \neg \left[ \b1_{\left\{\nu \le  \si_*\right\}}
   E_{Q^*} \neg\left[Y^{Q^*}_{\si_*}\big|\cF_{\nu_1}\right]+\b1_{\left\{\nu >  \si_*\right\}}
   E_{Q^*}  \neg\left[ Y^{Q^*}_{\si_*}\big|\cF_{\nu_1}\right]\right] = E_{Q^*}  \neg \big[Y^{Q^*}_{\si_*}\big]. \hspace{3.5cm}  \hb{\qed}
 \eeas

 \ms \no {\bf Proof of Lemma \ref{G_R_same}:} Fix $t \in [0, T]$. For any $\g  \in \cS_{\nu \vee t,T}$, we see from \eqref{eqn-c034} that
  \beas
      \wt{\G}_{\nu \vee t}
     =  \wt{\G}_\g       + \int_{\nu \vee t}^\g  f(s,    \th^{*,\nu}_s )   ds  +   \wt{K}_\g  - \wt{K}_{\nu \vee t} - \int_{\nu \vee t}^\g  \wt{\cZ}_s dB^{Q^{*,\nu}}_s, \q \pas
  \eeas
   Applying $E_{Q^{*,\nu}}\left[\,\cd\,|\cF_{\nu \vee t}\right]$ to both sides, we obtain
  \bea
   \wt{\G}_{\nu \vee t}             &=&  E_{Q^{*,\nu}}\left[\wt{\G}_\g   + \int_{\nu \vee t}^\g    f(s, \th^{*,\nu}_s)   ds +   \wt{K}_\g  - \wt{K}_{\nu \vee t} \,\Big|\,\cF_{\nu \vee t} \right]  \label{eqn-c050}\\
                   &\ge&     E_{Q^{*,\nu}}\left[Y_\g   + \int_{\nu \vee t}^\g    f(s, \th^{*,\nu}_s)   ds  \,\Big|\,\cF_{\nu \vee t} \right], \q \pas   \label{eqn-c054}
  \eea
  Let $\si^*_{\nu \vee t} \dfnn \inf\{s \in \left[\nu \vee t,T\right]: \wt{\G}_s = Y_s\}  \in \cS_{\nu \vee t,T}$. The flat-off condition satisfied by    $(\wt{\G},\wt{\cZ},\wt{K})$, and the continuity of $\wt{K}$, imply that
      \beas
    \q   0= \neg  \int_{ [\nu \vee t, \si^*_{\nu \vee t})}  \neg   \b1_{\{ \wt{\G}_s   >  Y_s  \}}   d\wt{K}_s
      =\neg  \int_{  [\nu \vee t, \si^*_{\nu \vee t} )} \neg   d\wt{K}_s      =   \underset{s \nearrow \si^*_{\nu \vee t}}{\lim}\wt{K}_s   -  \wt{K}_{\nu \vee t}
      =\wt{K}_{\si^*_{\nu \vee t}}  -  \wt{K}_{\nu \vee t}  , \q \pas
      \eeas
      Hence, taking $\g  =\si^*_{\nu \vee t} $ in \eqref{eqn-c050}, we obtain the $\pas $ property $\,
            \wt{\G}_{\nu \vee t}     =      E_{Q^{*,\nu}}\Big[ Y_{\si^*_{\nu \vee t}}
             + \int_{\nu \vee t}^{\si^*_{\nu \vee t}}   f(s, \th^{*,\nu}_s)   ds \,  \Big|\,\cF_{\nu \vee t} \Big] \,$   
              which, together with \eqref{eqn-c054} and \eqref{eqn-a03}, shows that
       \beas
         \wt{\G}_{\nu \vee t}  =\esup{\g  \in \cS_{\nu \vee t,T} } E_{Q^{*,\nu}}\Big[Y_\g   + \int_{\nu \vee t}^\g    f(s, \th^{*,\nu}_s)   ds  \,\Big|\,\cF_{\nu \vee t} \Big] = R^{Q^{*,\nu}}(\nu \vee t)= R^{{Q^{*,\nu}},0}_{\nu \vee t}, \q \pas
       \eeas
       Then the right-continuity of the processes $\wt{\G}$ and $R^{{Q^{*,\nu}},0}$ implies \eqref{eqn-c060}.  \qed

       \ms \no {\bf Proof of Theorem \ref{Thm-saddle-constr} :}  We shall show that $(Q^*, \t^{Q^*}(0))$ satisfies conditions (i)-(iii) of Theorem \ref{Thm-saddle}:

   \ss  \no 1)  It follows easily from Proposition \ref{Prop3.1} that  $Y_{\t^{Q^*}(0)} = R^{Q^*,0}_{\t^{Q^*}(0)}= R^{Q^*}\big(\t^{Q^*}(0)\big)$, \pas

  \ss  \no 2)  For any $k \in \hN$ and $Q \in \cQ^k_0$,  we can deduce from \eqref{eqn-c100}, the right-continuity of processes $ R^{Q^*,0} $ and $\wt{\G}$, as well as \eqref{eqn-c098}  that \pas
   \beas
   R^{Q^*,0}_t = \wt{\G}_t     \le       R^{Q, 0}_t   , \q    \fa t \in [0,T] .
   \eeas
   In particular, we have $  Y_{\t^Q(0)} \le  R^{Q^*,0}_{\t^Q(0)} =  R^{Q, 0}_{\t^Q(0)}=Y_{\t^Q(0)}  $, \pas ~
   Hence $  Y_{\t^Q(0)} =  R^{Q^*,0}_{\t^Q(0)} $, \pas, which implies further that $\t^{Q^*} (0) \le  \t^Q(0) $, \pas~Taking the essential infimum of right-hand-side over $Q \in \cQ^k_0$ and letting $\,k \to \infty\,$, we deduce that, in the notation of (\ref{defn_tau_nu}), we have $\,
     \t^{Q^*} (0) \le   \lmtd{k \to \infty}\, \einf{Q \in \cQ^k_0} \t^Q(0) =\t(0)\,$, $\, \pas$
     Then  \eqref{eqn-a21} shows     $   \,   V(0) \le E_Q\left[V^Q\big(\t^{Q^*} (0)\big)   \right]\,$   for any $\,Q \in \cQ_0\,$.

\smallskip
  \ss  \no 3)  For any $\nu \in \cS_{0, \t^{Q^*}(0)}$, and  since  $  \nu \le \t^{Q^*}(0) \le \t^{Q^*}(\nu)\,$ holds \pas,  one can deduce from \eqref{eqn-c100} and \eqref{eqn-a05}  that
 \beas
  \qq \q  V^{Q^*}(\nu)   &=&   R^{Q^* }(\nu)  + \int_0^\nu f(s, \th^*_s) ds
  =E_{Q^*}\left[R^{Q^*} \big( \t^{Q^*}(0) \big)  + \int_\nu^{\t^{Q^*}(0)}   f(s, \th^*_s) ds\Big| \cF_\nu\right]+ \int_0^\nu f(s, \th^*_s) ds \\
  &=& E_{Q^*}\left[R^{Q^*} \big( \t^{Q^*}(0) \big)  + \int_0^{\t^{Q^*}(0)}   f(s, \th^*_s) ds\Big| \cF_\nu\right]
     =E_{Q^*}\left[V^{Q^*} \big( \t^{Q^*}(0) \big)  \Big| \cF_\nu\right] ,  \q  \pas  \hspace{1.3cm} \hb{ \qed}
  \eeas

\medskip

\bibliographystyle{abbrvnat}
\bibliography{OS_CRM}

\begin{thebibliography}{22}
\providecommand{\natexlab}[1]{#1}
\providecommand{\url}[1]{\texttt{#1}}
\expandafter\ifx\csname urlstyle\endcsname\relax
  \providecommand{\doi}[1]{doi: #1}\else
  \providecommand{\doi}{doi: \begingroup \urlstyle{rm}\Url}\fi

\bibitem[Bayraktar and Yao(2009)]{OSNE}
E.~Bayraktar and S.~Yao.
\newblock Optimal stopping for nonlinear expectations.
\newblock Technical report, University of Michigan, 2009.
\newblock Available at http://arxiv.org/abs/0905.3601.

\bibitem[Bene{\v{s}}(1970)]{Benes_1970}
V.~E. Bene{\v{s}}.
\newblock Existence of optimal strategies based on specified information, for a
  class of stochastic decision problems.
\newblock \emph{SIAM J. Control}, 8:\penalty0 179--188, 1970.
\newblock ISSN 0363-0129.

\bibitem[Bion-Nadal(2009)]{Bion_2009}
J.~Bion-Nadal.
\newblock Time consistent dynamic risk processes.
\newblock \emph{Stochastic Process. Appl.}, 119\penalty0 (2):\penalty0
  633--654, 2009.
\newblock ISSN 0304-4149.

\bibitem[Cheridito et~al.(2006)Cheridito, Delbaen, and Kupper]{CDK-2006}
P.~Cheridito, F.~Delbaen, and M.~Kupper.
\newblock Dynamic monetary risk measures for bounded discrete-time processes.
\newblock \emph{Electron. J. Probab.}, 11:\penalty0 no. 3, 57--106
  (electronic), 2006.
\newblock ISSN 1083-6489.

\bibitem[Delbaen(2006)]{Delbaen_2006}
F.~Delbaen.
\newblock The structure of m-stable sets and in particular of the set of risk
  neutral measures.
\newblock In \emph{In memoriam {P}aul-{A}ndr\'e {M}eyer: {S}\'eminaire de
  {P}robabilit\'es {XXXIX}}, volume 1874 of \emph{Lecture Notes in Math.},
  pages 215--258. Springer, Berlin, 2006.

\bibitem[Delbaen et~al.(2009)Delbaen, Peng, and Rosazza-Gianin]{DPR_2009}
F.~Delbaen, S.~Peng, and E.~Rosazza-Gianin.
\newblock Representation of the penalty term of dynamic concave utilities.
\newblock Technical report, ETH, 2009.
\newblock Available at http://arxiv.org/abs/0802.1121.

\bibitem[El~Karoui(1981)]{El_Karoui_1981}
N.~El~Karoui.
\newblock Les aspects probabilistes du contr\^ole stochastique.
\newblock In \emph{Ninth {S}aint {F}lour {P}robability {S}ummer {S}chool---1979
  ({S}aint {F}lour, 1979)}, volume 876 of \emph{Lecture Notes in Math.}, pages
  73--238. Springer, Berlin, 1981.

\bibitem[El~Karoui et~al.(1997)El~Karoui, Kapoudjian, Pardoux, Peng, and
  Quenez]{EKPPQ-1997}
N.~El~Karoui, C.~Kapoudjian, E.~Pardoux, S.~Peng, and M.~C. Quenez.
\newblock Reflected solutions of backward {SDE}'s, and related obstacle
  problems for {PDE}'s.
\newblock \emph{Ann. Probab.}, 25\penalty0 (2):\penalty0 702--737, 1997.
\newblock ISSN 0091-1798.

\bibitem[Elliott(1982)]{Elliott_1982}
R.~J. Elliott.
\newblock \emph{Stochastic calculus and applications}, volume~18 of
  \emph{Applications of Mathematics (New York)}.
\newblock Springer-Verlag, New York, 1982.
\newblock ISBN 0-387-90763-7.

\bibitem[F{\"o}llmer and Penner(2006)]{Follmer_Penner_2006}
H.~F{\"o}llmer and I.~Penner.
\newblock Convex risk measures and the dynamics of their penalty functions.
\newblock \emph{Statist. Decisions}, 24\penalty0 (1):\penalty0 61--96, 2006.
\newblock ISSN 0721-2631.

\bibitem[F{\"o}llmer and Schied(2004)]{Follmer_Schied_2004}
H.~F{\"o}llmer and A.~Schied.
\newblock \emph{Stochastic finance}, volume~27 of \emph{de Gruyter Studies in
  Mathematics}.
\newblock Walter de Gruyter \& Co., Berlin, extended edition, 2004.
\newblock ISBN 3-11-018346-3.
\newblock An introduction in discrete time.

\bibitem[Karatzas and Shreve(1991)]{Kara_Shr_BMSC}
I.~Karatzas and S.~E. Shreve.
\newblock \emph{Brownian motion and stochastic calculus}, volume 113 of
  \emph{Graduate Texts in Mathematics}.
\newblock Springer-Verlag, New York, second edition, 1991.
\newblock ISBN 0-387-97655-8.

\bibitem[Karatzas and Shreve(1998)]{Kara_Shr_MF}
I.~Karatzas and S.~E. Shreve.
\newblock \emph{Methods of mathematical finance}, volume~39 of
  \emph{Applications of Mathematics (New York)}.
\newblock Springer-Verlag, New York, 1998.
\newblock ISBN 0-387-94839-2.

\bibitem[Karatzas and Zamfirescu(2006)]{Kara_Zam_2006}
I.~Karatzas and I.~M. Zamfirescu.
\newblock Martingale approach to stochastic control with discretionary
  stopping.
\newblock \emph{Appl. Math. Optim.}, 53\penalty0 (2):\penalty0 163--184, 2006.
\newblock ISSN 0095-4616.

\bibitem[Karatzas and Zamfirescu(2008)]{Kara_Zam_2008}
I.~Karatzas and I.~M. Zamfirescu.
\newblock Martingale approach to stochastic differential games of control and
  stopping.
\newblock \emph{Ann. Probab.}, 36\penalty0 (4):\penalty0 1495--1527, 2008.
\newblock ISSN 0091-1798.

\bibitem[Kazamaki(1994)]{expM_BMO}
N.~Kazamaki.
\newblock \emph{Continuous exponential martingales and {BMO}}, volume 1579 of
  \emph{Lecture Notes in Mathematics}.
\newblock Springer-Verlag, Berlin, 1994.
\newblock ISBN 3-540-58042-5.

\bibitem[Kl{\"o}ppel and Schweizer(2007)]{Kloppel_Schweizer_2007}
S.~Kl{\"o}ppel and M.~Schweizer.
\newblock Dynamic indifference valuation via convex risk measures.
\newblock \emph{Math. Finance}, 17\penalty0 (4):\penalty0 599--627, 2007.
\newblock ISSN 0960-1627.

\bibitem[Kobylanski et~al.(2002)Kobylanski, Lepeltier, Quenez, and
  Torres]{KLQT_RBSDE}
M.~Kobylanski, J.~P. Lepeltier, M.~C. Quenez, and S.~Torres.
\newblock Reflected {BSDE} with superlinear quadratic coefficient.
\newblock \emph{Probab. Math. Statist.}, 22\penalty0 (1, Acta Univ. Wratislav.
  No. 2409):\penalty0 51--83, 2002.
\newblock ISSN 0208-4147.

\bibitem[Lepeltier(1985)]{Lepeltier_1985}
J.-P. Lepeltier.
\newblock On a general zero-sum stochastic game with stopping strategy for one
  player and continuous strategy for the other.
\newblock \emph{Probab. Math. Statist.}, 6\penalty0 (1):\penalty0 43--50, 1985.
\newblock ISSN 0208-4147.

\bibitem[Neveu(1975)]{Neveu_1975}
J.~Neveu.
\newblock \emph{Discrete-parameter martingales}.
\newblock North-Holland Publishing Co., Amsterdam, revised edition, 1975.
\newblock Translated from the French by T. P. Speed, North-Holland Mathematical
  Library, Vol. 10.

\bibitem[Rockafellar(1997)]{rock}
R.~T. Rockafellar.
\newblock \emph{Convex analysis}.
\newblock Princeton Landmarks in Mathematics. Princeton University Press,
  Princeton, NJ, 1997.
\newblock ISBN 0-691-01586-4.
\newblock Reprint of the 1970 original, Princeton Paperbacks.

\bibitem[Williams(1991)]{DW_PwM_91}
D.~Williams.
\newblock \emph{Probability with martingales}.
\newblock Cambridge Mathematical Textbooks. Cambridge University Press,
  Cambridge, 1991.
\newblock ISBN 0-521-40455-X; 0-521-40605-6.

\end{thebibliography}

\end{document}